\documentclass[11pt,twoside]{article}
\usepackage{times}
\usepackage{amsmath,amssymb}
\usepackage{color}
\usepackage{graphicx}

\allowdisplaybreaks
 \def \no{\nonumber}

\newcommand{\be}{\begin{equation}}
\newcommand{\ee}{\end{equation}}
\newcommand{\bea}{\begin{eqnarray}}
\newcommand{\eea}{\end{eqnarray}}

\def\p{\partial}
\def\ve{\varepsilon}
\def\f{\frac}
\def\na{\nabla}
\def\la{\lambda}
\def\al{\alpha}
\def\t{\tilde}
\def\vp{\varphi}
\def\O{\Omega}
\def\o{\omega}
\def\th{\theta}

\def\g{\gamma}

\def\si{\sigma}

\def\dl{\delta}

\def\ds{\displaystyle}
\def\q{\quad}

\def\no{\nonumber}
\def\dP{\dot\Phi}

\def\beq{\begin{equation}}
\def\eeq{\end{equation}}
\def\ben{\begin{eqnarray}}
\def\een{\end{eqnarray}}
\def\bec{\begin{cases}}
\def\eec{\end{cases}}
\def\ss{\sum\limits}
\def\hP{\hat\Phi}
\def\tf{\tilde{f}}

\begin{document}
 \footskip=0pt
 \footnotesep=2pt
\let\oldsection\section
\renewcommand\section{\setcounter{equation}{0}\oldsection}
\renewcommand\thesection{\arabic{section}}
\renewcommand\theequation{\thesection.\arabic{equation}}
\newtheorem{claim}{\noindent Claim}[section]
\newtheorem{theorem}{\noindent Theorem}[section]
\newtheorem{lemma}{\noindent Lemma}[section]
\newtheorem{proposition}{\noindent Proposition}[section]
\newtheorem{definition}{\noindent Definition}[section]
\newtheorem{remark}{\noindent Remark}[section]
\newtheorem{corollary}{\noindent Corollary}[section]
\newtheorem{example}{\noindent Example}[section]

\title{\bf {The global existence and large time behavior of smooth compressible fluid  in an infinitely  expanding ball, I:
3D Euler equations}}
\author{Gang Xu$^{1}$, Huicheng Yin$^{2}$\footnote{Gang Xu (gxu@ujs.edu.cn) and  Huicheng Yin (huicheng@nju.edu.cn)
are supported by the National Natural Science Foundation of China (No.11025105, No.11571177)
and a project funded by the Priority Academic Program Development of Jiangsu Higher Education Institutions.
}\vspace{0.5cm}\\
\small 1.  Faculty of Science, Jiangsu University, Zhenjiang, Jiangsu
212013, China.\\
\small 2. School of Mathematical Sciences, Jiangsu Provincial Key Laboratory for Numerical Simulation\\
\small of Large Scale Complex Systems, Nanjing Normal University, Nanjing 210023, China.\\
}

\date{}
\maketitle

\centerline {\bf Abstract} \vskip 0.3 true cm

We concern with the global existence and large time behavior
of compressible fluids (including the inviscid gases, viscid gases, and Boltzmann gases)
in an infinitely expanding ball. Such a problem is one of
the interesting  models in studying the theory of global smooth solutions to multidimensional
compressible gases with time dependent boundaries and vacuum states at infinite time.
Due to  the  conservation of mass, the fluid in the  expanding ball becomes rarefied and eventually
tends to a vacuum state meanwhile  there
are no appearances of vacuum domains in any part of the expansive ball,
which is easily observed in finite time. In this paper,
as the first part of our three papers,
we will confirm this physical phenomenon for the compressible inviscid fluids
by  obtaining the exact lower and upper bound on the density function.

\vskip 0.3 true cm
{\bf Keywords:} Compressible Euler equations, expanding ball, global existence, degenerate,
weighted energy estimates, large time behavior\vskip 0.3 true cm


\section {Introduction}

In this paper, we consider the behavior of a compressible inviscid  fluid in a 3D expanding ball
given by  $\O_0=\{(t,x): t\ge 0, |x|=\sqrt{x_1^2+x_2^2+x_3^2}\le R_0(t)
\}$, where $R_0(t)\in C^{6}[0, \infty)$ satisfies $R_0(0)=1, R_0'(0)=0, R_0''(0)=0$, and $R_0(t)=1+L t$ for $t\ge 1$
with some positive constant $L$. From the expression of $\O_0$, we know that the  ball
$S^0_t=\{x: |x|\le R_0(t)\}$ at the time $t$ is artificially set by pulling out the initial unit ball $S^0=\{x:
|x|\le 1\}$ with
a smooth speed and acceleration (see Figure 1 below).
\begin{figure}[htbp]
\centering\includegraphics[width=7.5cm,height=6.5cm]{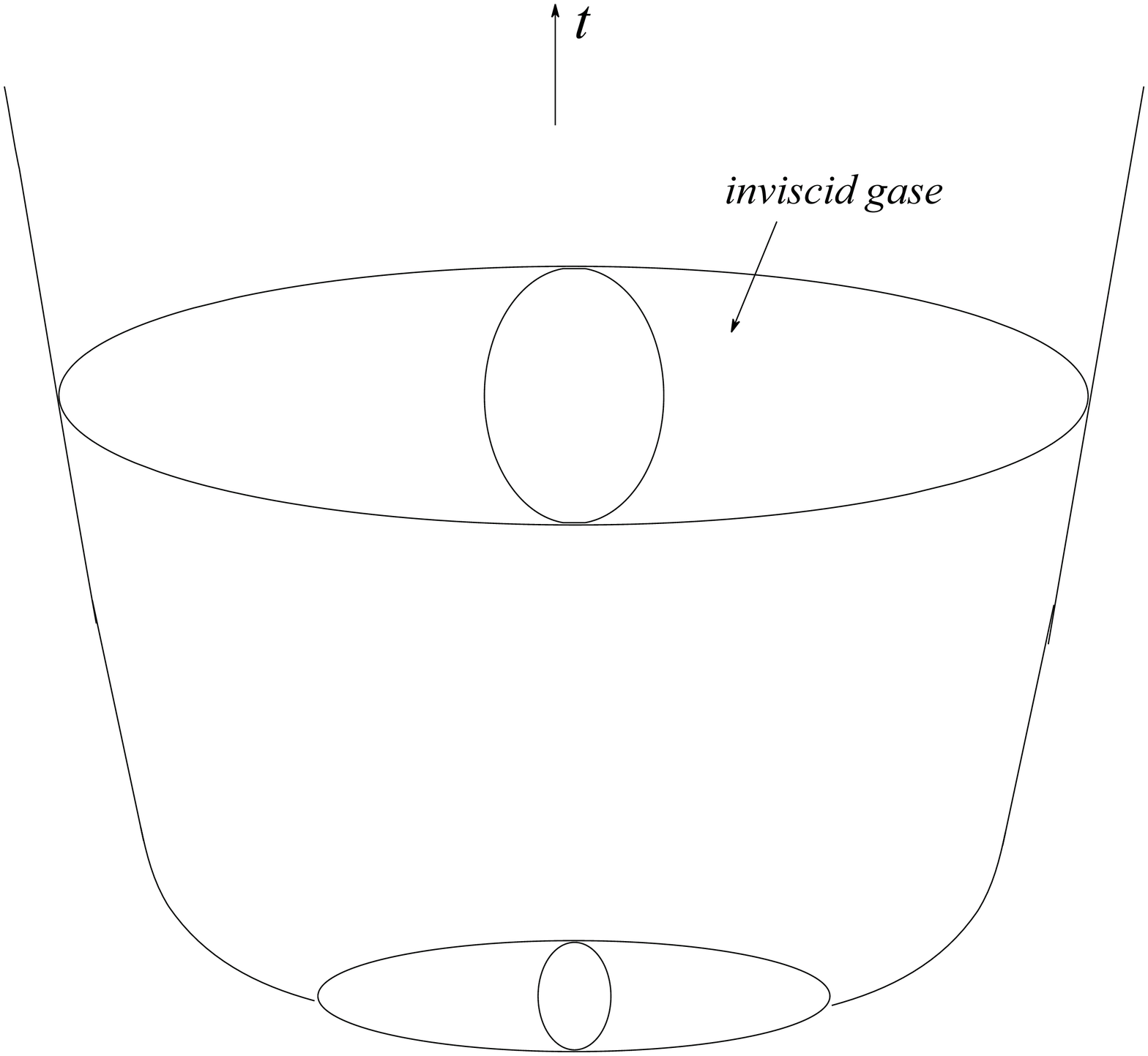}
\caption{\bf An inviscid flow in a 3D expanding ball}\label{fig:1}
\end{figure}

Suppose that the movement of fluid in the ball is governed by the
3D compressible isentropic Euler system:
\beq
\bec
\p_t\rho+\ss_{i=1}^3\p_i(\rho u_i)=0,\\
\p_t(\rho u_i)+\ss_{j=1}^3\p_j(\rho u_iu_j)+\p_iP=0,~~i=1,2,3,
\eec\label{Euler}\eeq
where $x=(x_1, x_2, x_3)$,  $ \rho, u=(u_1, u_2, u_3)$ and $P$ represent the density, velocity and pressure, respectively. Moreover,
assume the state equation $P=A\rho^{\g}$ holds with $A>0$ and $\gamma$ $(1<\g<\f{4}{3})$  being constants.
Without loss of generality, take $A=1$
and impose the following initial-boundary conditions
on (1.1)
\begin{equation}
\left\{
\begin{aligned}
&\rho(0,x)=\rho_0(x), \ \ u(0,x)=u_0(x),\quad\text{for $x\in S^0$},\\
&R'(t)=\ds\sum_{i=1}^3\ds\f{x_iu_i}{|x|},\qquad\qquad \qquad \quad \quad\text{for  $(t,x)\in \p \Omega_0=\{(t,x): t\ge 0, |x|=R_0(t)\}$},\\
\end{aligned}
\right.\label{1.2}
\end{equation}
where $\rho_0(x)\in H^4(S^0)$, $u_0(x)\in H_0^4(S^0)$, and $\rho_0(x)>0$ for $x\in S^0$.
Here, the boundary condition on $\p\O_0$ in (1.2) represents the solid wall.

The purpose of this paper is to prove the following theorem.

\begin{theorem}\label{bt1.1}
If $\rho_0(x)\in H^4(S^0)$, $u_0(x)\in H_0^4(S^0)$,  $rot\, u_0(x)\equiv 0$, and
the compatibility conditions on $\{(t,x): t=0, x\in \p S^0\}$ of $(\rho_0(x), u_0(x))$ hold,
then there exist a constant $h_0>0$, and a
small constant $\ve_0>0$ depending only on $h_0$, such that when
$\ds\sup_{0\le t\le 1, 1\le k\le 5}|R_0^{(k)}(t)|+\|\rho_0(x)-1\|_{H^4(S^0)}+\|u_0(x)\|_{H^4(S^0)}<\ve_0$,
$R_0(t)=1+Lt$ for $t\ge 1$, and $0<L<h_0$,
problem (1.1)-(1.2) with $1<\g<\f43$ admits a global solution
$(\rho, u)$ in $\O_0$ satisfying
\begin{align}
&(\rho(t,x),u(t,x))\in C([0,\infty), H^4(S^0_t)) \cap C^1([0,\infty), H^3(S^0_t)),\label{b1.4}\\
&rot \, u(t,x)\equiv 0,\qquad \qquad \qquad (t,x)\in\O_0,\label{b1.5}\\
&\ds\f{1}{2R^3(t)}\le \rho(t,x)\le\ds\f{3}{2R^3(t)}\qquad\text{for \quad $t\ge 1$},\label{b1.6}
\end{align}
and
\be
\sup \limits_{x\in S_t^0}\biggl(\bigl|u(t,x)-\ds\f{Lx}{R(t)}\bigr|+\bigl|R(t)\na(u(t,x)
-\ds\f{Lx}{R(t)})\bigr|\biggr)\rightarrow 0  \ as \ t \rightarrow +\infty.
\label{b1.7}
\ee
Here, $S^0_t=\{x: |x|\le R_0(t)\}$, and $R(t)=1+Lt$ for $t\ge 0$.
\end{theorem}

In the following, we give some remarks on the above theorem.

\vskip .1cm

{\bf Remark 1.1.} {\it The pointwise estimate $\rho(t,x)\sim \ds\f{1}{R^{3}(t)}$
for large $t$ can be expected by the conservation of mass.}

\vskip .1cm

{\bf Remark 1.2.} {\it  When the fluid is governed
 by the compressible Navier-Stokes equations and
the corresponding boundary condition is given by  $u(t,x)=\ds\f{R'(t)x}{R(t)}$
for $(t,x)\in \p \Omega_0$, the corresponding result was obtained
in [25]. The case for a rarefied gas in $\O_0$ governed by the Boltzmann equation
was also obtained in [26].}

\vskip .1cm

{\bf Remark 1.3.} {\it If the initial density contains vacuum, the local well-posedness of the
compressible Euler system have been extensively studied, cf.  [3], [7-8], [11-12], [16], [18]
and the references therein. In general, classical solutions will blow up in finite time as shown in [3] and [23]. For the problem considered
in this paper, the vacuum is the time asymptotic state
so that the classical solution exists globally in time.}

\vskip .1cm

{\bf Remark 1.4.} {\it When the initial velocity $u_0(x)$  is close to a linear field,  the authors in [10] and [19] proved the global existence
of smooth solution to the Cauchy problem of the compressible Euler system.
And this is different from the case considered in this paper.}

\vskip .1cm

{\bf Remark 1.5.} {\it If the  ball $S^0$  is pulled outwards rapidly, namely, when the number $L$ is large, part of the region inside the ball may
become vacuum  in finite time (see [6] and so on).}

\vskip .1cm

{\bf Remark 1.6.} {\it Note that $(\hat\rho(t,x), \hat u(t,x))=(\ds\f{1}{R^3(t)}, \ds\f{Lx}{R(t)})$ with $R(t)=1+Lt$
is a special solution to (1.1)-(1.2). In fact, Theorem 1.2 gives
 the stability of this special solution.
In addition,  the smallness of $L$ in Theorem 1.1 is only used to prove the local existence of the solution
to (1.1)-(1.2) and to obtain the smallness of $(\rho(1,x)-1, u(1,x))$
 in (6.11).}

\vskip .1cm

{\bf Remark 1.7.} {\it For the spherically symmetric solution,
one can relax
the restriction on $1<\g<\f43$ in Theorem 1.1 to $1<\g<\f53$. The main reasons
come from the simplified boundary conditions
 (5.66) and (5.112) in Section 5, and the assumption
(5.2) in the domain $\{(t,x): t>0, r>\f13R(t)\}$ can be replaced by
\begin{equation*}
\left\{
\begin{aligned}
&|Z\dP|\le M\ve R(t)^{-\si}\quad \text{if $\gamma\in(1,\f43)\cup(\f43,\f53)$};\qquad |Z\dP|\le M\ve \ln R(t)\quad
\text{if $\g=\f43$};\\
&|ZD_t\dot\Phi|\le M\ve R(t)^{-3(\g-1)},\\
\end{aligned}
\right.
\end{equation*}
where $\si=min\{0,3(\g-1)-1\}$. We omit the details of the
estimation for this special case.}
\\

To prove Theorem 1.1,  we first solve an unsteady potential flow equation in
the domain $\O=\{(t,x): t\ge 0, |x|\le R(t)
\}$ with the initial-boundary
conditions (1.2) (see the Figure 2 below).
Let $\Phi(t,x)$ be the potential of velocity $u=(u_1,u_2,u_3)$, i.e., $u_i=\p_i\Phi$ $(1\le i\le 3)$.
Then it follows from the Bernoulli law that
\beq
\p_t\Phi+\f{1}{2}|\nabla_x\Phi|^2+h(\rho)=B_0,\label{Bernoulli}
\eeq
where $h(\rho)=\ds\f{c^2(\rho)}{\g-1}$ is the specific enthalpy, $c(\rho)=\sqrt{P'(\rho)}$ is the local sound speed, $\nabla_x=(\p_1,\p_2,\p_3)$, $B_0=\ds\f{c^2(\rho_0)}{\g-1}$ is the Bernoulli constant of a static state with the constant density $\rho_0$.

\begin{figure}[htbp]
\centering\includegraphics[width=7.0cm,height=6.5cm]{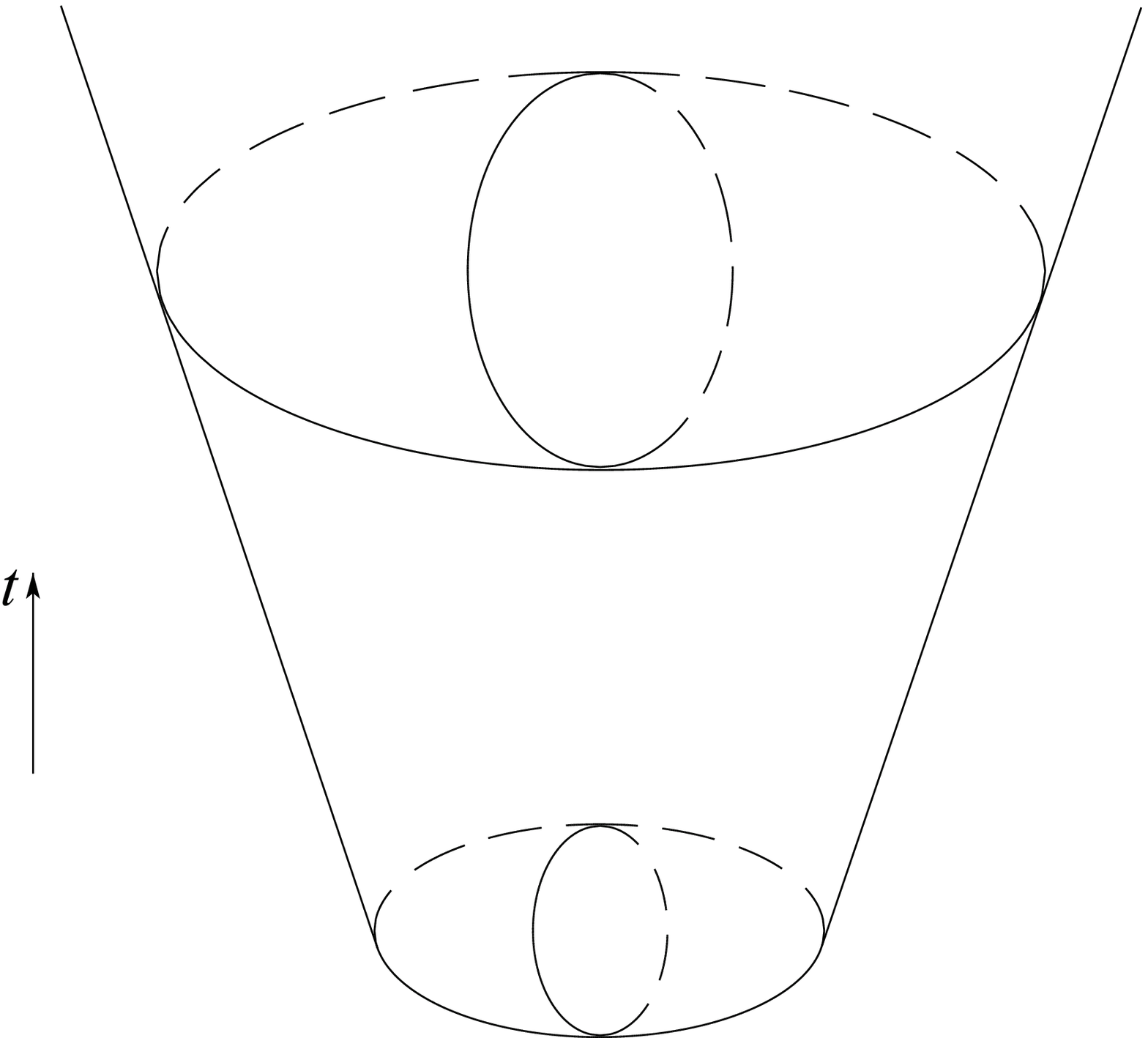}
\caption{}
\end{figure}

By (1.7) and the implicit function theorem with $h'(\rho)=\ds\f{c^2(\rho)}{\rho}>0$ for $\rho>0$, the density function $\rho(t,x)$ can be expressed as
\beq
\rho=h^{-1}\big(B_0-\p_t\Phi-\f{1}{2}|\nabla_x\Phi|^2\big)\equiv H(\nabla\Phi),\label{implicit}
\eeq
where $h^{-1}$ stands for the inverse function of $h(\rho)$, and $\nabla=(\p_t,\nabla_x)$.

Substituting (\ref{implicit}) into the  first equation in (\ref{Euler}) yields
\beq
\p_t(H(\nabla\Phi))+\ss_{i=1}^3\p_i(H(\nabla\Phi)\p_i\Phi)=0.\label{mass}
\eeq
In fact, for any $C^2$ solution $\Phi$, (\ref{mass}) can be rewritten
into the following second order quasilinear equation
\beq
\p_t^2\Phi+2\ss_{k=1}^3\p_k\Phi\p_{tk}^2\Phi+\ss_{i,j=1}^3\p_i\Phi\p_j\Phi\p_{ij}^2\Phi-c^2(\rho)\Delta \Phi=0.\label{potential}
\eeq
Denote the lateral boundary of $\O$ by $\p\O=\{(t,x): t\ge 0, |x|=R(t)\}$. Then  on $\p\O$,
\beq
\ss_{i=1}^3\p_i\Phi\cdot\f{x_i}{|x|}=L.\label{bd condition}
\eeq
Due to the geometric property of $\O$, it is convenient to work in the spherical coordinates $(r,\th,\vp)$:
\beq
(x_1,x_2,x_3)=(r\cos\th\sin\vp,r\sin\th\sin\vp,r\cos\vp),\label{spherical}
\eeq
where $r=\sqrt{x_1^2+x_2^2+x_3^2}$, $0\le\th\le 2\pi$, $0\le\vp\le\ds\f{\pi}{2}$.
Under the coordinate transformation (\ref{spherical}), (\ref{potential}) becomes
\ben
&&\p_t^2\Phi+2\p_r\Phi\p_{tr}^2\Phi+\f{2}{r^2}\p_{\vp}\Phi\p_{t\vp}^2\Phi+\f{2}{r^2\sin^2\vp}\p_{\th}\Phi\p_{t\th}^2\Phi
+\big((\p_r\Phi)^2-c^2(\rho)\big)\p_r^2\Phi\no\\
&&\quad +
\f{1}{r^2\sin^2\vp}\bigg(\f{(\p_{\th}\Phi)^2}{r^2\sin^2\vp}-c^2(\rho)\bigg)\p_{\th}^2\Phi
+\f{1}{r^2}\bigg(\f{(\p_{\vp}\Phi)^2}{r^2}-c^2(\rho)\bigg)\p_{\vp}^2\Phi+\f{2\p_r\Phi\p_{\th}\Phi}{r^2\sin^2\vp}\p_{r\th}^2\Phi
\no\\
&&\quad +\f{2}{r^2}\p_r\Phi\p_{\vp}\Phi\p_{r\vp}^2\Phi
+\f{2\p_{\th}\Phi\p_{\vp}\Phi}{r^4\sin^2\vp}\p_{\th\vp}^2\Phi
-\f{1}{r^3}\bigg(2r^2c^2(\rho)+(\p_{\vp}\Phi)^2+\f{(\p_{\th}\Phi)^2}{\sin^2{\vp}}\bigg)\p_r\Phi\no\\
&&\quad -\f{\cot\vp}{r^4}\bigg(r^2c^2(\rho)+\f{(\p_{\th}\Phi)^2}{\sin^2\vp}\bigg)\p_{\vp}\Phi=0.\label{potential-sphere}
\een
Note that some coefficients in (\ref{potential-sphere}) have strong singularities near $\vp=0$.
Consequently,  as in [14] and [24], we   rewrite (\ref{potential-sphere}) by introducing some
smooth vector fields tangent to the sphere $\Bbb S^2$.

Set
\beq
\bec
&Z_1=x_1\p_2-x_2\p_1=\p_{\th},\\
&Z_2=x_2\p_3-x_3\p_2=-\cot\vp\cos\th\p_{\th}-\sin\th\p_{\vp},\\
&Z_3=x_3\p_1-x_1\p_3=-\cot\vp\sin\th\p_{\th}+\cos\th\p_{\vp}.
\eec
\eeq
Then it follows from a direct computation that (\ref{potential-sphere})
can be written as
\ben
&&\p_t^2\Phi+2\p_r\Phi\p_{tr}^2\Phi+\f{2}{r^2}\ss_{i=1}^3Z_i\Phi\p_tZ_i\Phi+((\p_r\Phi)^2-c^2(\rho))\p_r^2\Phi
\no\\
&&\quad +\f{2\p_r\Phi}{r^2}\sum\limits_{i=1}^3Z_i\Phi\p_rZ_i\Phi
-\f{c^2(\rho)}{r^2}\sum\limits_{i=1}^3Z_i^2\Phi
+\f{1}{r^4}\sum\limits_{i,j=1}^3Z_i\Phi
Z_j\Phi
Z_iZ_j\Phi\no\\
&&\quad +\sum\limits_{i,j=1}^3\f{C_{ij}(\o)}{r^3}\p_r\Phi Z_i\Phi
Z_j\Phi+\sum\limits_{i,j,k=1}^3\f{C_{ijk}(\o)}{r^4}Z_i\Phi
Z_j\Phi
Z_k\Phi-\f{2c^2(\rho)}{r}\p_r\Phi=0,\label{potential-z}
\een
where $\o=\ds \f{x}{r}$,
$C_{ij}(\o)=C_{ij}(\ds \f{x}{r})$ and $C_{ijk}(\o)=C_{ijk}(\ds\f{x}{r})$
are smooth functions of their arguments.

Meanwhile, the boundary condition (\ref{bd condition}) becomes
\beq
\p_r\Phi=L,\qquad\qquad \text{on $\p\O$}.\label{bd-condition-z}
\eeq

In addition, we impose the following initial perturbation:
\beq
\Phi(0,x)=\f{1}{2}L|x|^2+\ve\Phi_0(x),\q\p_t\Phi(0,x)=-\f{1}{2}L^2|x|^2+\ve\Phi_1(x),\label{initial condition}
\eeq
where $\ve>0$ is a small parameter, $(\Phi_0(x), \Phi_1(x))\in (H^5(S^0), H^4(S^0))$, and the initial-boundary
value conditions (1.16)-(1.17) are compatible on $S^0$.
Note that the initial data (1.17) can be replaced by $(\Phi(0, x), \p_t\Phi(0,x))=(\ve\Phi_0(x),
\ve\Phi_1(x))$ when $L>0$ is small. On the other hand, due to $u_i=\p_i\Phi$ and (\ref{implicit}), the initial
conditions (\ref{initial condition}) can be realized by a small perturbation of the initial density and velocity of an
irrotational flow.

\begin{theorem}
Under the above assumptions on the initial and boundary
data, if $\g\in(1,\f43)$, then there exists a constant $\ve_0>0$ depending on $L, B_0$ and $\g$ such that problem (\ref{potential}) with (\ref{bd-condition-z})-(\ref{initial condition}) has a global solution $\Phi(t,x)\in C([0,\infty), H^5(S_t))\cap
{C^1}([0,\infty), H^4(S_t))$ for $\ve<\ve_0$, where $S_t=\{x: r\le R(t)\}$. Moreover, $\rho(t,x)>0$ and $\lim\limits_{t\rightarrow\infty}\rho(t,x)=0$ hold.
\end{theorem}
\vskip .1cm
{\bf Remark 1.8.} {\it The linearized operator
of the quasilinear wave equation (1.10) around the
special expanding solution has the approximate  form of
$$\p_t^2-\ds\f{\g}{(1+Lt)^{3(\g-1)}}
(\p_1^2+\p_2^2+\p_3^2)+\ds\f{3L(\g-1)}{1+Lt}\p_t.
$$
On the other hand,
if one considers the Cauchy problem of (1.1) for initial
data as a small perturbation of a uniform constant density $\rho_0$ and velocity $(0,0,q_0)$, that is,
\begin{equation}
\left\{
\begin{aligned}
&\p_t^2\Phi+2\ss_{k=1}^3\p_k\Phi\p_{tk}^2\Phi+\ss_{i,j=1}^3\p_i\Phi\p_j\Phi\p_{ij}^2\Phi-c^2(\rho)\Delta \Phi=0,\\
&\Phi(t,x)|_{t=0}=\ve\Phi_0(x),\quad \p_t\Phi(t,x)|_{t=0}=q_0+\ve\Phi_1(x),\qquad x
\in\Bbb R^3,
\end{aligned}
\right.\tag{1.18}
\end{equation}
where $\Phi_i(x)\in C_0^{\infty}(\Bbb R^3)$ ($i=0,1$),
then  (1.18) does not fulfill the ``null-condition'' introduced
in [4] and [13]. Therefore, according to the
results obtained in [1-2], [5], [20] and [27],  classical solution to
(1.18)  blows up and then  shock forms in finite time. Compared
this blowup result with  Theorem 1.1-1.2,
the global existence of  smooth solution to (1.10) together with a fixed
wall condition comes from the rarefaction property of
fluid.}\\


We now give some remarks on the proof of Theorem 1.2. Since the local solvability
of problem (1.10) together with (1.16)-(1.17) is known
as long as the vacuum does not appear, cf., for example [17], the
proof of Theorem 1.2 is based on the continuation argument.

First of all, note that
the linearized operator
$$\mathcal{L}=\p_t^2+\ds\f{2Lr}{R(t)}\p_{tr}^2
+\ds\sum_{i,j=1}^3\ds\f{L^2x_ix_j}{R^2(t)}\p_{ij}^2-\ds\f{\g}{R^{3(\g-1)}(t)}\Delta+
\ds\f{3L(\g-1)}{R(t)}(\p_t+\ds\f{Lr}{R(t)}\p_r),$$
cf. (3.2),  is  different from the corresponding one
in [24] which is
$$\p_t^2-\ds\f{1}{R^{2(\g-1)}(t)}(\p_1^2+\p_2^2)+\ds\f{2L(\g-1)}{R(t)}\p_t.$$
The key ingredients in
the analysis for the global existence and pointwise
estimate are to obtain some weighted energy estimates
by choosing some appropriate multiplier and  anisotropic weights.
For this, we need to

1.   show that the solution does not contain vacuum both on
the boundary and inside the region;

2.  obtain the different time decay rates of the density and velocity of
 the solution  when it tends
to vacuum state at infinite time;

3. fully use the Neumann-type boundary condition (1.16) on $\Phi$
by applying   the material derivative
$D_t=\p_t+\ds\f{Lr}{R(t)}\p_r$ on the solution, and estimate
the radial derivatives and angular derivatives of $\Phi$,
together with some weighted Sobolev interpolation inequalities given
in [15].

\vskip 0.1 true cm

The rest of the paper is organized as follows. In the next section, we
will  give some basic
properties of the background solution
and  some preliminary  weighted Sobolev interpolation
inequalities. In Section 3, we will reformulate problem (1.10) together with (1.16)-(1.17) by
decomposing its solution as a sum of the background solution and a
small perturbation so that its linearization can be studied
clearly. In Section 4, we will establish a uniform weighted
energy estimate for the corresponding linear problem, where an
appropriate multiplier is constructed. In Section 5, the uniform higher
order weighted estimates of
$\dot\Phi$ are obtained by the  analysis on the radial derivatives and angular derivatives of
$\dot\Phi$, where the domain decomposition
technique is  applied. In the last section, we complete the proof of Theorem~1.2 by applying the Sobolev embedding theorem and the continuation argument,
and then Theorem 1.1 follows from
 Theorem 1.2 and the local existence result on the problem (1.1)-(1.2).

\section{Background solution and some preliminaries}

In this section,  we analyze the background solution to (\ref{potential}) with (\ref{bd-condition-z})-(\ref{initial condition}) when the initial data (\ref{initial condition}) are
\beq
\hat{\Phi}(0,x)=\f{1}{2}L|x|^2,\q\p_t\hat{\Phi}(0,x)=-\f{1}{2}L^2|x|^2.\label{initial background}
\eeq
In this case, the density $\rho(x)$ and velocity $u(t,x)=\nabla_x\Phi(t,x)$ in $\O$ take the form: $\rho(t,x)=\hat{\rho}(t,r)$, $u(t,x)=\f{\ds x}{\ds r}\hat{U}(t,r)$. Consequently, problem (\ref{potential}) with (\ref{bd-condition-z}) and (\ref{initial background}) is equivalent to
\beq
\bec
&r^2\p_t\hat{\rho}+\p_r(r^2\hat{\rho}\hat{U})=0,\\
&\p_t(r^2\hat{\rho}\hat{U})+\p_r(r^2\hat{\rho}\hat{U})+r^2\p_rP=0,\\
&\hat{\rho}(0,r)=1,\q \hat{U}(0,r)=Lr.
\eec\label{Euler-background}
\eeq
One can easily check that (\ref{Euler-background}) has a solution
\beq
\hat{\rho}(t,r)=\ds\f{1}{R^3(t)},\q \hat{U}(t,r)=\f{\ds Lr}{\ds R(t)}.\label{background}
\eeq
Then for $1<\g<\ds\f{4}{3}$, it follows from $u_i=\p_i\Phi$, (\ref{Bernoulli}) and (\ref{background}) that (\ref{potential}) with  (\ref{bd-condition-z}) and (\ref{initial background}) has a  solution
\beq
\hat{\Phi}(t,r)=\f{\ds \g}{\ds (\g-1)(4-3\g)L}+B_0t+\f{\ds Lr^2}{\ds 2(1+Lt)}-\f{\ds \g}{\ds (\g-1)(4-3\g)L}(1+Lt)^{4-3\g},\label{background-Phi}
\eeq
where $B_0=\ds\f{\g}{\g-1}$.

Next, we include the weighted Sobolev interpolation inequality from
 [15] that will be used  in Lemma 2.4.

\begin{lemma}\label{interpolation}Suppose $s, \tau, p, \al, \beta, q, a$ are real
numbers, and $j\ge0, m>0$ are integers, satisfying
\beq
\bec
&\ds p,q\ge1, \f{j}{m}\le a\le1, s>0,\\
&\ds \f{1}{s}+\f{\tau}{n}>0,
\f{1}{p}+\f{\al}{n}>0, \f{1}{q}+\f{\beta}{n}>0,\\
&\ds m-j-\f{n}{p} \quad \text {is not a non-negative integer.}
\eec\label{lemma2.1-1}\eeq
There exists a positive constant $C$ such that the following
inequality holds for all $v\in C_0^{\infty}(\Bbb R^n)$:
\beq
\big||x|^{\tau}\na_x^jv\big|_{L^s}\le
C\big||x|^{\al}\na_x^mv\big|_{L^p}^a\big||x|^{\beta}v\big|_{L^q}^{1-a},\label{lemma2.1}
\eeq
if and only if the following conditions hold:
\beq
\f{1}{s}+\f{\tau-j}{n}=a(\f{1}{p}+\f{\al-m}{n})+(1-a)(\f{1}{q}+\f{\beta}{n})\quad\text{with} \quad\tau\le a\al+(1-a)\beta;\label{lemma2.1-2}
\eeq
if $\ds\f{1}{q}+\f{\beta}{n}=\f{1}{p}+\f{\al-m}{n}$, then
\beq
a(\al-m)+(1-a)\beta+j\le\tau; \label{lemma2.1-3}
\eeq
if $\ds a=\f{j}{m}$, then
\beq\tau=a\al+(1-a)\beta.\label{lemma2.1-4}\eeq
\end{lemma}

\begin{corollary}For the domain $Q$ (\text{$=\O_0$ or $\O$}) defined in
Section 1, if a function $u(y)\in C^m(\bar Q)$ and
\beq
u|_{|y|\ge T}\equiv0,\label{restriction}
\eeq
where $T>0$ is a constant, $y=(t, x)$ and $|y|=\sqrt{t^2+|x|^2}$, then

(i) (\ref{lemma2.1}) still holds under the restriction (\ref{lemma2.1-1})
and (\ref{lemma2.1-2})-(\ref{lemma2.1-4}), moreover, the constant $C$ on the right hand side of (\ref{lemma2.1})
is independent of $T$.

(ii) When $m=2, 1<\g<\f53$ and $0<\dl<3\g$,  for $\mu>0$,
\ben
&&\big||y|^{\f{2\mu+2-\dl}4}\nabla u\big|_{L^4(Q)}\le C\big||y|^{\f{\mu+1-\dl}2}\nabla^2u\big|_{L^2(Q)}^{\f12}\big||y|^{\f{\mu+1}2}u\big|_{L^{\infty}(Q)}^{\f12},\\
&&\big||y|^{\f{\mu+3-\dl}2}\nabla u\big|_{L^4(Q)}\le C\big||y|^{\f{\mu+3-\dl}2}\nabla^2u\big|_{L^2(Q)}^{\f12}\big||y|^{\f{\mu+3-\dl}2}u\big|_{L^{\infty}(Q)}^{\f12},\\
&&\big||y|^{\f{2\mu+2-3(\g-1)-\dl}4}\nabla u\big|_{L^4(Q)}\le C\big||y|^{\f{\mu+1-3(\g-1)}2}\nabla^2u\big|_{L^2(Q)}^{\f12}\big||y|^{\f{\mu+1-\dl}2}u\big|_{L^{\infty}(Q)}^{\f12},\\
&&\big||y|^{\f{3\g-3-\dl}2}\nabla u\big|_{L^4(Q)}\le C\big||y|^{\f{3\g-3-\dl}2}\nabla^2u\big|_{L^2(Q)}^{\f12}\big||y|^{\f{3\g-3-\dl}2}u\big|_{L^{\infty}(Q)}^{\f12},\\
&&\big||y|^{\f{6\g-6-\dl}2}\nabla u\big|_{L^4(Q)}\le C\big||y|^{\f{6\g-6-\dl}2}\nabla^2u\big|_{L^2(Q)}^{\f12}\big||y|^{\f{6\g-6-\dl}2}u\big|_{L^{\infty}(Q)}^{\f12},
\een
where $\na=\na_y$ and $\na^2=\na_y^2$.
\end{corollary}
{\bf Proof.} (i) The proof is completely parallel to that of Lemma \ref{interpolation} (one can check the details in [15]), then
we omit it here.

(ii) In (\ref{lemma2.1-1}) and (\ref{lemma2.1-2})-(\ref{lemma2.1-4}) of Lemma \ref{interpolation}, set $s=4,p=2,q=\infty, a=\f12$ and $j=1,m=2,n=4$, one can conclude that
(2.11) and (2.12) come from (\ref{lemma2.1}) with $\tau=\f{2\mu+2-\dl}4$, $\al=\f{\mu+1-\dl}2$, $\beta=\f{\mu+1}2$ and $\tau=\f{\mu+3-\dl}2$, $\al=\f{\mu+3-\dl}2$, $\beta=\f{\mu+3-\dl}2$, respectively.

(2.13) and (2.14) are from (\ref{lemma2.1}) with $\tau=\f{2\mu+2-3(\g-1)-\dl}4$, $\al=\f{\mu+1-3(\g-1)}2$, $\beta=\f{\mu+1-\dl}2$ and $\tau=\al=\beta=\f{3\g-3-\dl}2$, respectively.

(2.15) is from (\ref{lemma2.1}) with $\tau=\al=\beta=\f{6\g-6-\dl}2$. And this
completes the proof of the corollary. $\hfill\square$\\

As in [24], in order to apply Lemma \ref{interpolation} or Corollary 2.3 to derive some weighted Sobolev inequalities in $Q$ without
the restriction (\ref{restriction}), we need
to establish an extension result as stated in

\begin{lemma} Set $y=(t,x)$, $Q_T=\{y: 0<t<T,|x|\le R(t)\}$ for $T>1$ and $D_{S}=\{y: (t+\f{1}{L})^2+|x|^2\le S^2, t\ge 0, |x|\le R(t)\}$
for $S>1+\f{1}{L}$. If $u(t,x)\in C^{3}(\bar Q_T)$ and $R(t)^{\beta}\na^{\alpha}_yu\in L^2(Q_T)$ $(|\al|\le 3)$ with some $\beta\in\Bbb R$,
then there exists an extension $Eu\in C^3(\bar D_{\f{9}{8}H})$  of $u$ such that $Eu=u$ in $Q_T$,
$Eu|_{t\ge \f{9}{8}H}\equiv 0$ and
\beq
|R(t)^{\beta}Eu|_{L^\infty(D_{\f{9}{8}H})}\le C|R(t)^{\beta}u|_{L^\infty(Q_{T})},\quad |R(t)^{\beta}\na^{\al}_yEu|_{L^2(D_{\f{9}{8}H})}\le C \ds\sum_{|\nu|\le |\al|}|R(t)^{\beta-|\al|+|\nu|}\na^{\nu}_yu|_{L^2(D_T)},\label{extension}
\eeq
where $H=\sqrt{R(T)^2+(T+\f{1}{L})^2}$, and $C>0$ is independent of $T$.
\end{lemma}

{\bf Proof.}
Let $\t E$ be an extension operator defined by
$$
(\t Eu)(t,x)=\bec u(t,x),\qquad \quad \q\q\q\q\q\q\q\q\q 0\le t\le T,\quad |x|\le R(t);\\
\ss_{j=1}^4\la_j u(T+\ds\f{j(T-t)}{4L^2+4L+1},x),\q\q T<t\le H-\f{1}{L},\quad |x|\le \sqrt{H^2-(t+\f{1}{L})^2},\\
\eec
$$
where $\ss_{j=1}^4(-\ds\f{j}{4L^2+4L+1})^k\la_j=1$ for $k=0, 1, 2, 3$.

Note that $H-\f{1}{L}<T+LR(T)$ and $T>1$ hold. This means that for $T<t\le H-\f{1}{L}$ and $0\le j\le3$,
$$T\ge T+\f{j(T-t)}{4L^2+4L+1}\ge T-\f{4LR(T)}{4L^2+4L+1}> \f{T}{4L^2+4L+1}.$$
In addition,
$$1\le \f{t}{T+\f{j(T-t)}{4L^2+4L+1}}\le (4L^2+4L+1)(L^2+L+1)\q \text{for}~~ T<t\le H-\f{1}{L}~~ \text{and}~~ 0\le j\le3.$$
Consequently, one has that for $(t,x)\in D_{S}$,
\beq
C_1R(t)\le \sqrt{(t+\f{1}{L})^2+|x|^2}\le C_2R(t).
\eeq
This yields for $|\al|\le3$
\beq\label{extension1}
|R(t)^{\beta}\tilde{E}u|_{L^{\infty}(D_H)}\le C|R(t)^{\beta}u|_{L^{\infty}(Q_T)},
\q|R(t)^{\beta}\nabla_y^{\alpha}\tilde{E}u|_{L^2(D_H)}\le C|R(t)^{\beta}\nabla_y^{\alpha}u|_{L^2(Q_T)}.
\eeq
Next, we construct the extension operator $E$ starting from the operator $\t E$.
In terms of the geometric property of $D_T$, it is convenient to use the spherical coordinates. Let $x_0=t+\f{1}{L}$, and
$$s^2=\ss_{i=0}^3x_i^2,\q \omega_i=\f{x_i}{s},\q i=0,1,2,3.$$
Denoting by $\tilde{u}(s,\omega)=u(s\omega)$ with $\o=(\o_0, ..., \o_3)$.
Let $\bar{E}$ be an extension operator defined by
\beq
(\bar{E}u)(s,\omega)=
\bec \tilde{E}\tilde{u}(s,\omega),\q\q\q\q\q\q\q\q\q (t,x)\in D_H,\\
\ss_{j=1}^4\nu_j\tilde{E}\tilde{u}(H+j(H-s),\omega),\q (t,x)\in D_{\f98H}\backslash D_H,
\eec
\eeq
where $\ss_{j=1}^4(-1)^k\nu_j=1$ for $k=0,1,2,3$.

Note that
$$1\le \f{s}{H+j(H-s)}\le\f94\q\q \text{for}~~H\le s\le \f98H~~ \text{and}~~0\le j\le3.$$
This together with (\ref{extension1})  yields
\beq\label{extension2}
\bec
|R(t)^{\beta}\bar{E}u|_{L^{\infty}(D_{\f98H})}\le C|R(t)^{\beta}\tilde{E}u|_{L^{\infty}(D_H)} \le C|R(t)^{\beta}u|_{L^{\infty}(Q_T)},\\
|R(t)^{\beta}\nabla_{t,x}^{\alpha}\bar{E}u|_{L^2(D_{\f98H})}\le C|R(t)^{\beta}\nabla_{t,x}^{\alpha}\tilde{E}u|_{L^2(D_H)}\le C|R(t)^{\beta}\nabla_{t,x}^{\alpha}u|_{L^2(Q_T)}.
\eec
\eeq
Choose a $C^{\infty}$ smooth function $\eta(s)$ with $\eta(s)=1$ for $s\le 1$ and $\eta(s)=0$ for $s\ge \f98$ and set
$$Eu(t,x)=\eta(\f{s}{H})\bar{E}u.$$
Then $Eu$ satisfies (\ref{extension}), and thus  Lemma 2.2 is proved. $\hfill\square$

\begin{remark}
Lemma 2.2 implies  that Corollary 2.1 still holds without assumption (\ref{restriction}).
\end{remark}
The $Z-$fileds introduced in (1.14), as shown in [14], have the following properties.
\begin{lemma}\label{prop-of-z}
\ben
&(i)& [Z_1, Z_2]=Z_3, [Z_2, Z_3]=Z_1, [Z_3,Z_1]=Z_2.\no\\
&(ii)&[Z_i, \p_r]=0,  Z_ir=0, [Z_i,\Delta]=0.\no\\
&(iii)&\ds\na_x f\cdot\na_x g=\p_r f\cdot\p_r
g+\f{1}{r^2}\ds\sum_{i=1}^{3}Z_i f\cdot Z_i g\quad \text{for any $C^1$ smooth
functions $f$ and $g$}.\no\\
&(iv)&  |Z v|\leq r|\na_x v|\quad\text{for any $C^1$ smooth function $v$},
\text{here and below $Z\in\{Z_1, Z_2, Z_3\}$}.\no\\
&(v)&  \ds\p_1=\f{x_1}{r}\p_r+\f{x_2}{r^2}Z_1-\f{x_3}{r^2}Z_3;
\quad \p_2=\f{x_2}{r}\p_r+\f{x_3}{r^2}Z_2-\f{x_1}{r^2}Z_1;\quad
\p_3=\f{x_3}{r}\p_r+\f{x_1}{r^2}Z_3-\f{x_2}{r^2}Z_2.\no
\een
\end{lemma}

\begin{remark} If $u\in C^m(\Bbb R^3)$ with $m\in\Bbb N$, then by Lemma \ref{prop-of-z} we have
$|\na_x^mu|\sim |\p_r^mu|+\ds\f{|\p_r^{m-1}Zu|}{r}+\ds\f{|\p_r^{m-2}Z^2u|}{r^2}+...+\ds\f{|Z^mu|}{r^m}$.
\end{remark}

As a direct consequence of Remark 2.1, we have the following weighted
inequalities which will be used often in Section 5.
\begin{lemma}
Define $D_t=\p_t+\ds\f{Lr}{R(t)}\p_r$ and $S_0=R(t)D_t$. If $1<\g<\f53, 0<\dl\le \f35(\g-1)$,
and $u(y)\in C^4(\bar{Q}_T)$ (\text{$Q_T=\O_0\cap\{t\le T\}$ or $\O\cap\{t\le T\}$})
with $y=(t,x)$,
then there exists a generic positive constant $C$ independent of $T$ such that
\ben
&(i)&|R(t)^{\f{2\mu+2-\dl}4}\nabla_yS_0^2u|_{L^4(Q_T)}\le C\bigg(\ss_{k=0}^2|R(t)^{\f{\mu+1-\dl}2-k}\nabla_y^{2-k}S_0^2u|_{L^2(Q_T)}^{\f12}\bigg)|R(t)^{\f{\mu+1}2}S_0^2u|_{L^{\infty}(Q_T)}^{\f12};\no\\
&(ii)&|R(t)^{\f{\mu+3-\dl}2}\nabla_y^2S_0u|_{L^4(Q_T)}\le C\bigg(\ss_{k=0}^2|R(t)^{\f{\mu+3-\dl}2-k}\nabla_y^{3-k}S_0u|_{L^2(Q_T)}^{\f12}\bigg)|R(t)^{\f{\mu+3-\dl}2}\nabla_yS_0u|_{L^{\infty}(Q_T)}^{\f12};\no\\
&(iii)&|R(t)^{\f{2\mu+5-3\g-\dl}4}\nabla_yZu|_{L^4(Q_T)}\le C\bigg(\ss_{k=0}^1|R(t)^{\f{\mu+4-3\g}2-k}\nabla_y^{2-k}Zu|_{L^2(Q_T)}^{\f12}\bigg)|R(t)^{\f{\mu+1-\dl}2}Zu|_{L^{\infty}(Q_T)}^{\f12};\no\\
&(iv)&|R(t)^{\f{3\g-3-\dl}2}\nabla_y^2Zu|_{L^4(Q_T)}\no\\
&&\le C\bigg(\ss_{k=0}^2|R(t)^{\f{3\g-3-\dl}2-k}\nabla_y^{3-k}Zu|_{L^2(Q_T)}^{\f12}\bigg)|R(t)^{\f{3\g-3-\dl}2}\nabla_{t,x}Zu|_{L^{\infty}(Q_T)}^{\f12};\no\\
&(v)&|R(t)^{\f{6\g-6-\dl}2}\nabla_yZD_tu|_{L^4(Q_T)}\no\\
&&\le C\bigg(\ss_{k=0}^2|R(t)^{\f{6\g-6-\dl}2-k}\nabla_y^{2-k}ZD_tu|_{L^2(Q_T)}^{\f12}\bigg)|R(t)^{\f{6\g-6-\dl}2}ZD_tu|_{L^{\infty}(Q_T)}^{\f12};\no\\
&(vi)&|R(t)^{\f{3\g-5-\dl}2}\nabla_yZ^2u|_{L^4(Q_T)}+|R(t)^{\f{3\g-5-\dl}2}Z^2\nabla_yu|_{L^4(Q_T)}\no\\
&&\le C\bigg(|R(t)^{\f{3\g-3-\dl}2}\nabla_y^2Zu|_{L^4(Q_T)}+|R(t)^{\f{2\mu+5-3\g-\dl}4}\nabla_yZu|_{L^4(Q_T)}\bigg);\no\\
&(vii)&|R(t)^{\f{6\g-8-\dl}2}Z^2D_tu|_{L^4(Q_T)}\le |R(t)^{\f{6\g-6-\dl}2}\nabla_yZD_tu|_{L^4(Q_T)}.\no
\een
\end{lemma}

{\bf Proof.} Let $E$ be the extension operator given in Lemma 2.2. Then we have

(i)
\ben
&&|R(t)^{\f{2\mu+2-\dl}4}\nabla_yS_0^2u|_{L^4(Q_T)}\no\\
&\le&C|R(t)^{\f{2\mu+2-\dl}4}\nabla_yE(S_0^2u)|_{L^4(Q_T)}\no\\
&\le&C|R(t)^{\f{\mu+1-\dl}2}\nabla_y^{2}E(S_0^2u)|_{L^2(Q_T)}^{\f12}|R(t)^{\f{\mu+1}2}E(S_0^2u)|_{L^{\infty}(Q_T)}^{\f12}
\qquad\quad\text{(Applying (2.11) for $E(S_0^2u)$)}\no\\
&\le&C\bigg(\ss_{k=0}^2|R(t)^{\f{\mu+1-\dl}2-k}\nabla_y^{2-k}S_0^2u|_{L^2(Q_T)}^{\f12}\bigg)|R(t)^{\f{\mu+1}2}S_0^2u|_{L^{\infty}(Q_T)}^{\f12}. \qquad \quad \text{(By Lemma 2.2)}\no
\een

(ii)
\ben
&&|R(t)^{\f{\mu+3-\dl}2}\nabla_y^2S_0u|_{L^4(Q_T)}\no\\
&\le&C|R(t)^{\f{\mu+3-\dl}2}\nabla_yE(\nabla_yS_0u)|_{L^4(Q_T)}\no\\
&\le&|R(t)^{\f{\mu+3-\dl}2}\nabla_y^{2}E(\nabla_yS_0u)|_{L^2(Q_T)}^{\f12}|R(t)^{\f{\mu+3-\dl}2}E(\nabla_yS_0u)|_{L^{\infty}(Q_T)}^{\f12}
~~\text{(Applying (2.12) for $E(\nabla_yS_0u)$)} \no\\
&\le&C\bigg(\ss_{k=0}^2|R(t)^{\f{\mu+3-\dl}2-k}\nabla_y^{3-k}S_0u|_{L^2(Q_T)}^{\f12}\biggr)|R(t)^{\f{\mu+3-\dl}2}
\nabla_yS_0u|_{L^{\infty}(Q_T)}^{\f12}.\qquad\quad\text{(By Lemma 2.2)}\no
\een

(iii)
\ben
&&|R(t)^{\f{2\mu+5-3\g-\dl}4}\nabla_yZu|_{L^4(Q_T)}\no\\
&\le&C|R(t)^{\f{2\mu+5-3\g-\dl}4}\nabla_yE(Zu)|_{L^4(Q_T)}\no\\
&\le&C|R(t)^{\f{\mu+4-3\g}2}\nabla_y^2E(Zu)|_{L^2(Q_T)}^{\f12}|R(t)^{\f{\mu+1-\dl}2}E(Zu)|_{L^{\infty}(Q_T)}^{\f12}\qquad\text{(Applying (2.13) for
$E(Zu)$)} \no\\
&\le&C\bigg(\ss_{k=0}^1|R(t)^{\f{\mu+4-3\g}2-k}\nabla_y^{2-k}Zu|_{L^2(Q_T)}^{\f12}\bigg)|R(t)^{\f{\mu+1-\dl}2}Zu|_{L^{\infty}(Q_T)}^{\f12}.
\qquad\quad\text{(By Lemma 2.2)}\no
\een

(iv)
\ben
&&|R(t)^{\f{3\g-3-\dl}2}\nabla_y^2Zu|_{L^4(Q_T)}\no\\
&\le&|R(t)^{\f{3\g-3-\dl}2}\nabla_yE(\nabla_yZu)|_{L^4(Q_T)}\no\\
&\le&C|R(t)^{\f{3\g-3-\dl}2}\nabla_y^{2}E(\nabla_yZu)|_{L^2(Q_T)}^{\f12}|R(t)^{\f{3\g-3-\dl}2}E(\nabla_yZu)|_{L^{\infty}(Q_T)}^{\f12}
\quad\text{(Applying (2.14) for $E(\nabla_yZu)$)} \no\\
&\le&C\bigg(\ss_{k=0}^2|R(t)^{\f{3\g-3-\dl}2-k}\nabla_y^{3-k}Zu|_{L^2(Q_T)}^{\f12}\bigg)|R(t)^{\f{3\g-3-\dl}2}
\nabla_yZu|_{L^{\infty}(Q_T)}^{\f12}.\qquad \quad\text{(By Lemma 2.2)}\no
\een

(v)\ben
&&|R(t)^{\f{6\g-6-\dl}2}\nabla_yZD_tu|_{L^4(Q_T)}\no\\
&\le&C|R(t)^{\f{6\g-6-\dl}2}\nabla_yE(ZD_tu)|_{L^4(Q_T)}\no\\
&\le&C|R(t)^{\f{6\g-6-\dl}2}\nabla_y^{2}E(ZD_tu)|_{L^2(Q_T)}^{\f12}|R(t)^{\f{6\g-6-\dl}2}E(ZD_tu)|_{L^{\infty}(Q_T)}^{\f12}
\quad\text{(Applying (2.15) for $E(ZD_tu)$)} \no\\
&\le&C\bigg(\ss_{k=0}^2|R(t)^{\f{6\g-6-\dl}2-k}\nabla_y^{2-k}ZD_tu|_{L^2(Q_T)}^{\f12}\bigg)
|R(t)^{\f{6\g-6-\dl}2}ZD_tu|_{L^{\infty}(Q_T)}^{\f12}.\qquad\quad\text{(By Lemma 2.2)}\no
\een

(vi) and (vii) follow from the definition of $Z-$fields in (1.14) and Lemma 2.3. $\hfill\square$

\section{Reformulation of  problem (\ref{potential}) with (1.16) and (\ref{initial condition})}

Firstly, we state a local solvability result on problem (\ref{potential}) with (\ref{bd condition})
and (\ref{initial condition}).

\begin{lemma}\label{local-existence}
 There exists a $T_0>0$ such that the problem
(\ref{potential}) with (1.16) and (\ref{initial condition}) has a local solution
$\Phi(t,x)\in C([0, T_0], H^5(S_t))\cap C^1([0, T_0], H^4(S_t))$ with $S_t=\{x: |x|\le R(t)\}$. Moreover,
$$||\Phi(t,x)-\hat{\Phi}(t,r)||_{C([0, T_0], H^5(S_t))}+||\Phi(t,x)-\hat{\Phi}(t,r)||_{C^1([0, T_0], H^4(S_t))}\le C\ve,$$
where $\hat\Phi(t,r)$ is given in (\ref{background-Phi}).
\end{lemma}

{\bf Proof.} The quasilinear equation (\ref{potential}) is
strictly hyperbolic with respect to $t$. Thus, by the standard
Picard iteration  as in [17], one can derive that Lemma \ref{local-existence} holds.
\qquad \qquad\quad\quad $\square$

Next, we reformulate (\ref{potential}) with (1.16)-(\ref{initial condition}).

Let
$\dot{\Phi}=\Phi-\hat{\Phi}$. Then  (\ref{potential}) can be reduced to
\beq\label{linear-equ}
\mathcal{L}\dot{\Phi}=\dot
f\qquad\text{in $\O$},
\eeq
where
\beq
\bec
\mathcal{L}\dot\Phi=\p_t^2\dP+2\ss_{i=1}^3\p_i\hat\Phi\p_{ti}^2\dP+\ss_{i,j=1}^3\p_i\hP\p_j\hP\p_{ij}^2\dP
+2\ss_{i=1}^3\p_{ti}^2\hP\p_i\dP -\hat{c}^2\Delta\dP\\
\qquad\quad +\ds\f{3L(\g-1)}{R(t)}(\p_t\dP+\ss_{i=1}^3\p_i\hat\Phi\p_i\dP),\\
\dot{f}=\ss_{i=1}^3f_{0i}\p_{ti}^2\dP+\ss_{1\le i\neq j\le3}f_{ij}\p_{ij}^2\dP+\ss_{i=1}^3f_{ii}\p_i^2\dP+f_0,
\eec\label{L-in-t-x}
\eeq
with
$$
\bec
f_{0i}=-2\p_i\dP,\\
f_{ij}=-\p_i\dP\p_j\dP-2\p_i\hat\Phi\p_j\dP,\\
f_{ii}=-(\p_i\dP)^2-(\g-1)(\p_t\dP+\ss_{j=1}^3\p_j\hat\Phi\p_j\dP+\ds\f{1}{2}\ss_{j=1}^3(\p_j\dP)^2)-2\p_i\hat\Phi\p_i\dP,\\
f_0=-\ds\f{(3\g-1)L}{2R(t)}\ss_{i=1}^3(\p_i\dP)^2.
\eec$$

For later analysis, we use $Z$-fields to rewrite $\mathcal{L}\dP$ as follows:
\ben
\mathcal{L}\dP&=&\p_t^2\dP+2\p_r\hat\Phi\p_{tr}^2\dP+((\p_r\hP)^2-\hat{c}^2)\p_r^2\dP-\f{\hat{c}^2}{r^2}\ss_{i=1}^3Z_i^2\dP +\f{\g-1}{r}(r\p_r^2\hP+2\p_r\hP)\p_t\dP\no\\
&&+\f{1}{r}\big(2r\p_{tr}^2\hP+(\g+1)r\p_r\hP\p_r^2\hP+2(\g-1)(\p_r\hP)^2-2\hat{c}^2\big)\p_r\dP.\label{linear-equ-z}
\een

On the lateral boundary $\p\O$ of $\O$, $\dP$ satisfies
\beq\label{linear-boundary}
\p_r\dP=0.
\eeq

In addition, we have the following initial data of $\dP$ from (\ref{initial condition})
\beq\label{linear-initial-data}
\dP(0,x)=\ve\Phi_0(x),\q \p_t(0,x)=\ve\Phi_1(x).
\eeq

{\bf Remark 3.1.} {\it By (2.3) or (2.4), we know that $\mathcal{L}\dot\Phi$ in (3.2) has the following form at the $t-$axis:
$$\p_t^2\dP-\ds\f{\g}{R^{3(\g-1)}(t)}(\p_1^2+\p_2^2+\p_3^2)\dP+\ds\f{3L(\g-1)}{R(t)}\p_t\dP,\eqno(3.6)$$
which is strictly hyperbolic but degenerate as
$t\to\infty$. On the other hand, the operator $\mathcal{L}=\p_t^2+\ds\f{2Lr}{R(t)}\p_{tr}^2
+\ds\sum_{i,j=1}^3\ds\f{L^2x_ix_j}{R^2(t)}\p_{ij}^2-\ds\f{\g}{R^{3(\g-1)}(t)}\Delta+
\ds\f{3L(\g-1)}{R(t)}(\p_t+\ds\f{Lr}{R(t)}\p_r)$ in (3.2) is different from
[24], where the corresponding linear operator is
$\p_t^2-\ds\f{1}{R^{2(\g-1)}(t)}(\p_1^2+\p_2^2)+\ds\f{2L(\g-1)}{R(t)}\p_t$.
Recently, with respect to the semilinear wave equations
$$\p_t^2u-\Delta u+\ds\f{\mu}{(1+t)^{\al}}\p_tu=f(u),\eqno(3.7)$$
where $\mu>0$ and $\al>0$ are suitable constants, there
have been extensive works on the global existence or blowup results for different nonlinear function $f(u)$, cf.
[9], [21-22] and the references therein. However, for the critical exponent $\al=1$ and near the critical value
$\mu=1$, there are still some open questions on the blowup or global existence of solutions to (3.7). Here,
(3.6) corresponds to the critical case of (3.7) when $L=1$ and $\g$ is close to $\f43$.}

\section{The first-order weighted energy estimate and reformulation of (3.2)-(3.3)}
In this section, we derive the weighted energy estimate of
$\nabla_{t,x}\dot\Phi$ for the linear part (\ref{linear-equ-z})
together with (\ref{linear-boundary})-(\ref{linear-initial-data}).

Set $\O_T=\O\cap\{0<t<T\}$ and $B_T=\p\O\cap\{0<t<T\}$. Then, we have

\begin{theorem} Let $\dot\Phi\in C^{2}(\bar \Omega_T)$ satisfy the boundary
condition (\ref{linear-boundary}) and initial data condition (\ref{linear-initial-data}). Then for $0<\g<\f53$, there exists a multiplier
$\mathcal{M}\dot\Phi=R(t)^{\mu}a(t)D_t\dP$ such that for fixed
constant $\mu=6\g-9$ we have
\ben
&&R(T)^{\mu}\int_{S_T}(D_t\dot\Phi)^2dS+T^{\mu-3(\g-1)}\int_{S_T}(\nabla_x\dP)^2dS\no\\
&&\quad +C\int_{\Omega_T}\big(R(t)^{\mu-1-\dl}(D_t\dot\Phi)^2+R(t)^{\mu-1-3(\g-1)}(\nabla_x\dot\Phi)^2\big)dtdx\no\\
&&\le \int_{\Omega_T}\mathcal{L}\dot\Phi\cdot\mathcal{M}\dot\Phi
dtdx
+C\ve^2,\label{first-order-estimate}
\een
where $D_t=\p_t+\ss_{i=1}^3\p_i\hP\p_i=\p_t+\ds\f{Lr}{R(t)}\p_r$
is the material derivative,
$C>0$ is a generic positive constant depending only on the initial data, and $\dl>0$ is a small fixed constant.
\end{theorem}

\begin{remark} The choice of $\mu=6\g-9$ in (\ref{first-order-estimate}) is
 necessary because of the following two reasons:
First, to guarantee the positivity of $III$ in (4.2),  $\mu\le 6\g-9$; Second, by the Bernoulli law (\ref{Bernoulli}),
we have $c^2(\rho)=
c^2(\hat\rho)-(\g-1)D_t\dP-\f{\g-1}{2}|\nabla_x\dP|^2$.
Notice that only the estimate of $|D_t\dP|\le C\ve  R(t)^{-\f{\mu+3}{2}}$ and $|\nabla_x\dot\Phi|\le
C\ve R(t)^{-\f{\mu-3(\g-1)+3}{2}}$ can be obtained as shown  in Sections 5 and
6. On the other hand, $c^2(\hat\rho(t))=\g R(t)^{-3(\g-1)}$ holds. Therefore, in order to guarantee the absence of vacuum for any finite time $t$ in $\O$, we need to choose
the constant $\mu$ such that $-\ds\f{\mu+3}{2}\le
-3(\g-1)$, that is, $\mu\ge 6\g-9.$ In combination,  $\mu=6\g-9$.
\end{remark}

{\bf Proof.} Choosing $\mathcal{M}\dP=R(t)^{\mu}(a(t,r)\p_t\dP+b(t,r)\p_r\dP)$, where the non-negative functions $a(t,r)$ and $b(t,r)$
will be determined later, then
\beq\int_{\O_T}\mathcal{L}\dP\cdot\mathcal{M}\dP dtdx=I+II+III,\eeq
where
\begin{align*}
&I=\int_{B_T}\ds\f{R(t)^{\mu}}{2\sqrt{1+L^2}}(La(t,r)-b(t,r))\bigg((\p_t\dP)^2-\f{\hat{c}^2}{r^2}\ss_{i=1}^3(Z_i\dP)^2\\ &\quad +2((L^2-\hat{c}^2)a(t,r)-Lb(t,r))\p_t\dP\p_r\dP+(L(L^2-\hat{c}^2)a(t,r)-(L^2+\hat{c}^2)b(t,r))(\p_r\dP)^2\bigg)dS,\\
&II\equiv II_1-II_2,\\
&III=\int_{\O_T}\biggl(A(t,r)(\p_t\dP)^2+B(t,r)\p_t\dP\p_r\dP+C(t,r)(\p_r\dP)^2+\ds\f{1}{r^2}D(t,r)\ss_{i=1}^3(Z_i\dP)^2\biggr)dtdx,
\end{align*}
with
\begin{align*}
&II_1=\int_{S_T}R(t)^{\mu}\bigg(\ds\f{1}{2}a(t,r)(\p_t\dP)^2+b(t,r)\p_t\dP\p_r\dP+(b(t,r)\p_r\hP-\f12((\p_r\hP)^2-\hat{c}^2)a(t,r))(\p_r\dP)^2\\
&\qquad +\ds\f{\hat{c}^2}{2r^2}a(t,r)\ss_{i=1}^3(Z_i\dP)^2\bigg)dS,\\
&II_2=\int_{S^0}\bigg(\ds\f{1}{2}a(t,r)(\p_t\dP)^2+b(t,r)\p_t\dP\p_r\dP
+(b(t,r)\p_r\hP-\f12((\p_r\hP)^2-\hat{c}^2)a(t,r))(\p_r\dP)^2\\
&\qquad +\ds\f{\hat{c}^2}{2r^2}a(t,r)\ss_{i=1}^3(Z_i\dP)^2\bigg)dS,\\
&A(t,r)=-\ds\f12\p_t\bigg(R(t)^{\mu}a(t,r)\bigg)-R(t)^{\mu}\p_r\bigg(a(t,r)\p_r\hP\bigg)-2R(t)^{\mu}r^{-1}a(t,r)\p_r\hP\\
&\q\q+r^{-2} \p_r\bigg(\ds\f12r^2R(t)^{\mu}b(t,r)\bigg)+(\g-1)(\p_r^2\hP+\f2r\p_r\hP)R(t)^{\mu}a(t,r),\\
&B(t,r)=-R(t)^{\mu}\p_r\bigg(((\p_r\hP)^2-\hat{c}^2)a(t,r)\bigg)
-\ds\f2rR(t)^{\mu}a(t,r)\bigg((\p_r\hP)^2-\hat{c}^2\bigg)-\p_t\bigg(R(t)^{\mu}b(t,r)\bigg)\\
&\q\q+R(t)^{\mu}a(t,r)\bigg(2\p_{tr}\hP+(\g+1)\p_r\hP\p_r^2\hP+\f{2(\g-1)}{r}(\p_r\hP)^2-\f2r\hat{c}^2\bigg) \\ &\q\q+(\g-1)\bigg(\p_r^2\hP+\f2r\p_r\hP\bigg)R(t)^{\mu}b(t,r),\\
&C(t,r)=\ds\f12\p_t\bigg(((\p_r\hP)^2-\hat{c}^2)R(t)^{\mu}a(t,r)\bigg)
-r^{-2}\p_r\bigg(\f12R(t)^{\mu}b(t,r)r^2((\p_r\hP)^2-\hat{c}^2)\bigg)\\
&\q\q-\p_t\bigg(R(t)^{\mu}b(t,r)\p_r\hP\bigg)+R(t)^{\mu}b(t,r)\bigg(2\p_{tr}^2\hP+(\g+1)\p_r\hP\p_r^2\hP
+\ds\f{2(\g-1)}{r}(\p_r\hP)^2-\f2r\hat{c}^2\bigg),\\
&D(t,r)=-\ds\f1{2}\p_t\bigg(\hat{c}^2R(t)^{\mu}a(t,r)\bigg)-\p_r\bigg(\f12R(t)^{\mu}b(t,r)\hat{c}^2\bigg).
\end{align*}
In view of the boundary condition (\ref{linear-boundary}), we have
$$
I=\int_{B_T}\f{1}{2\sqrt{1+L^2}}(La(t,R(t))-b(t,R(t)))\big((\p_t\dP)^2-\f{\hat{c}^2}{r^2}\ss_{i=1}^3(Z_i\dP)^2\big)dS.
$$
To guarantee $I\ge 0$, it requires that on the boundary $r=R(t)$,
\beq\label{a-b-boundary}
b(t,R(t))=La(t,R(t)).
\eeq
In this case,
\beq\label{estimate-on-boundary}
I=0.\eeq
Next, we consider $II$. To fulfill $II_1>0$, it requires on $t=T$
\beq
\bec
a(t,r)>0,\\
b(t,r)^2-a(t,r)\bigg(2b(t,r)\p_r\hP-((\p_r\hP)^2-\hat{c}^2)a(t,r)\bigg)\le 0.
\eec
\eeq
This means that on $\{t=T\}$,
\beq\label{a-b-T}
\bec
a(t,r)>0,\\
\p_r\hP-\hat{c}\le \f{\ds b(t,r)}{\ds a(t,r)}\le\p_r\hP+\hat{c}.
\eec
\eeq
Thus, combining (\ref{a-b-boundary}) and (\ref{a-b-T}) yields
\beq
b(t,r)=\p_r\hP\cdot a(t,r).
\eeq
On the other hand, by $\hat{c}^2=\g R(t)^{-2(\g-1)}$, we have
\ben
II_1&=&\int_{S_T}\f12R(t)^{\mu}a(T,r)\bigg((\p_t\dP+\p_r\hP\p_r\dP)^2
+\hat{c}^2\big((\p_r\dP)^2+\f{1}{r^2}\ss_{i=1}^3(Z_i\dP)^2\big)\bigg)dS\no\\
&\ge&C\int_{S_T}a(T,r)\bigg(R(T)^{\mu}(D_t\dP)^2
+R(T)^{\mu-3(\g-1)}\big((\p_r\dP)^2+\f{1}{r^2}\ss_{i=1}^3(Z_i\dP)^2\big)\bigg)dS.\label{estimate-on-T}
\een
It follows from Lemma \ref{local-existence} that
\beq
II_2\le C\ve^2.
\eeq
Finally, we deal with $III$. In fact, we only need to choose $a(t,r)\equiv a(t)$.
In this case, direct computation yields
$$\bec
A(t,r)=\ds\f12R(t)^{\mu-1}\bigg(L(6\g-9-\mu)a(t)-R(t)a'(t)\bigg),\\
B(t,r)=LrR(t)^{\mu-2}\bigg(L(6\g-9-\mu)a(t)-R(t)a'(t)\bigg),\\
C(t,r)=\ds\f12R(t)^{\mu-3\g-3}\bigg(\big(L\g(-4+3\g-\mu) R(t)^5+L^2r^2R(t)^{3\g}(-9+6\g-\mu)\big)a(t)\\
\q\q\q\q\q-\big(L^2r^2R(t)^{1+3\g}+\g R(t)^6\big)a'(t)\bigg),\\
D(t,r)=\ds\f{\g}{2}R(t)^{\mu-3\g+2}\bigg(L(3\g-4-\mu)a(t)-R(t)a'(t)\bigg).
\eec$$
In order to have $III>0$, we choose
\beq\label{positivity}
A(t,r)>0,\q B(t,r)^2-4A(t,r)C(t,r)<0,\q D(t,r)>0.
\eeq
From this, we naturally set
\beq
\mu=6\g-9,\q a(t)>0\q \text{and}\q a'(t)<0.
\eeq
If we choose
\beq
a(t)=1+R(t)^{-\dl}\q \text{with}\q \dl>0,
\eeq
then
\ben
III&=&\int_{\O_T}\biggl\{\f{\dl}{2}R(t)^{\mu-1-\dl}(\p_t\dP+\p_r\hP\p_r\dP)^2\no\\
&&\quad +\f{\g}{2}R(t)^{\mu-1-3(\g-1)}\big(L(5-3\g)a(t)-R(t)a'(t)\big)
\bigg((\p_r\dP)^2+\f{1}{r^2}\ss_{i=1}^3(Z_i\dP)^2\bigg)\biggr\}dtdx\no\\
&\ge& C\int_{\O_T}\biggl\{R(t)^{\mu-1-\dl}(D_t\dP)^2+R(t)^{\mu-1-3(\g-1)}\bigg((\p_r\dP)^2
+\f{1}{r^2}\ss_{i=1}^3(Z_i\dP)^2\bigg)\biggr\}dtdx.\label{estimate-on-domain}
\een
Substituting (4.4), (4.8)-(4.9) and (4.13) into (4.2) yields
\ben
&&R(T)^{\mu}\int_{S_T}(D_t\dot\Phi)^2dS+R(T)^{\mu-3(\g-1)}\int_{S_T}\big((\p_r\dP)^2+\f{1}{r^2}\ss_{i=1}^3(Z_i\dP)^2\big)dS\no\\
&&\quad +C\int_{\Omega_T}\bigg(R(t)^{\mu-1-\dl}(D_t\dot\Phi)^2+R(t)^{\mu-1-3(\g-1)}\big((\p_r\dP)^2
+\f{1}{r^2}\ss_{i=1}^3(Z_i\dP)^2\big)\bigg)dtdx\no\\
&&\le \int_{\Omega_T}\mathcal{L}\dot\Phi\cdot\mathcal{M}\dot\Phi
dtdx
+C\ve^2.
\een
This together with Remark 2.2 give (\ref{first-order-estimate}),
and it completes the proof of the lemma.\hfill$\square$

Note that the material derivative $D_t$ plays a crucial role in the energy estimate (4.1). Then it is necessary
that
we rewrite equation (\ref{L-in-t-x}) and formula (\ref{linear-equ-z}) by using the operator $D_t$ as follows:

\beq
\bec
\mathcal{L}\dot\Phi=D_t^2\dP-\hat{c}^2\Delta\dP+\ds\f{3L(\g-1)}{R(t)}D_t\dP,\\
\dot{f}=\ss_{i=1}^3\tilde{f}_{0i}\p_i D_t\dP+\ss_{1\le i\neq j\le3}\tilde{f}_{ij}\p_{ij}^2\dP+\ss_{i=1}^3\tilde{f}_{ii}\p_i^2\dP+\tilde{f}_0,
\eec\label{L-in-Dt}
\eeq
with
$$
\bec
\tilde{f}_{0i}=-2\p_i\dP,\\
\tilde{f}_{ij}=-\p_i\dP\p_j\dP,\\
\tilde{f}_{ii}=-(\p_i\dP)^2-(\g-1)(D_t\dP+\ds\f{1}{2}\ss_{j=1}^3(\p_j\dP)^2),\\
\tilde{f}_0=-\ds \f{(3\g-5)L}{2R(t)}\ss_{i=1}^3(\p_i\dP)^2.
\eec$$
Especially, in the spherical coordinates, (4.15) has the form
\beq\label{L-in-Z}
\mathcal{L}\dot\Phi=D_t^2\dP-\hat{c}^2\big(\p_r^2\dP+\f{1}{r^2}\ss_{i=1}^3Z_i^2\dP+\ds\f{2}{r}\p_r\dP\big)+\ds\f{3L(\g-1)}{R(t)}D_t\dP
=\dot f,
\eeq
where
\begin{align*}
\dot{f}&=-2\p_r\dP D_t\p_r\dP+\f{2}{r^2}\ss_{i=1}^3Z_i\dP D_tZ_i\dP-\bigg((\g-1)D_t\dP+\f{\g+1}{2}(\p_r\dP)^2+\f{\g-1}{2r^2}\ss_{i=1}^3(Z_i\dP)^2\bigg)\p_r^2\dP\no\\
&-\f{2}{r^2}\p_r\dP\ss_{i=1}^3Z_i\dP\p_rZ_i\dP -\f{\g-1}{r^2}\bigg(D_t\dP+\f{1}{2}(\p_r\dP)^2+\f{1}{2r^2}\ss_{i=1}^3(Z_i\dP)^2\bigg)\ss_{i=1}^3Z_i^2\dP\no\\
& -\f{1}{r^4}\ss_{i,j=1}^3Z_i\dP Z_j\dP Z_iZ_j\dP-\f{2(\g-1)}{r}\p_r\dP D_t\dP-\f{(3\g-1)L}{2R(t)}(\p_r\dP)^2-\f{3(\g-1)L}{2r^2R(t)}\ss_{i=1}^3(Z_i\dP)^2\no\\
&+\f{\g-1}{r}\big((\p_r\dP)^2 -\f{1}{r^2}\ss_{i=1}^3(Z_i\dP)^2\big)-\ss_{i,j=1}^3\f{C_{ij}}{r^3}\big(\f{Lr}{R(t)}+\p_r\dP\big)Z_i\dP Z_j\dP-\ss_{i,j,k=1}^3\f{C_{ijk}}{r^4}Z_i\dP Z_j\dP Z_k\dP\label{f-in-Z}.
\end{align*}
Here, we point out that the precise expression of $\dot f$ is useful in the higher order
energy estimates of $\dot\Phi$, cf. Setion 5.

\section{Higher order weighted energy estimates of $\dP$}

In this section, we will derive the higher order energy
estimates of  $\dot\Phi$ to (\ref{linear-equ}) with (\ref{linear-boundary})-(\ref{linear-initial-data}). For this, we need to take care of the difficulties
coming from
the Neumann boundary condition (\ref{linear-boundary}),
the asymptotic degeneracy of some coefficients in (\ref{L-in-Z}), and the different decay rates of $D_t\dP$ and
$\nabla_x\dP$.

\begin{theorem}
Let $\dP\in C^4(\bar{\O}_T)$ be the solution to (\ref{linear-equ}) with (\ref{linear-boundary})-(\ref{linear-initial-data}),
and assume with some $\delta>0$,
\beq\label{apriori-bound}
\bec
|\nabla_x\dot\Phi|\le M\ve R(t)^{-3(\g-1)+\f{\dl}2},\quad |R(t)^{l-1}D_t^l\dot\Phi|\le M\ve R(t)^{-3(\g-1)},\quad\text{for $l=1,2$},\\
|R(t)\nabla_xD_t\dot\Phi|\le M\ve R(t)^{-3(\g-1)+\f{\dl}2},\quad |R(t)\nabla_x^2\dot\Phi|\le M\ve R(t)^{-\f{3(\g-1)-\dl}2},
\eec
\eeq
in addition,  for $r>\ds\f{1}{3}R(t)$, further assume
\beq\label{apriori-better-bound}
|Z\dot\Phi|\le M\ve R(t)^{1-3(\g-1)},|ZD_t\dot\Phi|\le M\ve R(t)^{-3(\g-1)},|\nabla_xZ\dot\Phi|\le M\ve R(t)^{-\f{3(\g-1)}2}.
\eeq
Then for sufficiently small $\ve>0$ and $0\le k\le2$, we have
\begin{align}
&\int_{S_T}\bigg(R(T)^{\mu+2k}(\nabla_{t,x}^kD_t\dot\Phi)^2+R(T)^{\mu-3(\g-1)+2k}(\nabla_{t,x}^k\nabla_x\dot\Phi)^2\bigg)dS\no\\
&\quad +\int_{\Omega_T}\bigg(R(t)^{\mu-1-\dl+2k}(\nabla_{t,x}^kD_t\dot\Phi)^2+R(t)^{\mu-1-3(\g-1)+2k}(\nabla_{t,x}^k\nabla_x\dot\Phi)^2\bigg)dtdx\no\\
&\le C\ve^2,\label{estimate-to-3}
\end{align}
and
\begin{align}
&\int_{S_T}\bigg(R(T)^{\mu-\dl+6}(\nabla_{t,x}^3D_t\dot\Phi)^2+R(T)^{\mu-\dl-3(\g-1)+6}(\nabla_{t,x}^3\nabla_x\dot\Phi)^2\bigg)dS\no\\
&\quad +\int_{\Omega_T}\bigg(R(t)^{\mu+5-\dl}(\nabla_{t,x}^3D_t\dot\Phi)^2+R(t)^{\mu-\dl+5-3(\g-1)}(\nabla_{t,x}^3\nabla_x\dot\Phi)^2\bigg)dtdx\no\\
&\le C\ve^2,\label{estimate-to-4}
\end{align}
where $\mu=6\g-9$, $0<\dl\le\f{3(\g-1)}{5}$, $C>0$ is independent of $M$,
and the domains $\O_T, S_T$ are defined in Section 4.
\end{theorem}

In order to prove Theorem 5.1, we will apply the induction  on $k$ in (\ref{estimate-to-3})-(\ref{estimate-to-4}) to establish
the following estimates respectively:

(i) $D_tS^k\dP$ and $\nabla_xS^k\dP$ with $S=S_0^{l_1}Z^{l_2}(S_0=R(t)D_t)$ and $1\le k=l_1+l_2\le 3$ (in this case, all
the tangent derivatives of $\na_x\dP$
up to the third order are estimated, where the tangent derivative
means the one of boundary $B_T$ );

(ii) $D_tS_1\dP$ and $\nabla_x S_1\nabla_x\dP$ with $S_1=r\p_r$ (in this case, together with the case $k=1$ in (i),
all the second order derivatives $\na_{t,x}^2\dot\Phi$ are analyzed);

(iii) $D_tSS_1\dP$, $\nabla_xSS_1\dP$, $D_tS_1^2\dP$ and $\nabla_xS_1^3\dP$ (in this case, together with the case $k=2$ in (i),
all the estimates of third order derivatives $\na_{t,x}^2\dot\Phi$ are given);

(iv) $D_tS^2S_1\dP$, $\nabla_xS^2S_1\dP$, $D_tSS_1^2\dP$, $\nabla_xSS_1^2\dP$, $D_tS_1^3\dP$ and $\nabla_xS_1^3\dP$ (in this case,
together with the case $k=3$ in (i),
all the fourth order derivatives $\na_{t,x}^4\dot\Phi$ are estimated).

These estimates will be given in Lemma 5.2-Lemma 5.5 respectively.

Firstly, we establish the tangent derivative estimates of $\dP$ under the suitable induction assumption. Set $S_0=R(t)D_t$ and
$S^m=S_0^{l_1}Z^{l_2}$ with $m=l_1+l_2$, which  are tangent to the boundary $\p\O$. Then we have

\begin{lemma}{\bf (Tangent derivative estimates)} Under the assumptions of Theorem 5.1, if (\ref{estimate-to-3}) holds for $0\le k \le m$ with $1\le m\le2$, then
\begin{align}
&R(T)^{\mu}\int_{S_T}(D_tS^m\dP)^2dS+R(T)^{\mu-3(\g-1)}\int_{S_T}(\nabla_xS^m\dP)^2dS\no\\
&\quad +\int_{\Omega_T}\bigg(R(t)^{\mu-1-\dl}(D_tS^m\dP)^2+R(t)^{\mu-1-3(\g-1)}(\nabla_xS^m\dP)^2\bigg)dtdx\no\\
&\le C\ve^2+C\ve\int_{\Omega_T}\ss_{l=0}^m\bigg(R(t)^{\mu-1-\dl+2l}(\nabla_{t,x}^lD_t\dot\Phi)^2
+R(t)^{\mu-1-3(\g-1)+2l}(\nabla_{t,x}^l\nabla_x\dot\Phi)^2\bigg)dtdx,\label{estimate-tangent}
\end{align}
where $0<\dl\le \ds\f{3(\g-1)}{5}$.

Especially, for $m=0$, the following estimate holds
\begin{align}
&R(T)^{\mu}\int_{S_T}(D_t\dP)^2dS+R(T)^{\mu-3(\g-1)}\int_{S_T}(\nabla_x\dP)^2dS\no\\
&\quad +\int_{\Omega_T}\bigg(R(t)^{\mu-1-\dl}(D_t\dP)^2+R(t)^{\mu-1-3(\g-1)}(\nabla_x\dP)^2\bigg)dtdx\no\\
&\le C\ve^2.\no
\end{align}

For $m=3$, if (5.2)-(5.3) hold, then
\begin{align}
&R(T)^{\mu}\int_{S_T}(D_tS^3\dP)^2dS+R(T)^{\mu-3(\g-1)}\int_{S_T}(\nabla_xS^3\dP)^2dS\no\\
&\quad +\int_{\Omega_T}\bigg(R(t)^{\mu-1-\dl}(D_tS^3\dP)^2+R(t)^{\mu-1-3(\g-1)}(\nabla_xS^3\dP)^2\bigg)dtdx\no\\
&\le C\ve^2+C\ve\int_{\Omega_T}\ss_{l=0}^2\bigg(R(t)^{\mu-1-\dl+2l}(\nabla_{t,x}^lD_t\dot\Phi)^2
+R(t)^{\mu-1-3(\g-1)+2l}(\nabla_{t,x}^l\nabla_x\dot\Phi)^2\bigg)dtdx\no\\
&\quad +C\ve\int_{\Omega_T}\bigg(R(t)^{\mu+5-\dl}(\nabla_{t,x}^3D_t\dot\Phi)^2+R(t)^{\mu+5-3(\g-1)-\dl}
(\nabla_{t,x}^3\nabla_x\dot\Phi)^2\bigg)dtdx.
\end{align}

\end{lemma}

\begin{remark}
For the case when $m=0$ in (5.5), we do not require any induction assumption.
\end{remark}
\begin{remark}
It is noted that the normal derivatives of $\dP$
are  included on the right hand side of (\ref{estimate-tangent}) and (5.6), which implies that we have not
obtained the close estimates on the tangent derivative estimates of $\dP$. However,
since the coefficients of normal derivatives of $\dP$  in (\ref{estimate-tangent}) are small, then together
with some subsequent normal  derivative estimates, we can derive (\ref{estimate-to-3}).
\end{remark}

{\bf Proof.} Note that on $B_T$
\beq\label{tangent-bd-condition}
S^m\bigg(\ss_{i=1}^3 x_i\p_i\dP\bigg)=S^m\big(r\p_r\dP\big)=0.
\eeq

This, together with Theorem 4.1 and (\ref{linear-initial-data}), yields
\ben
&&R(T)^{\mu}\int_{S_T}(D_tS^m\dP)^2dS+R(T)^{\mu-3(\g-1)}\int_{S_T}(\nabla_xS^m\dP)^2dS\no\\
&&\quad +\int_{\Omega_T}\bigg(R(t)^{\mu-1-\dl}(D_tS^m\dP)^2+R(t)^{\mu-1-3(\g-1)}(\nabla_xS^m\dP)^2\bigg)dtdx\no\\
&&\le C\int_{\O_T}\mathcal{L}S^m\dP\cdot\mathcal{M}S^{m}\dP dtdx+C\ve^2.\label{estimate-S_0^m}
\een

We divide the estimate (5.5) into three cases: (1) $S^m=S_0^m$, (2) $S^m=Z^m$, (3) $S^m=S_0^{l_1}Z^{l_2}
(1\le l_1,l_2\le m-1, l_1+l_2=m)$.

\vskip 0.3cm
{\bf Case (1)\quad $S^m=S_0^m$.}
\vskip 0.3cm

We now derive an explicit representation of $\mathcal{L}S_0^m\dP$ for later use.
By direct computation, we have
\beq
S_0\mathcal{L}=\mathcal{L}S_0-2L\mathcal{L}+3(\g-1)\hat{c}^2\Delta.
\eeq
By induction, for $1\le m\le3$, we obtain
\beq
\mathcal{L}S_0^m=S_0^m\mathcal{L}+B_{1m}+B_{2m}
\eeq
with
\beq
B_{1m}=\ss_{0\le l\le m-1}C_{lm}S_0^l\mathcal{L}+\ss_{0\le l\le m-2}C'_{lm}S_0^l(\hat{c}^2\Delta),~~B_{2m}=-3m(\g-1)S_0^{m-1}(\hat{c}^2\Delta),
\eeq
where $B_{2m}\dP$ contains the $(m+1)-$th order (the highest order) derivatives of $\dot\Phi$,
but $B_{1m}\dP$ only includes the terms $\na_{t,x}^{\al}\dP$ with $|\al|\le m$ (the lower order derivatives of $\dP$)
and the term $\na_{t,x}^{m+1}\dP$ with small coefficients.

In addition, from equation (\ref{L-in-Dt}), we have that for $0\le m\le3$
\beq
S_0^m\mathcal{L}\dP=I_1^m+I_2^m+I_3^m,
\eeq
where
\ben
&&I_1^m=\ss_{i=1}^3\tf_{0i}\p_iD_tS_0^m\dP+\ss_{1\le i\neq j\le3}\tf_{ij}\p_{ij}^2S_0^m\dP+\ss_{i=1}^3\tf_{i}^2\p_{i}^2S_0^m\dP,\no\\
&&I_2^m=\ss_{i=1}^3\tf_{0i}[S_0^m,\p_iD_t]\dP+\ss_{0\le i\neq j\le3}\tf_{ij}[S_0^m,\p_{ij}^2]\dP
+\ss_{i=1}^3\tf_{ii}[S_0^m,\p_{i}^2]\dP,\no\\
&&I_3^m=\ss_{1\le l\le m}C_{lm}\bigg\{\ss_{l_1+l_2=l,l_1\ge1}\tilde{C}_{l_1l_2}\bigg(\ss_{i=1}^3(S_0^{l_1}\tf_{0i})S_0^{l_2}(\p_iD_t\dP) +\ss_{1\le i,j\le3}(S_0^{l_1}\tf_{ij})S_0^{l_2}(\p_{ij}^2\dP)\bigg)\bigg\}+S_0^m\tf_0.\no
\een

Based on the above preparation, we now treat $\int_{\O_T}\mathcal{L}S_0^m\dP\cdot\mathcal{M}S_0^m\dP dtdx$ on the right hand side of (\ref{estimate-S_0^m})
as follows. This procedure is divided into  five parts.

\vskip 0.3cm
{\bf Part 1. Estimate on $\int_{\O_T}I_1^m\cdot\mathcal{M}S^m\dP dtdx.$}
\vskip 0.3cm

Note that for $C^1$-smooth functions $g$, one has
$$
D_tg=\p_tg+\ss_{i=1}^3\p_i\big(\f{L}{R(t)}x_ig\big)-\f{3L}{R(t)}g.
$$
Then
\beq
\int_{\O_T}D_tgdtdx=\int_{S_T}gdS-\int_{S^0}gdS-\int_{\O_T}\f{3L}{R(t)}gdtdx.
\eeq
In addition, for $m\le3$,  direct computation yields

\begin{align}
&I_1^m\cdot\mathcal{M}S_0^m\dP\no\\
&=D_t\bigg(-\f12a(t)R(t)^{\mu}\ss_{i,j=1}^3\big(\tf_{ij}\p_iS_0^m\dP\p_jS_0^m\dP\big)\bigg)\no\\
&\quad+\ss_{i=1}^3\p_i\bigg(\f12a(t)R(t)^{\mu}\tf_{0i}(D_tS_0^m\dP)^2+a(t)R(t)^{\mu}D_tS_0^m\dP\ss_{j=1}^3\tf_{ij}\p_jS_0^m\dP\bigg)\no\\
&\quad-\f12a(t)R(t)^{\mu}(D_tS_0^m\dP)^2\ss_{i=1}^3\p_if_{0i}-a(t)LR(t)^{\mu-1}\ss_{i,j=1}^3\big(\tf_{ij}\p_iS_0^m\dP\p_jS_0^m\dP\big)\no\\ &\quad-a(t)R(t)^{\mu}D_tS_0^m\dP\ss_{i,j=1}^3\big(\p_i\tf_{ij}\p_jS_0^m\dP\big)
+\f12\ss_{i,j=1}^3\big(D_t(a(t)R(t)^{\mu}\tf_{ij})\p_iS_0^m\dP\p_jS_0^m\dP\big).
\end{align}

On the other hand, by the expressions of $\tf_{ij}$ and $\tf_{0i}$, and (\ref{tangent-bd-condition}) on $B_T$, we observe that
\begin{align}
&\ss_{i=1}^3x_i\cdot\bigg(\f12a(t)R(t)^{\mu}\tf_{0i}(D_tS_0^m\dP)^2+a(t)R(t)^{\mu}D_tS_0^m\dP\ss_{j=1}^3\tf_{ij}\p_jS_0^m\dP\bigg)\no\\
&=-a(t)R(t)^{\mu}(D_tS_0^m\dP)^2\ss_{i=1}^3\big(x_i\p_i\dP\big)-a(t)R(t)^{\mu}D_tS_0^m\dP\ss_{i,j=1}^3\big(x_i\p_i\dP\p_j\dP\big)\no\\
&\quad -(\g-1)a(t)R(t)^{\mu}D_tS_0^m\dP\big(D_t\dP+\f12\ss_{k=1}^3(\p_k\dP)^2\big)\ss_{i=1}^3x_i\p_iS_0^m\dP\no\\
&=0.
\end{align}
Thus, by (5.14) together with (5.13) and (5.15), it follows from the expressions of $\tf_{ij},\tf_{0i},\tf_0$
and assumption (\ref{apriori-bound}) that
\begin{align}
&|\int_{\Omega_T}I_1^m\cdot\mathcal{M}S_0^m\dP dtdx|\no\\
&\le C\ve^2+C\ve\bigg(R(T)^{\mu}\int_{S_T}(D_tS_0^m\dP)^2dS+R(T)^{\mu-3(\g-1)}\int_{S_T}(\nabla_xS_0^m\dP)^2dS\no\\
&\quad+\int_{\O_T}\big(R(t)^{\mu-1-\dl}(D_tS_0^m\dP)^2+R(t)^{\mu-1-3(\g-1)}(\nabla_xS_0^m\dP)^2dtdx\bigg),
\end{align}
where we have used
\begin{align}
&|a(t)R(t)^{\mu}(D_tS_0^m\dP)^2\ss_{i=1}^3\p_i\tf_{0i}|\no\\
&\le C\ve R(t)^{\mu-1-\f{3(\g-1)-\dl}{2}}(D_tS_0^m\dP)^2\no\\
&\le C\ve R(t)^{\mu-1-\dl}(D_tS_0^m\dP)^2\q \text{for}\q0<\dl\le\g-1.\no
\end{align}

\vskip 0.3cm
{\bf Part 2. Estimate on $\int_{\O_T}I_2^m\cdot\mathcal{M}S_0^m\dP dtdx$.}
\vskip 0.3cm

It follows from the expressions of $\tf_{ij},\tf_{0i},\tf_0$ and (\ref{apriori-bound}) that
$$|I_2^m|\le C\ve\bigg(R(t)^{-3(\g-1)-1+\f{\dl}2}\ss_{0\le l\le m-1}|\nabla_xS_0^{l+1}\dP|+R(t)^{-3(\g-1)}\ss_{0\le l\le m-1}|\nabla_x^2S_0^l\dP|\bigg).$$
This derives that for $m\le2$
\beq
\int_{\O_T}|I_2^m\cdot\mathcal{M}\dP|dtdx\le C\ve\int_{\Omega_T}\ss_{l=0}^m\bigg(R(t)^{\mu-1-\dl+2l}(\nabla_{t,x}^lD_t\dot\Phi)^2
+R(t)^{\mu-1-3(\g-1)+2l}(\nabla_{t,x}^l\nabla_x\dot\Phi)^2\bigg)dtdx,
\eeq
and for $m=3$
\ben
&&\int_{\O_T}|I_2^m\cdot\mathcal{M}\dP|dtdx\no\\
&\le& C\ve\int_{\Omega_T}\ss_{l=0}^2\bigg(R(t)^{\mu-1-\dl+2l}(\nabla_{t,x}^lD_t\dot\Phi)^2
+R(t)^{\mu-1-3(\g-1)+2l}(\nabla_{t,x}^l\nabla_x\dot\Phi)^2\bigg)dtdx\no\\
&&+C\ve\int_{\Omega_T}\bigg(R(t)^{\mu+5-\dl}(\nabla_{t,x}^3D_t\dot\Phi)^2
+R(t)^{\mu+5-3(\g-1)-\dl}(\nabla_{t,x}^3\nabla_x\dot\Phi)^2\bigg)dtdx.
\een

\vskip 0.3cm
{\bf Part 3. Estimate on $\int_{\O_T}I_3^m\cdot\mathcal{M}S_0^m\dP dtdx$.}
\vskip 0.3cm

First, we treat the case of $\int_{\O_T}I_3^m\cdot\mathcal{M}S_0^m\dP dx$ with $m\le2$.

For $m\le2$, as in Part 2, it follows from the expression of $f_i$ and assumption (\ref{apriori-bound}) that
$$|I_3^m|\le C\ve\bigg(R(t)^{-3(\g-1)-1+\f{\dl}2}\ss_{0\le l\le m}|\nabla_xS_0^l\dP|
+R(t)^{-3(\g-1)}\ss_{0\le l\le m}|\nabla_x^2S_0^l\dP|\bigg),$$
which yields for $m\le2$
\beq
\int_{\O_T}|I_3^m\cdot\mathcal{M}S^m\dP|dx\le C\ve\ss_{l=0}^m\int_{\O_T}\bigg(R(t)^{\mu-1-\dl+2l}(\nabla_{t,x}^lD_t\dot\Phi)^2
+R(t)^{\mu-1-3(\g-1)+2l}(\nabla_{t,x}^l\nabla_x\dot\Phi)^2\bigg)dtdx.
\eeq

Next we deal with $\int_{\O_T}I_3^3\cdot\mathcal{M}S_0^3\dP dtdx$.
Note that the difficult terms in $I_{3}^3$ are the those that
 include the products of third order derivatives
of $\dot\Phi$ because there are no related weighted $L^{\infty}$ estimates in (\ref{apriori-bound}). For the convenience, we decompose $I_{3}^3$ into $J_1$ and $J_2$
by using $S_0^2=LR(t)D_t+R(t)^2D_t^2$, where only $J_2$ contains the products of third order derivatives
of $\dot\Phi$. Namely,
\beq
I_{3}^3=J_1+J_2,
\eeq
where
\ben
J_1&=&\ss_{1\le l\le 2}C_{lm}\bigg\{\ss_{l_1+l_2=l,l_1\ge1}\tilde{C}_{l_1l_2}\bigg(\ss_{i=1}^3(S_0^{l_1}\tf_{0i})S_0^{l_2}(\p_iD_t\dP) +\ss_{1\le i,j\le3}(S_0^{l_1}\tf_{ij})S_0^{l_2}(\p_{ij}^2\dP)\bigg)\bigg\}\no\\
&&+C_{33}\ss_{(l_1,l_2)\neq(2,1)}\tilde{C}_{l_1l_2}\bigg(\ss_{i=1}^3(S_0^{l_1}\tf_{0i})S_0^{l_2}(\p_iD_t\dP) +\ss_{1\le i,j\le3}(S_0^{l_1}\tf_{ij})S_0^{l_2}(\p_{ij}^2\dP)\bigg)\no\\
&&+C_{33}\tilde{C}_{21}\bigg(\ss_{i=1}^3(LS_0\tf_{0i})S_0^{2}(\p_iD_t\dP) +\ss_{1\le i,j\le3}(LS_0\tf_{ij})S_0^{2}(\p_{ij}^2\dP)\bigg)+S_0^3\tf_0,
\een
and
$$J_2=C_{33}\tilde{C}_{21}\bigg(\ss_{i=1}^3(R(t)^2D_t^2\tf_{0i})S_0(\p_iD_t\dP)
+\ss_{1\le i,j\le3}(R(t)^2D_t^2\tf_{ij})S_0(\p_{ij}^2\dP)\bigg).$$

By assumption (\ref{apriori-bound}) and the expressions of $\tf_{0i},\tf_{ij},\tf_0$,  direct computation yields
\beq|J_1|\le C\ve\bigg(R(t)^{-3(\g-1)-1+\f{\dl}2}\ss_{0\le l\le 3}|\nabla_xS_0^l\dP|
+R(t)^{-3(\g-1)}\ss_{0\le l\le 2}|\nabla_x^2S_0^l\dP|\bigg).\eeq

On the other hand, by the expression of $f_{ij}$,
 we further decompose $J_2$  into $J_2=J_{21}+J_{22}$ so that only
$J_{22}$ contains the product terms of third order derivatives
of $\dot\Phi$. More precisely,
$$J_2=J_{21}+J_{22}$$
with
\ben
J_{22}&=&C_{33}\tilde{C}_{21}\bigg(-2R(t)^{-1}\ss_{i=1}^3(\p_iS_0^2\dP)^2
-4\ss_{i,j=1}^3\big(\p_i\dP\p_jS_0^2\dP\p_{ij}^2S_0\dP\big)-(\g-1)D_tS_0^2\dP\ss_{i=1}^3\p_i^2S_0\dP\no\\
&&-(\g-1)\ss_{i,j=1}^3\big(\p_i\dP\p_iS_0^2\dP\p_j^2S_0\dP\big)\bigg),\no
\een
and
\beq
|J_{21}|\le C\ve\bigg(R(t)^{-3(\g-1)-1+\f{\dl}2}\ss_{0\le l\le 3}|\nabla_xS_0^l\dP|
+R(t)^{-3(\g-1)}\ss_{0\le l\le 2}|\nabla_x^2S_0^l\dP|\bigg).
\eeq
Combining (5.22) and  (5.23)  yields
\ben
&&\int_{\O_T}|(J_1+J_{21})\cdot\mathcal{M}S_0^3\dP|dtdx\no\\
&\le& C\ve\ss_{l=0}^2\int_{\O_T}\bigg(R(t)^{\mu-1-\dl+2l}(\nabla_{t,x}^lD_t\dot\Phi)^2
+R(t)^{\mu-1-3(\g-1)+2l}(\nabla_{t,x}^l\nabla_x\dot\Phi)^2\bigg)dtdx\no\\
&&+\int_{\O_T}\bigg(R(t)^{\mu+5-\dl}(\nabla_{t,x}^lD_t\dot\Phi)^2
+R(t)^{\mu+5-3(\g-1)-\dl}(\nabla_{t,x}^l\nabla_x\dot\Phi)^2\bigg)dtdx.
\een

Finally, we estimate $\int_{\O_T}|J_{22}\cdot\mathcal{M}S_0^3\dP|dtdx$.
To overcome the difficulty induced by the lack of weighted $L^{\infty}$ estimates of
$|\na_{t,x}^3 \dot\Phi|$ in $J_{22}$, we will use the interpolation inequalities in Corollary 2.1 and Lemma 2.4.
In fact, by (5.1) and the expression of $J_{22}$, it is
sufficient to estimate some typical terms
in $\int_{D_T}|J_{22}\cdot\mathcal {M}S_0^3\dP|dtdx$ as follows:
\vskip 0.3 true cm
{\bf (A) Estimate on $|R(t)^{\mu-1}(\p_iS_0^2\dP)^2D_tS_0^3\dP|_{L^1(\O_T)}$.}

\begin{align}
&|R(t)^{\mu-1}(\p_iS_0^2\dP)^2D_tS_0^3\dP|_{L^1(\O_T)}\no\\
&\le |R(t)^{\f{2\mu+2-\dl}4}\p_iS_0^2\dP|_{L^4}^2|R(t)^{\f{\mu-1-\dl}2}D_tS_0^3\dP|_{L^2}\no\\
&\le C\ve\bigg(\ss_{k=0}^2|R(t)^{\f{\mu+3-\dl}2-k}\nabla_{t,x}^{3-k}S_0\dP|_{L^2}^{\f12}\bigg)^2 |R(t)^{\f{\mu-1-\dl}2}D_tS_0^3\dP|_{L^2}
\qquad \text{(Applying Lemma 2.4 (i) for $\dP$)}\no\\
&\le C\ve\int_{\Omega_T}\ss_{l=0}^2\bigg(R(t)^{\mu-1-\dl+2l}(\nabla_{t,x}^lD_t\dot\Phi)^2
+R(t)^{\mu-1-3(\g-1)+2l}(\nabla_{t,x}^l\nabla_x\dot\Phi)^2\bigg)dtdx\no\\
&\quad +C\ve\int_{\Omega_T}\bigg(R(t)^{\mu+5-\dl}(\nabla_{t,x}^3D_t\dot\Phi)^2+R(t)^{\mu+5-3(\g-1)-\dl}(\nabla_{t,x}^3\nabla_x\dot\Phi)^2\bigg)dtdx.
\end{align}

\vskip 0.3 true cm
{\bf (B) Estimate on $|R(t)^{\mu}\p_i\dP\p_jS_0^2\dP\p_{ij}S_0\dP D_tS_0^3\dP|_{L^1(\O_T)}$.}

\begin{align}
&|R(t)^{\mu}\p_i\dP\p_jS_0^2\dP\p_{ij}^2S_0\dP D_tS_0^3\dP|_{L^1(\O_T)}\no\\
&\le C\ve|R(t)^{\mu-3(\g-1)+\f{\dl}2}\p_jS_0^2\dP\p_{ij}^2S_0\dP D_tS_0^3\dP|_{L^1(\O_T)}\qquad\qquad
\text{(By assumption (5.1))}\no\\
&\le C\ve|R(t)^{\f{2\mu+2-\dl}4}\p_jS_0^2\dP|_{L^4}|R(t)^{\f{\mu+3-\dl}2}\p_{ij}^2S_0\dP|_{L^4}|R(t)^{\f{\mu-1-\dl}2}D_tS_0^3\dP|_{L^2}\no\\
&\le C\ve\bigg(\ss_{k=0}^2|R(t)^{\f{\mu+3-\dl}2-k}\nabla_{t,x}^{3-k}S_0\dP|_{L^2}^{\f12}\bigg)
\bigg(\ss_{k=0}^2|R(t)^{\f{\mu+3-\dl}2-k}\nabla_{t,x}^{3-k}S_0\dP|_{L^2}^{\f12}\bigg) |R(t)^{\f{\mu-1-\dl}2}D_tS_0^3\dP|_{L^2}\no\\
&\qquad\qquad \qquad \qquad \qquad \qquad \qquad\qquad \text{(Applying Lemma 2.4 (i)-(ii) for $\dP$)}\no\\
&\le C\ve\int_{\Omega_T}\ss_{l=0}^2\bigg(R(t)^{\mu-1-\dl+2l}(\nabla_{t,x}^lD_t\dot\Phi)^2
+R(t)^{\mu-1-3(\g-1)+2l}(\nabla_{t,x}^l\nabla_x\dot\Phi)^2\bigg)dtdx\no\\
&\quad+C\ve\int_{\Omega_T}\bigg(R(t)^{\mu+5-\dl}(\nabla_{t,x}^3D_t\dot\Phi)^2+R(t)^{\mu+5-3(\g-1)-\dl}(\nabla_{t,x}^3\nabla_x\dot\Phi)^2\bigg)dtdx.
\end{align}

\vskip 0.3 true cm
{\bf (C) Estimate on $|R(t)^{\mu}D_tS_0^2\dP\p_{i}^2S_0\dP D_tS_0^3\dP|_{L^1(\O_T)}$.}

\begin{align}
&|R(t)^{\mu}D_tS_0^2\dP\p_{i}^2S_0\dP D_tS_0^3\dP|_{L^1(\O_T)}\no\\
&\le C|R(t)^{\f{2\mu+2-\dl}4}D_tS_0^2\dP|_{L^4}|R(t)^{\f{\mu+3-\dl}2}\p_{i}^2S_0\dP|_{L^4}|R(t)^{\f{\mu-1-\dl}2}D_tS_0^3\dP|_{L^2}\no\\
&\le C\ve \bigg(\ss_{k=0}^2|R(t)^{\f{\mu+3-\dl}2-k}\nabla_{t,x}^{3-k}S_0\dP|_{L^2}^{\f12}\bigg)
\bigg(\ss_{k=0}^2|R(t)^{\f{\mu+3-\dl}2-k}\nabla_{t,x}^{3-k}S_0\dP|_{L^2}^{\f12}\bigg) |R(t)^{\f{\mu-1-\dl}2}D_tS_0^3\dP|_{L^2}\no\\
&\qquad\qquad \qquad \qquad \qquad \qquad \qquad\qquad\qquad \text{(Applying Lemma 2.4 (i)-(ii) for $\dP$)}\no\\
&\le C\ve\int_{\Omega_T}\ss_{l=0}^2\bigg(R(t)^{\mu-1-\dl+2l}(\nabla_{t,x}^lD_t\dot\Phi)^2
+R(t)^{\mu-1-3(\g-1)+2l}(\nabla_{t,x}^l\nabla_x\dot\Phi)^2\bigg)dtdx\no\\
&\quad+C\ve\int_{\Omega_T}\bigg(R(t)^{\mu+5-\dl}(\nabla_{t,x}^3D_t\dot\Phi)^2
+R(t)^{\mu+5-3(\g-1)-\dl}(\nabla_{t,x}^3\nabla_x\dot\Phi)^2\bigg)dtdx.
\end{align}
Combining (5.25)-(5.27) gives
\begin{align}
&\int_{\O_T}|J_{22}\cdot\mathcal{M}S_0^3\dP|dtdx\no\\
&\le C\ve\int_{\Omega_T}\ss_{l=0}^2\bigg(R(t)^{\mu-1-\dl+2l}(\nabla_{t,x}^lD_t\dot\Phi)^2
+R(t)^{\mu-1-3(\g-1)+2l}(\nabla_{t,x}^l\nabla_x\dot\Phi)^2\bigg)dtdx\no\\
&\quad +C\ve\int_{\Omega_T}\bigg(R(t)^{\mu+5-\dl}(\nabla_{t,x}^3D_t\dot\Phi)^2
+R(t)^{\mu+5-3(\g-1)-\dl}(\nabla_{t,x}^3\nabla_x\dot\Phi)^2\bigg)dtdx.
\end{align}
This, together with (5.19) and (5.24), yields for $m\le2$
\beq
\int_{\O_T}|I_3^m\cdot\mathcal{M}S_0^m\dP|dtdx\le C\ve\ss_{l=0}^m\bigg(|R(t)^{\f{\mu-1-\dl+2l}2}\nabla_{t,x}^lD_t\dot\Phi|_{L^2}^2
+|R(t)^{\f{\mu-1-3(\g-1)+2l}2}\nabla_{t,x}^l\nabla_x\dot\Phi|_{L^2}^2\bigg),
\eeq
and for $m=3$
\ben
&&\int_{\O_T}|I_3^m\cdot\mathcal{M}S_0^3\dP|dtdx\no\\
&\le&C\ve\int_{\Omega_T}\ss_{l=0}^2\bigg(R(t)^{\mu-1-\dl+2l}(\nabla_{t,x}^lD_t\dot\Phi)^2
+R(t)^{\mu-1-3(\g-1)+2l}(\nabla_{t,x}^l\nabla_x\dot\Phi)^2\bigg)dtdx\no\\
&&+C\ve\int_{\Omega_T}\bigg(R(t)^{\mu+5-\dl}(\nabla_{t,x}^3D_t\dot\Phi)^2
+R(t)^{\mu+5-3(\g-1)-\dl}(\nabla_{t,x}^3\nabla_x\dot\Phi)^2\bigg)dtdx.
\een

\vskip 0.3cm
{\bf Part 4. Estimate on $\int_{\O_T}B_{1m}\dP\cdot\mathcal{M}S_0^m\dP dtdx$.}
\vskip 0.3cm

First, from the expressions of $B_{1m}$ and $\dot{f}$, and by (\ref{apriori-bound}), we have
\begin{align}
|B_{11}\dP|&\le C\ve \bigg(R(t)^{-3(\g-1)+\f{\dl}2}|\nabla_x D_t\dP|+R(t)^{-3(\g-1)}|\nabla_x^2\dP|
+R(t)^{-3(\g-1)-1+\f{\dl}2}|\nabla_x\dP|\bigg),\\
|B_{1m}\dP|&\le C\ve \ss_{0\le l\le m-1}\bigg(R(t)^{-3(\g-1)+\f{\dl}2}|\nabla_xS_0^l D_t\dP|+R(t)^{-3(\g-1)}|\nabla_x^2S_0^l\dP|+R(t)^{-3(\g-1)-1+\f{\dl}2}|\nabla_xS_0^l\dP|\no\\
&\quad + R(t)^{-\f{3(\g-1)-\dl}2-1}|D_tS_0^l\dP|\bigg)+CR(t)^{-3(\g-1)}\ss_{0\le l\le m-2}|\Delta S_0^l\dP|,\q m=2,3.
\end{align}
Since (\ref{estimate-to-3}) holds for $l\le m-1$, we  have from (5.29)-(5.30) that, for $m\le2$,
\begin{align}
&\int_{\O_T}|B_{1m}\dP\cdot\mathcal{M}S_0^m\dP|dtdx\no\\
&\le C\ve\ss_{l=0}^m\int_{\O_T}\bigg(R(t)^{\mu-1-\dl+2l}(\nabla_{t,x}^lD_t\dot\Phi)^2
+R(t)^{\mu-1-3(\g-1)+2l}(\nabla_{t,x}^l\nabla_x\dot\Phi)^2\bigg)dtdx\no\\
&\quad +C\ve\bigg(\int_{\O_T}R(t)^{\mu-1-\dl}(D_tS_0^m\dP)^2dtdx\bigg)^{\f12},
\end{align}
and for $m=3$,
\begin{align}
&\int_{\O_T}|B_{1m}\dP\cdot\mathcal{M}S_0^m\dP|dtdx\no\\
&\le C\ve\int_{\Omega_T}\ss_{l=0}^2\bigg(R(t)^{\mu-1-\dl+2l}(\nabla_{t,x}^lD_t\dot\Phi)^2
+R(t)^{\mu-1-3(\g-1)+2l}(\nabla_{t,x}^l\nabla_x\dot\Phi)^2\bigg)dtdx\no\\
&\quad +C\ve\int_{\Omega_T}\bigg(R(t)^{\mu+5-\dl}(\nabla_{t,x}^3D_t\dot\Phi)^2
+R(t)^{\mu+5-3(\g-1)-\dl}(\nabla_{t,x}^3\nabla_x\dot\Phi)^2\bigg)dtdx.\no\\
&\quad +C\ve\bigg(\int_{\O_T}R(t)^{\mu-1-\dl}(D_tS_0^3\dP)^2dtdx\bigg)^{\f12}.
\end{align}

\vskip 0.3cm
{\bf Part 5. Estimate on $\int_{\O_T}B_{2m}\dP\cdot\mathcal{M}S_0^m\dP dtdx$.}
\vskip 0.3cm

Note that $B_{2m}\dot\Phi$ contains the $(m+1)-$th (the highest order) derivatives of $\dP$
and then $B_{2m}\dot\Phi\cdot\mathcal{M}S^m\Phi$ contains the term $\na_{t,x}^{\al}\dP\na_{t,x}^{\beta}\dP$
($|\al|=|\beta|=m+1$). Thanks to the definite sign
of the coefficient $-3m(\g-1)<0$ in the $B_{2m}$ and the  (\ref{linear-equ}),  $\na_{t,x}^{\al}\dP\na_{t,x}^{\beta}\dP$
with $|\al|=|\beta|=m+1$  can be controlled as  follows.

Since
$$
D_tS_0^m\dP=\ss_{0\le l\le m-1}C_{lm}D_tS_0^l\dP+R(t)S_0^{m-1}D_t^2\dP,
$$
it follows from (\ref{linear-equ}) that
\beq
D_tS_0^m\dP=\ss_{0\le l\le m-1}C_{lm}D_tS_0^l\dP+R(t)S_0^{m-1}\big(\hat{c}^2\Delta\dP-\f{3L(\g-1)}{R(t)}D_t\dP+\mathcal{L}\dP\big).
\eeq
Direct computation yields
\begin{align}
&\int_{\O_T}B_{2m}\dP\cdot\mathcal{M}S_0^m\dP dtdx\no\\
&=\int_{\O_T}-3m(\g-1)R(t)^{\mu}a(t)\bigg\{S_0^{m-1}(\hat{c}^2\Delta\dP)^2\no\\
&\q+\bigg(\ss_{0\le l\le m-1}C_{lm}D_tS_0^l\dP+R(t)S_0^{m-1}\big(\mathcal{L}\dP-\f{3L(\g-1)}{R(t)}D_t\dP\big)\bigg)S_0^{m-1}(\hat{c}^2\Delta\dP)\bigg\}dtdx\no\\
&\le-3m(\g-1)\int_{\O_T}R(t)^{\mu}a(t)\bigg(\ss_{0\le l\le m-1}C_{lm}D_tS_0^l\dP\no\\
&\q \qquad \qquad \qquad +R(t)S_0^{m-1}\big(\mathcal{L}\dP-\f{3L(\g-1)}{R(t)}D_t\dP\big)\bigg)S_0^{m-1}(\hat{c}^2\Delta\dP)dtdx.\no\\
\end{align}
Note that $\ss_{0\le l\le m-1}C_{lm}D_tS_0^l\dP-R(t)S_0^{m-1}(\f{3L(\g-1)}{R(t)}D_t\dP)$ only contains at most $m-$th
order derivatives of $\dP$, then we have by (\ref{apriori-bound}) that, for $l\le m-1$,
\beq
\int_{\O_T}R(t)^{\mu-1-\dl}\bigg(\ss_{0\le l\le m-1}C_{lm}D_tS_0^l\dP-R(t)S_0^{m-1}(\f{3L(\g-1)}{R(t)}D_t\dP)\bigg)^2dtdx\le C\ve^2.
\eeq
On the other hand, we have
\beq
|S_0^{m-1}\mathcal{L}\dP|\le C\ve\bigg(\ss_{0\le l\le m}R(t)^{-3(\g-1)-1+\f{\dl}2}|\nabla_xS_0^l\dP|+\ss_{0\le l\le m-1}R(t)^{-3(\g-1)}|\nabla_x^2S_0^l\dP|\bigg).
\eeq
Therefore, putting (5.37)-(5.38) into (5.36) yields for $m\le2$
\begin{align}
&\int_{\O_T}|B_{2m}\dP\cdot\mathcal{M}S_0^m\dP| dtdx\no\\
&\le C\ve^2+C\ve\ss_{l=0}^m\int_{\O_T}\bigg(R(t)^{\mu-1-\dl+2l}(\nabla_{t,x}^lD_t\dot\Phi)^2
+R(t)^{\mu-1-3(\g-1)+2l}(\nabla_{t,x}^l\nabla_x\dot\Phi)^2\bigg)dtdx,
\end{align}
and for $m=3$
\ben
&&\int_{\O_T}|B_{2m}\dP\cdot\mathcal{M}S_0^m\dP|dtdx\no\\
&\le&C\ve^2+ C\ve\int_{\Omega_T}\ss_{l=0}^2\bigg(R(t)^{\mu-1-\dl+2l}(\nabla_{t,x}^lD_t\dot\Phi)^2
+R(t)^{\mu-1-3(\g-1)+2l}(\nabla_{t,x}^l\nabla_x\dot\Phi)^2\bigg)dtdx\no\\
&&+C\ve\int_{\Omega_T}\bigg(R(t)^{\mu+5-\dl}(\nabla_{t,x}^3D_t\dot\Phi)^2
+R(t)^{\mu+5-3(\g-1)-\dl}(\nabla_{t,x}^3\nabla_x\dot\Phi)^2\bigg)dtdx.
\een

Consequently, putting (5.17)-(5.18), (5.29), (5.32), (5.37)-(5.38) into (5.8), we  complete the proofs of (5.4) and (5.6)
for $S^m=S_0^m$.

For the case of $m=0$, (5.5) comes directly from Theorem 4.1, (5.17)-(5.18) and (5.29).

\vskip 0.3cm
{\bf Case (2)\q $S^m=Z^m$.}
\vskip 0.3cm
By direct computation, we have
\beq
Z\mathcal{L}=\mathcal{L}Z,
\eeq
that implies
\beq
\mathcal{L}Z^m=Z^m\mathcal{L}.
\eeq
In addition, from equation (\ref{L-in-Dt}), we have that for $0\le m\le3$,
\beq
Z^m\mathcal{L}\dP=K_1^m+K_2^m+K_3^m,
\eeq
where
\ben
&&K_1^m=\ss_{i=1}^3\tf_{0i}\p_iD_tZ^m\dP+\ss_{1\le i\neq j\le3}\tf_{ij}^2\p_{ij}Z^m\dP+\ss_{i=1}^3\tf_{ii}\p_{i}^2Z^m\dP,\no\\
&&K_2^m=\ss_{i=1}^3\tf_{0i}[Z^m,\p_iD_t]\dP+\ss_{0\le i\neq j\le3}\tf_{ij}[Z^m,\p_{ij}^2]\dP+\ss_{i=1}^3\tf_{ii}[Z^m,\p_{i}^2]\dP,\no\\
&&K_3^m=\ss_{1\le l\le m}C_{lm}\bigg\{\ss_{l_1+l_2=l,l_1\ge1}\tilde{C}_{l_1l_2}\bigg(\ss_{i=1}^3(Z^{l_1}\tf_{0i})Z^{l_2}(\p_iD_t\dP) +\ss_{1\le i,j\le3}(Z^{l_1}\tf_{ij})Z^{l_2}(\p_{ij}^2\dP)\bigg)\bigg\}+Z^m\tf_0.\no
\een

Based on the above preparation, we now estimate $\int_{\O_T}\mathcal{L}Z^m\dP\cdot\mathcal{M}Z^m\dP dtdx$ on the right hand side of (\ref{estimate-S_0^m}) in the following three steps.

\vskip 0.3cm
{\bf Part 1. Estimate on $\int_{\O_T}K_1^m\cdot\mathcal{M}S^m\dP dtdx.$}
\vskip 0.3cm

Replacing $S_0$ in (5.16) by $Z$,  we have as in (5.16),
\begin{align}
&|\int_{\Omega_T}K_1^m\cdot\mathcal{M}Z^m\dP dtdx|\no\\
&\le C\ve^2+C\ve\bigg(R(T)^{\mu}\int_{S_T}(D_tZ^m\dP)^2dS+R(T)^{\mu-3(\g-1)}\int_{S_T}(\nabla_xZ^m\dP)^2dS\no\\
&\quad +\int_{\O_T}\big(R(t)^{\mu-1-\dl}(D_tZ^m\dP)^2+R(t)^{\mu-1-3(\g-1)}(D_tZ^m\dP)^2dtdx\bigg).
\end{align}

\vskip 0.3cm
{\bf Part 2. Estimate on $\int_{\O_T}K_2^m\cdot\mathcal{M}Z^m\dP dtdx.$}
\vskip 0.3cm

Direct computation yields
$$
|K_2^m|\le C \bigg(\ss_{i=1}^3|\tf_{0i}|\ss_{0\le l\le m-1}|\nabla_xZ^lD_t\dP|+\ss_{i,j=1}^3|\tf_{ij}|\ss_{0\le l\le m-2}|\nabla_x^2Z^{l}\dP|\bigg).
$$
Then it follows from the expressions of $\tf_{ij},\tf_{0i}$ and (\ref{apriori-bound}) that
$$
|K_2^m|\le C\ve \bigg(R(t)^{-3(\g-1)+\f{\dl}2}\ss_{0\le l\le m-1}|\nabla_xZ^lD_t\dP|+R(t)^{-3(\g-1)}\ss_{0\le l\le m-2}|\nabla_x^2Z^{l}\dP|\bigg),
$$
which gives that for $m\le2$,
\begin{align}
&\int_{\O_T}|K_2^m\cdot\mathcal{M}Z^m\dP|dtdx\no\\
&\le C\ve\ss_{l=0}^m\int_{\O_T}\bigg(R(t)^{\mu-1-\dl+2l}(\nabla_{t,x}^lD_t\dot\Phi)^2+R(t)^{\mu-1-3(\g-1)+2l}(\nabla_{t,x}^l\nabla_x\dot\Phi)^2\bigg)dtdx,
\end{align}
and for $m=3$,
\begin{align}
&\int_{\O_T}|K_2^m\cdot\mathcal{M}Z^m\dP|dtdx\no\\
&\le C\ve\int_{\Omega_T}\ss_{l=0}^2\bigg(R(t)^{\mu-1-\dl+2l}(\nabla_{t,x}^lD_t\dot\Phi)^2
+R(t)^{\mu-1-3(\g-1)+2l}(\nabla_{t,x}^l\nabla_x\dot\Phi)^2\bigg)dtdx\no\\
&\quad +C\ve\int_{\Omega_T}\bigg(R(t)^{\mu+5-\dl}(\nabla_{t,x}^3D_t\dot\Phi)^2
+R(t)^{\mu+5-3(\g-1)-\dl}(\nabla_{t,x}^3\nabla_x\dot\Phi)^2\bigg)dtdx.
\end{align}

\vskip 0.3cm
{\bf Part 3. Estimate on $\int_{\O_T}K_3^m\cdot\mathcal{M}Z^m\dP dtdx.$}
\vskip 0.3cm

For $m=1$, it follows from the expressions of $\tf_{ij},\tf_{0i},\tf_0$ and (\ref{apriori-bound}) that
$$|K_3^1|\le C\ve\bigg(R(t)^{-3(\g-1)-1+\f{\dl}2}(|\nabla_x\dP|+|\nabla_xZ\dP|)+R(t)^{-3(\g-1)}|\nabla_x^2\dP|\bigg).$$
This yields
\beq
\int_{\O_T}|K_3^1\cdot\mathcal{M}\dP|dtdx\le C\ve\ss_{l=0}^1\int_{\O_T}\bigg(R(t)^{\mu-1-\dl+2l}(\nabla_{t,x}^lD_t\dot\Phi)^2+R(t)^{\mu-1-3(\g-1)+2l}(\nabla_{t,x}^l\nabla_x\dot\Phi)^2\bigg)dtdx.
\eeq

For $m=2$,
$$K_3^2=K_{31}^2+K_{32}^2,$$
where
\begin{align}
K_{31}^2&=\ss_{1\le l\le 2}C_{l2}\bigg\{\ss_{l_1+l_2=l,l_1\ge1}\tilde{C}_{l_1l_2}\bigg(\ss_{i=1}^3(Z^{l_1}\tf_{0i})Z^{l_2}(\p_iD_t\dP) +\ss_{1\le i,j\le3}(Z^{l_1}\tf_{ij})Z^{l_2}(\p_{ij}\dP)\bigg)\bigg\}\no\\
&\quad -\f{3(\g-5)L}{2R(t)}\ss_{i=1}^3\p_i\dP Z^2\p_i\dP-\f{3(\g-5)}{R(t)}\ss_{i=1}^3([Z,\p_i]\dP)^2,\no\\
K_{32}^2&=-\f{3(\g-5)}{R(t)}\ss_{i=1}^3(\p_iZ\dP)^2.\no
\end{align}
Therefore, it follows from the expressions of $\tf_{ij},\tf_{0i},\tf_0$ and (\ref{apriori-bound}) that
\ben
&&|K_{31}^2|\le C\ve\bigg(R(t)^{-3(\g-1)-1+\dl}\ss_{0\le l\le2}|\nabla_xZ^l\dP|+R(t)^{-\f{3(\g-1)}2-1}\ss_{0\le l\le2}|D_tZ^l\dP|\no\\
&&\q\q\q\q\q+R(t)^{-\f{3(\g-1)-\dl}2}\ss_{i=1}^3|Z\p_iD_t\dP|+R(t)^{-3(\g-1)+\f{\dl}2}|Z\nabla_x^2\dP|\bigg).\no
\een
This implies
\begin{align}
&\int_{\O_T}|K_{31}^2\cdot\mathcal{M}Z^2\dP|dtdx\no\\
&\le C\ve\ss_{l=0}^2\int_{\O_T}\bigg(R(t)^{\mu-1-\dl+2l}(\nabla_{t,x}^lD_t\dot\Phi)^2
+R(t)^{\mu-1-3(\g-1)+2l}(\nabla_{t,x}^l\nabla_x\dot\Phi)^2\bigg)dtdx.
\end{align}

Next, we deal with $\int_{\O_T}|K_{32}^2\cdot\mathcal{M}Z^2\dP|dtdx$.
By H\"older inequality and Lemma 2.4, we have  for $0<\dl\le\g-1$ that
\begin{align}
&|R(t)^{\mu-1}(\nabla_xZ\dP)^2D_tZ^2\dP|_{L^1(\O_T)}\no\\
&\le C|R(t)^{\f{2\mu+5-3\g-\dl}4}\nabla_xZ\dP|_{L^4(\O_t)}^2|R(t)^{\f{\mu-1-\dl}2}D_tZ\dP|_{L^2(\O_T)}\no\\
&\le C\ve\ss_{l=0}^2\bigg(|R(t)^{\f{\mu-1-\dl+2l}2}\nabla_{t,x}^lD_t\dot\Phi|_{L^2}
+|R(t)^{\f{\mu-1-3(\g-1)+2l}2}\nabla_{t,x}^l\nabla_x\dot\Phi|_{L^2}\bigg) |R(t)^{\f{\mu-1-\dl}2}D_tZ^2\dP|_{L^2}\no\\
&\qquad\qquad \qquad \qquad \qquad \qquad\qquad\qquad \qquad\qquad \text{(Applying Lemma 2.4 (iii) for $\dP$)}\no\\
&\le C\ve\ss_{l=0}^2\bigg(|R(t)^{\f{\mu-1-\dl+2l}2}\nabla_{t,x}^lD_t\dot\Phi|_{L^2}^2
+|R(t)^{\f{\mu-1-3(\g-1)+2l}2}\nabla_{t,x}^l\nabla_x\dot\Phi|_{L^2}^2\bigg),
\end{align}
which yields
\beq
\int_{\O_T}|K_{32}^2\cdot\mathcal{M}Z^2\dP|dtdx\le C\ve\ss_{l=0}^2\int_{\O_T}\bigg(R(t)^{\mu-1-\dl+2l}(\nabla_{t,x}^lD_t\dot\Phi)^2+R(t)^{\mu-1-3(\g-1)+2l}(\nabla_{t,x}^l\nabla_x\dot\Phi)^2\bigg)dtdx.
\eeq
Combining (5.47) and (5.49) gives
\beq
\int_{\O_T}|K_{3}^2\cdot\mathcal{M}Z^2\dP|dtdx\le C\ve\ss_{l=0}^2\int_{\O_T}\bigg(R(t)^{\mu-1-\dl+2l}(\nabla_{t,x}^lD_t\dot\Phi)^2+R(t)^{\mu-1-3(\g-1)+2l}(\nabla_{t,x}^l\nabla_x\dot\Phi)^2\bigg)dtdx.
\eeq

Next, we consider the case of $m=3$.
Due to the lack of $L^{\infty}$ assumptions on the third order derivatives of $\dot\Phi$
in (\ref{apriori-bound}), we will decompose $K_3^3$ as follows
$$K_3^3=K_{31}^3+K_{32}^3,$$
where $K_{31}^3$ is linear with respect to the third and fourth order derivatives of $\dot\Phi$,
$K_{32}^3$ contains the products of two third order
derivatives of $\dot\Phi$.  $K_{31}^3$ and $K_{32}^3$ have the following
expressions:
\begin{align}
K_{31}^3&=\ss_{1\le l\le 2}C_{l2}\bigg\{\ss_{l_1+l_2=l,l_1\ge1}\tilde{C}_{l_1l_2}\bigg(\ss_{i=1}^3(Z^{l_1}\tf_{0i})Z^{l_2}(\p_iD_t\dP) +\ss_{1\le i,j\le3}(Z^{l_1}\tf_{ij})Z^{l_2}(\p_{ij}^2\dP)\bigg)\bigg\}\no\\
&+C_{33}\ss_{(l_1,l_2)=(1,2)~\text{or}~(3,0)}\tilde{C}_{l_1l_2}\bigg(\ss_{i=1}^3(Z^{l_1}\tf_{0i})Z^{l_2}(\p_iD_t\dP) +\ss_{1\le i,j\le3}(Z^{l_1}\tf_{ij})Z^{l_2}(\p_{ij}^2\dP)\bigg)\no\\
&+C_{33}C_{21}\bigg(\ss_{i=1}^3Z^2\tf_{0i}[Z,\p_i]D_t\dP+\ss_{i,j=1}^3\big(Z^2\tf_{ij}[Z,\p_{ij}^2]\dP-2Z\p_i\dP Z\p_j\dP\p_{ij}^2Z\dP\big)\no\\
&-2\ss_{i=1}^3[Z^2,\p_i]\dP\p_iZD_t\dP-(\g-1)\ss_{i,j=1}^3\big(\p_j\dP[Z^2,\p_j]\dP\p_{i}^2Z\dP+(Z\p_j\dP)^2\p_{i}^2Z\dP\big)\bigg)\no\\
&
-\f{6(\g-5)L}{R(t)}\ss_{i=1}^3[Z,\p_i]\dP Z^2\p_i\dP,\no
\end{align}
and
\begin{align}
K_{32}^3&=-C_{33}\tilde{C}_{21}\bigg(-2\ss_{i=1}^3\p_iZ^2\dP \p_iZD_t\dP-2\ss_{1\le i,j\le3}\p_iZ^2\dP\p_j\dP\p_{ij}^2Z\dP\no\\
&\quad-(\g-1)\big(D_tZ^2\dP+\ss_{j=1}^3\p_j\dP\p_jZ^2\dP\big)\ss_{i=1}^3\p_{i}^2Z\dP\bigg)-\f{6(\g-5)L}{R(t)}\ss_{i=1}^3\p_iZ\dP Z^2\p_i\dP.\no
\end{align}
Then it follows from the expressions of $\tf_{ij},\tf_{0i},\tf_0$ and (\ref{apriori-bound}) that
\ben
&&|K_{31}^3|\le C\ve\bigg(R(t)^{-3(\g-1)-1+\dl}\ss_{0\le l\le3}|\nabla_xZ^l\dP|+R(t)^{-\f{3(\g-1)+\dl}2-1}\ss_{0\le l\le3}|D_tZ^l\dP|\no\\
&&\q\q\q\q\q+R(t)^{-\f{3(\g-1)-\dl}2}\ss_{0\le l\le2}\ss_{i=1}^3|Z^l\p_iD_t\dP|+R(t)^{-3(\g-1)+\dl}\ss_{0\le l\le2}|Z^l\nabla_x^2\dP|\bigg),\no
\een
which gives
\begin{align}
&\int_{\O_T}|K_{31}^3\cdot\mathcal{M}Z^2\dP|dtdx\no\\
&\le C\ve\int_{\Omega_T}\ss_{l=0}^2\bigg(R(t)^{\mu-1-\dl+2l}(\nabla_{t,x}^lD_t\dot\Phi)^2
+R(t)^{\mu-1-3(\g-1)+2l}(\nabla_{t,x}^l\nabla_x\dot\Phi)^2\bigg)dtdx\no\\
&\quad+C\ve\int_{\Omega_T}\bigg(R(t)^{\mu+5-\dl}(\nabla_{t,x}^3D_t\dot\Phi)^2
+R(t)^{\mu+5-3(\g-1)-\dl}(\nabla_{t,x}^3\nabla_x\dot\Phi)^2\bigg)dtdx.
\end{align}
We now turn to estimate $\int_{\O_T}|K_{32}^3\cdot\mathcal{M}S_0^3\dP|dtdx$.
Again, to overcome the  lack of weighted $L^{\infty}$ estimates of
$|\na_{t,x}^3 \dot\Phi|$ in $K_{32}^3$, we use the interpolation inequalities in Corollary 2.1 and Lemma 2.4.
In fact, by (5.1) and the expression of $K_{32}^3$, it suffices to
estimate the following typical terms
in $\int_{D_T}|K_{32}^3\cdot\mathcal {M}Z^3\dP|dtdx$:
\vskip 0.3 true cm
{\bf (A) Estimate on $|R(t)^{\mu}\p_iZ^2\dP \p_iZD_t\dP D_tZ^3\dP|_{L^1(\O_T)}$.}

\begin{align}
&|R(t)^{\mu}\p_iZ^2\dP \p_iZD_t\dP D_tZ^3\dP|_{L^1(\O_T)}\no\\
&\le C|R(t)^{\f{3\g-5-\dl}2}\p_iZ^2\dP|_{L^4}|R(t)^{\f{6\g-6-\dl}2}\p_iZD_t\dP|_{L^4}|R(t)^{\f{\mu-1-\dl}2}D_tZ^3\dP|_{L^2}\no\\
&\le C\ve\bigg(\ss_{k=0}^1|R(t)^{\f{\mu+4-3\g}2-k}\nabla_{t,x}^{2-k}Z\dP|_{L^2}^{\f12}
+\ss_{k=0}^2|R(t)^{\f{3\g-3-\dl}2-k}\nabla_{t,x}^{3-k}Z\dP|_{L^2}^{\f12}\bigg) |R(t)^{\f{\mu-1-\dl}2}D_tZ^3\dP|_{L^2}\no\\
&\qquad\qquad \qquad \qquad \qquad \qquad \qquad\qquad \text{(Applying Lemma 2.4 (v)-(vi) for $\dP$)}\no\\
&\le C\ve\int_{\Omega_T}\ss_{l=0}^2\bigg(R(t)^{\mu-1-\dl+2l}(\nabla_{t,x}^lD_t\dot\Phi)^2+R(t)^{\mu-1-3(\g-1)
+2l}(\nabla_{t,x}^l\nabla_x\dot\Phi)^2\bigg)dtdx\no\\
&\quad+C\ve\int_{\Omega_T}\bigg(R(t)^{\mu+5-\dl}(\nabla_{t,x}^3D_t\dot\Phi)^2+R(t)^{\mu
+5-3(\g-1)-\dl}(\nabla_{t,x}^3\nabla_x\dot\Phi)^2\bigg)dtdx.
\end{align}

\vskip 0.3 true cm
{\bf (B) Estimate on $|R(t)^{\mu}\p_iZ^2\dP\p_j\dP \p_{ij}^2Z\dP D_tZ^3\dP|_{L^1(\O_T)}$.}
\begin{align}
&|R(t)^{\mu}\p_iZ^2\dP\p_j\dP \p_{ij}^2Z\dP D_tZ^3\dP|_{L^1(\O_T)}\no\\
&\le C\ve|R(t)^{\f{3\g-5-\dl}2}\p_iZ^2\dP|_{L^4}|R(t)^{\f{3\g-3-\dl}2}\p_{ij}^2Z\dP|_{L^4}|R(t)^{\f{\mu-1-\dl}2}D_tZ^3\dP|_{L^2}
\qquad\text{(By assumption (5.1))}\no\\
&\le C\ve\bigg(\ss_{k=0}^1|R(t)^{\f{\mu+4-3\g}2-k}\nabla_{t,x}^{2-k}Z\dP|_{L^2}^{\f12}
+\ss_{k=0}^2|R(t)^{\f{3\g-3-\dl}2-k}\nabla_{t,x}^{3-k}Z\dP|_{L^2}^{\f12}\bigg) |R(t)^{\f{\mu-1-\dl}2}D_tZ^3\dP|_{L^2}\no\\
&\qquad\qquad \qquad \qquad \qquad \qquad \qquad\qquad \text{(Applying Lemma 2.4 (iv) and (vii) for $\dP$)}\no\\
&\le C\ve\int_{\Omega_T}\ss_{l=0}^2\bigg(R(t)^{\mu-1-\dl+2l}(\nabla_{t,x}^lD_t\dot\Phi)^2
+R(t)^{\mu-1-3(\g-1)+2l}(\nabla_{t,x}^l\nabla_x\dot\Phi)^2\bigg)dtdx\no\\
&\quad+C\ve\int_{\Omega_T}\bigg(R(t)^{\mu+5-\dl}(\nabla_{t,x}^3D_t\dot\Phi)^2
+R(t)^{\mu+5-3(\g-1)-\dl}(\nabla_{t,x}^3\nabla_x\dot\Phi)^2\bigg)dtdx.
\end{align}

\vskip 0.3 true cm
{\bf (C) Estimate on $|R(t)^{\mu}D_tZ^2\dP\p_{i}^2Z\dP D_tZ^3\dP|_{L^1(\O_T)}$.}

\begin{align}
&|R(t)^{\mu}D_tZ^2\dP\p_{i}^2Z\dP D_tZ^3\dP|_{L^1(\O_T)}\no\\
&\le
C|R(t)^{\f{6\g-8-\dl}2}D_tZ^2\dP|_{L^4}|R(t)^{\f{3\g-3-\dl}2}\p_{i}^2Z\dP|_{L^4}|R(t)^{\f{\mu-1-\dl}2}D_tZ^3\dP|_{L^2}\no\\
&\le
C\ve\bigg(\ss_{k=0}^2|R(t)^{\f{3\g-3-\dl}2-k}\nabla_{t,x}^{3-k}Z\dP|_{L^2}^{\f12}
+\ss_{k=0}^2|R(t)^{\f{6\g-6-\dl}2-k}\nabla_{t,x}^{2-k}ZD_t\dP|_{L^2}^{\f12}\bigg) |R(t)^{\f{\mu-1-\dl}2}D_tZ^3\dP|_{L^2}\no\\
&\qquad\qquad \qquad \qquad \qquad \qquad \qquad\qquad \text{(Applying Lemma 2.4 (iv) and (vii) for $\dP$)}\no\\
&\le
C\ve\int_{\Omega_T}\ss_{l=0}^2\bigg(R(t)^{\mu-1-\dl+2l}(\nabla_{t,x}^lD_t\dot\Phi)^2
+R(t)^{\mu-1-3(\g-1)+2l}(\nabla_{t,x}^l\nabla_x\dot\Phi)^2\bigg)dtdx\no\\
&\quad +C\ve\int_{\Omega_T}\bigg(R(t)^{\mu+5-\dl}(\nabla_{t,x}^3D_t\dot\Phi)^2
+R(t)^{\mu+5-3(\g-1)-\dl}(\nabla_{t,x}^3\nabla_x\dot\Phi)^2\bigg)dtdx.
\end{align}

\vskip 0.3 true cm
{\bf (D) Estimate on $|R(t)^{\mu-1}\p_iZ\dP Z^2\p_i\dP D_tZ^3\dP|_{L^1(\O_T)}$.}

\begin{align}
&|R(t)^{\mu-1}\p_iZ\dP Z^2\p_i\dP D_tZ^3\dP|_{L^1(\O_T)}\no\\
&=|R(t)^{\f{5\dl-3(\g-1)}4}\cdot R(t)^{\f{9\g-13-\dl}4}\p_iZ\dP\cdot R(t)^{\f{3\g-5-\dl}2}Z^2\p_{i}\dP \cdot R(t)^{\f{\mu-1-\dl}2}D_tZ^3\dP|_{L^1(\O_T)}\no\\
&\qquad \qquad\qquad \qquad \qquad \qquad  \text{(here we have used the restriction of $\dl\le\f{3}{5}(\g-1)$)}\no\\
&\le
C|R(t)^{\f{9\g-13-\dl}4}\p_iZ\dP|_{L^4}|R(t)^{\f{3\g-5-\dl}2}Z^2\p_{i}\dP|_{L^4}|R(t)^{\f{\mu-1-\dl}2}D_tZ^3\dP|_{L^2}\no\\
&\le
C\ve\bigg(\ss_{k=0}^1|R(t)^{\f{\mu+4-3\g}2-k}\nabla_{t,x}^{2-k}Z\dP|_{L^2}^{\f12}
+\ss_{k=0}^2|R(t)^{\f{3\g-3-\dl}2-k}\nabla_{t,x}^{3-k}Z\dP|_{L^2}^{\f12}\bigg) |R(t)^{\f{\mu-1-\dl}2}D_tZ^3\dP|_{L^2}\no\\
&\qquad \qquad\qquad \qquad \qquad \qquad\text{(Applying Lemma 2.4 (iii) and (vi) for $\dP$)}\no\\
&\le
C\ve\int_{\Omega_T}\ss_{l=0}^2\bigg(R(t)^{\mu-1-\dl+2l}(\nabla_{t,x}^lD_t\dot\Phi)^2
+R(t)^{\mu-1-3(\g-1)+2l}(\nabla_{t,x}^l\nabla_x\dot\Phi)^2\bigg)dtdx\no\\
&\quad+C\ve\int_{\Omega_T}\bigg(R(t)^{\mu+5-\dl}(\nabla_{t,x}^3D_t\dot\Phi)^2
+R(t)^{\mu+5-3(\g-1)-\dl}(\nabla_{t,x}^3\nabla_x\dot\Phi)^2\bigg)dtdx.
\end{align}

Substituting (5.53)-(5.56) into $\int_{\O_T}|K_{32}^3\cdot\mathcal{M}Z^3\dP|dtdx$ yields
\begin{align}
&\int_{\O_T}|K_{32}^3\cdot\mathcal{M}Z^3\dP|dtdx\no\\
&\le
C\ve\int_{\Omega_T}\ss_{l=0}^2\bigg(R(t)^{\mu-1-\dl+2l}(\nabla_{t,x}^lD_t\dot\Phi)^2
+R(t)^{\mu-1-3(\g-1)+2l}(\nabla_{t,x}^l\nabla_x\dot\Phi)^2\bigg)dtdx\no\\
&\quad +C\ve\int_{\Omega_T}\bigg(R(t)^{\mu+5-\dl}(\nabla_{t,x}^3D_t\dot\Phi)^2
+R(t)^{\mu+5-3(\g-1)-\dl}(\nabla_{t,x}^3\nabla_x\dot\Phi)^2\bigg)dtdx.
\end{align}

Then putting (5.44)-(5.46), (5.51)-(5.52) and (5.57) into (5.8) yields (5.4) and (5.6) for the case $S^m=Z^m$.

\vskip 0.3cm
{\bf Case (3)\quad $S^m=S_0^{l_1}Z^{l_2}$ $(1\le l_1, l_2\le m-1, l_1+l_2=m)$.}

Since the estimation for this case is similar to
the Case (1) and Case (2), we omit the detail.

Combining the Cases (1)-(3), we complete the proof of Lemma 5.1. $\hfill\square$

Based on Lemma 5.1, we now derive a series of estimates on the higher order
derivatives of $\dP$. For convenience, set $S_1=r\p_r=\ss_{i=1}^3x_i\p_i$.

\begin{lemma}{\bf (The second order radial derivative estimates)} Under the assumptions of Theorem 5.1, we have
\begin{align}
&R(T)^{\mu}\int_{S_T}(D_tS_1\dP)^2dS+R(T)^{\mu-3(\g-1)}\int_{S_T}(\nabla_xS_1\dP)^2dS\no\\
&\quad+\int_{\Omega_T}\bigg(R(t)^{\mu-1-\dl}(D_tS_1\dP)^2+R(t)^{\mu-1-3(\g-1)}(\nabla_xS_1\dP)^2\bigg)dtdx\no\\
&\le C\ve^2+C\ve\int_{\Omega_T}\ss_{l=0}^1\bigg(R(t)^{\mu-1-\dl+2l}(\nabla_{t,x}^lD_t\dot\Phi)^2
+R(t)^{\mu-1-3(\g-1)+2l}(\nabla_{t,x}^l\nabla_x\dot\Phi)^2\bigg)dtdx\no\\
&\quad +C\ve\bigg(\int_{\Omega_T}\ss_{l=0}^1\bigg(R(t)^{\mu-1-\dl+2l}(\nabla_{t,x}^lD_t\dot\Phi)^2
+R(t)^{\mu-1-3(\g-1)+2l}(\nabla_{t,x}^l\nabla_x\dot\Phi)^2\bigg)dtdx\bigg)^{\f12}.
\end{align}
\end{lemma}

\begin{remark} Lemma 5.2, together with Lemma 5.1 for $m=1$ yield
\ben
&&\ss_{0\le l_1+l_2\le1}\bigg(R(T)^{\mu}\int_{S_T}(D_tS^{l_1}S_1^{l_2}\dP)^2dS+R(T)^{\mu-3(\g-1)}\int_{S_T}(\nabla_xS^{l_1}S_1^{l_2}\dP)^2dS\bigg)\no\\
&&+\ss_{0\le l_1+l_2\le1}\int_{\Omega_T}\bigg(R(t)^{\mu-1-\dl}(D_tS^{l_1}S_1^{l_2}\dP)^2+R(t)^{\mu-1-3(\g-1)}(\nabla_xS^{l_1}S_1^{l_2}\dP)^2\bigg)dtdx\no\\
&\le&C\ve^2+C\ve\int_{\Omega_T}\ss_{l=0}^1\bigg(R(t)^{\mu-1-\dl+2l}(\nabla_{t,x}^lD_t\dot\Phi)^2
+R(t)^{\mu-1-3(\g-1)+2l}(\nabla_{t,x}^l\nabla_x\dot\Phi)^2\bigg)dtdx\no\\
&&+C\ve\bigg(\int_{\Omega_T}\ss_{l=0}^1\bigg(R(t)^{\mu-1-\dl+2l}(\nabla_{t,x}^lD_t\dot\Phi)^2
+R(t)^{\mu-1-3(\g-1)+2l}(\nabla_{t,x}^l\nabla_x\dot\Phi)^2\bigg)dtdx\bigg)^{\f12}.\no
\een
\end{remark}

{\bf Proof of Lemma 5.2.} Noting $\mathcal{L}S_1\dP=S_1\mathcal{L}\dP-2\hat{c}^2\Delta\dP$, then it follows from Theorem 4.1 that
\ben
&&R(T)^{\mu}\int_{S_T}(D_tS_1\dP)^2dS+R(T)^{\mu-3(\g-1)}\int_{S_T}(\nabla_xS_1\dP)^2dS\no\\
&&+\int_{\Omega_T}\bigg(R(t)^{\mu-1-\dl}(D_tS_1\dP)^2+R(t)^{\mu-1-3(\g-1)}(\nabla_xS_1\dP)^2\bigg)dtdx\no\\
&\le&C\ve^2+C\int_{\O_T}\big(S_1\dot{f}-2\hat{c}^2\Delta\dP\big)\cdot\mathcal{M}S_1\dP dtdx.
\een

Since $S_1$ and $Z$ are expressed in terms of  $x_i\p_j$ $(1\le i,j\le 3)$,  similar
to the  proofs of (5.43)-(5.45) and (5.47), we have
\begin{align}
&|\int_{\Omega_T}S_1\dot{f}\cdot\mathcal{M}S_1\dP dtdx|\no\\
&\le C\ve^2+C\ve\bigg(R(T)^{\mu}\int_{S_T}(D_tS_1\dP)^2dS+R(T)^{\mu-3(\g-1)}\int_{S_T}(\nabla_xS_1\dP)^2dS\bigg)\no\\
&\quad +C\ve\int_{\Omega_T}\ss_{l=0}^1\bigg(R(t)^{\mu-1-\dl+2l}(\nabla_{t,x}^lD_t\dot\Phi)^2
+R(t)^{\mu-1-3(\g-1)+2l}(\nabla_{t,x}^l\nabla_x\dot\Phi)^2\bigg)dtdx.
\end{align}

In addition, thanks to (5.4) for the case of $S^m=S_0$, it follows
\begin{align}
&\int_{\O_T}|\hat{c}^2\Delta\dP\cdot\mathcal{M}S_1\dP| dtdx\no\\
&\le C\int_{\O_T}R(t)^{\mu-3(\g-1)}|\Delta\dP \nabla_xS_0\dP| dtdx\no\\
&\le C\big(\int_{\O_T}R(t)^{\mu+1-3(\g-1)}(\Delta\dP)^2dtdx\big)^{\f12}
\big(\int_{\O_T}R(t)^{\mu-1-3(\g-1)}(\nabla_xS_0\dP)^2dtdx\big)^{\f12}\no\\
&\le C\ve^2+ C\ve\bigg(\ss_{l=0}^1\int_{\Omega_T}R(t)^{\mu-1-\dl+2l}(\nabla_{t,x}^lD_t\dot\Phi)^2
+R(t)^{\mu-1-3(\g-1)+2l}(\nabla_{t,x}^l\nabla_x\dot\Phi)^2 dtdx\bigg)^{\f12}.
\end{align}
Finally, substituting (5.60)-(5.61) into (5.59) yields (5.58),
and this completes the proof of the lemma. $\hfill\square$

\begin{lemma}{\bf (The third order radial derivative estimates)} Under the assumptions of Theorem 5.1,
we have
\ben
&&R(T)^{\mu}\int_{S_T}(D_tS_1^2\dP)^2dS+R(T)^{\mu-3(\g-1)}\int_{S_T}(\nabla_xS_1^2\dP)^2dS\no\\
&&+\int_{\Omega_T}\bigg(R(t)^{\mu-1-\dl}(D_tS_1^2\dP)^2+R(t)^{\mu-1-3(\g-1)}(\nabla_xS_1^2\dP)^2\bigg)dtdx\no\\
&\le&C\ve^2+C\ve\int_{\Omega_T}\bigg(\ss_{l=0}^3 R(t)^{\mu-1-\dl+2l}(\nabla_{t,x}^lD_t\dot\Phi)^2
+\ss_{l=0}^2 R(t)^{\mu-1-3(\g-1)+2l}(\nabla_{t,x}^l\nabla_x\dot\Phi)^2\bigg)dtdx\no\\
&&+C\ve\bigg(\int_{\Omega_T}R(t)^{\mu+3-\dl}(\nabla_{t,x}^2D_t\dot\Phi)^2
+R(t)^{\mu+3-3(\g-1)}(\nabla_{t,x}^2\nabla_x\dot\Phi)^2 dtdx\bigg)^{\f12}.
\een
\end{lemma}

\begin{remark}
Under the assumptions of Theorem 5.1, as in the proof of Lemma 5.1, we have
\ben
&&R(T)^{\mu}\int_{S_T}(D_tSS_1\dP)^2dS+R(T)^{\mu-3(\g-1)}\int_{S_T}(\nabla_xSS_1\dP)^2dS\no\\
&&\quad +\int_{\Omega_T}\bigg(R(t)^{\mu-1-\dl}(D_tSS_1\dP)^2+R(t)^{\mu-1-3(\g-1)}(\nabla_xSS_1\dP)^2\bigg)dtdx\no\\
&\le&C\ve^2+C\ve\int_{\Omega_T}\ss_{l=0}^2\bigg(R(t)^{\mu-1-\dl+2l}(\nabla_{t,x}^lD_t\dot\Phi)^2
+R(t)^{\mu-1-3(\g-1)+2l}(\nabla_{t,x}^l\nabla_x\dot\Phi)^2\bigg)dtdx\no\\
&&+C\ve\bigg(\ss_{l=0}^2\int_{\Omega_T}R(t)^{\mu-1-\dl+2l}(\nabla_{t,x}^lD_t\dot\Phi)^2
+R(t)^{\mu-1-3(\g-1)+2l}(\nabla_{t,x}^l\nabla_x\dot\Phi)^2 dtdx\bigg)^{\f12}.
\een

This, together with (5.62) and Lemma 5.1 for $m=2$ and Remark 2.2, yields the following estimate
\ben
&&\ss_{0\le l_1+l_2\le2}\bigg(R(T)^{\mu}\int_{S_T}(D_tS^{l_1}S_1^{l_2}\dP)^2dS+R(T)^{\mu-3(\g-1)}\int_{S_T}(\nabla_xS^{l_1}S_1^{l_2}\dP)^2dS\bigg)\no\\
&&\quad +\ss_{0\le l_1+l_2\le2}\int_{\Omega_T}\bigg(R(t)^{\mu-1-\dl}(D_tS^{l_1}S_1^{l_2}\dP)^2+R(t)^{\mu-1-3(\g-1)}(\nabla_xS^{l_1}S_1^{l_2}\dP)^2\bigg)dtdx\no\\
&\le&C\ve^2+C\ve\int_{\Omega_T}\ss_{l=0}^2\bigg(R(t)^{\mu-1-\dl+2l}(\nabla_{t,x}^lD_t\dot\Phi)^2
+R(t)^{\mu-1-3(\g-1)+2l}(\nabla_{t,x}^l\nabla_x\dot\Phi)^2\bigg)dtdx\no\\
&&+C\ve\int_{\Omega_T}\bigg(R(t)^{\mu+5-\dl}(\nabla_{t,x}^3D_t\dot\Phi)^2
+R(t)^{\mu+5-3(\g-1)-\dl}(\nabla_{t,x}^3\nabla_x\dot\Phi)^2\bigg)dtdx\no\\
&&+C\ve\bigg(\int_{\Omega_T}\ss_{l=0}^2\bigg(R(t)^{\mu-1-\dl+2l}(\nabla_{t,x}^lD_t\dot\Phi)^2
+R(t)^{\mu-1-3(\g-1)+2l}(\nabla_{t,x}^l\nabla_x\dot\Phi)^2\bigg)dtdx\no\\
&&+\int_{\Omega_T}\bigg(R(t)^{\mu+5-\dl}(\nabla_{t,x}^3D_t\dot\Phi)^2
+R(t)^{\mu+5-3(\g-1)-\dl}(\nabla_{t,x}^3\nabla_x\dot\Phi)^2\bigg)dtdx\bigg)^{\f12}.
\een

\end{remark}

{\bf Proof of Lemma 5.3.} In order to derive the third order radial derivative estimate, we need to study a higher order boundary condition
on $B_T$.
Differentiating (\ref{L-in-Z}) with respect to $r$ and applying (\ref{tangent-bd-condition}) yield on
$B_T$
\beq
S_1^3\dP-S_1^2\dP-2\ss_{i=1}^3Z_i^2\dP-\f{G}{A}=0,
\eeq
where
\begin{align}
A&=\hat{c}^2-(\g-1)D_t\dP-\f{\g-1}{2R(t)^2}\ss_{i=1}^3(Z_i\dP)^2,\no\\
G&=-\f{\g-1}{R(t)^2}\ss_{i=1}^3(Z_i\dP)^2\bigg(S_1^2\dP+\ss_{i=1}^3Z_i^2\dP\bigg)
-\f{3(\g-1)L}{R(t)}\ss_{i=1}^3(Z_i\dP)^2-\f{2L}{R(t)}\ss_{i,j=1}^3C_{ij}Z_i\dP Z_j\dP\no\\
&\quad -\f{4}{R(t)^2}\ss_{i,j=1}^3Z_i\dP Z_j\dP Z_iZ_j\dP+\f{1}{R(t)^2}S_1^2\dP\ss_{i,j=1}^3C_{ij}Z_i\dP Z_j\dP-\f{4}{R(t)^2}\ss_{i,j,k=1}^3C_{ijk}Z_i\dP Z_j\dP Z_k\dP.\no
\end{align}
For convenience, we rewrite (5.65) as follows
\beq
S_1^3\dP-S_1^2\dP-2\ss_{i=1}^3Z_i^2\dP-\chi\big(\f{r}{R(t)}\big)\f{G}{A}=0\qquad \text{on $B_T$},
\eeq
where $\chi(s)$ is a smooth cut-off function
\beq\chi(s)=\bec\q\q\q1, &\text{for}\q \f23\le s\le1,\\
\q\q\q0,&\text{for}\q 0\le s\le \f13,\\
\text{smooth connection}, &\text{for}\q\f13\le s\le \f23.
\eec\eeq
On the other hand, direct computation yields
\beq
\mathcal{L} (S_1^2-S_1)\dP=(S_1^2-S_1)\mathcal{L}\dP+4\hat{c}^2\Delta S_1\dP-2\hat{c}^2\Delta\dP.
\eeq
As in (4.2), we have
\begin{align}
&\int_{\Omega_T}\mathcal{L}(S_1^2-S_1)\dP\cdot \mathcal{M}(S_1^2-S_1)\dP dtdx\no\\
&=-\int_{B_T}\f{R(t)^{\mu}}{\sqrt{1+L^2}}\hat{c}^2a(t)D_t(S_1^2-S_1)\dP\cdot\p_r(S_1^2-S_1)\dP dS\no\\
&\quad+\bigg(\int_{S_T}-\int_{S^0}\bigg)\f12R(t)^{\mu}a(t)\big(D_t(S_1^2-S_1)\dP+\hat{c}^2|\nabla_x(S_1^2-S_1)\dP|^2\big)dS\no\\
&\quad+\int_{\O_T}\biggl\{\f{\dl}2R(t)^{\mu-1-\dl}\big(D_t(S_1^2-S_1)\dP\big)^2+\f{\g}2R(t)^{\mu-1-3(\g-1)}\big(L(5-3\g)a(t)\no\\
&\quad -R(t)a'(t)\big)|\nabla_x(S_1^2-S_1)\dP|^2\biggr\}dtdx.\no\\
\end{align}
Substituting (5.68) into (5.69), it follows
\begin{align}
&R(T)^{\mu}\int_{S_T}(D_t(S_1^2-S_1)\dP)^2dS+R(T)^{\mu-3(\g-1)}\int_{S_T}(\nabla_x(S_1^2-S_1)\dP)^2dS\no\\
&\quad+\int_{\Omega_T}\bigg(R(t)^{\mu-1-\dl}(D_t(S_1^2-S_1)\dP)^2+R(t)^{\mu-1-3(\g-1)}(\nabla_x(S_1^2-S_1)\dP)^2\bigg)dtdx\no\\
&\le
C\ve^2+\int_{\O_T}\bigg((S_1^2-S_1)\dot f+4\hat{c}^2\Delta S_1\dP-2\hat{c}^2\Delta\dP\bigg)\cdot\mathcal{M}(S_1^2-S_1)\dP dtdx\no\\
&\quad+\int_{B_T}\f{R(t)^{\mu}}{\sqrt{1+L^2}}\hat{c}^2a(t)D_t(S_1^2-S_1)\dP\cdot\p_r(S_1^2-S_1)\dP dS.
\end{align}
In view of (5.66) and divergence theorem, we have
\begin{align}
&\int_{B_T}\f{R(t)^{\mu}}{\sqrt{1+L^2}}\hat{c}^2a(t)D_t(S_1^2-S_1)\dP\cdot\p_r(S_1^2-S_1)\dP dS\no\\
&=\int_{\O_T}R(t)^{\mu}\hat{c}^2a(t)\ss_{i=1}^3\p_i\bigg[\p_i(D_tS_1\dP-D_t\dP)\bigg(\chi\big(\f{r}{R(t)}\big)\f{G}{A}
+2\ss_{i=1}^3Z_i^2\dP\bigg)\bigg]dtdx\no\\
&=\int_{\O_T}R(t)^{\mu}\hat{c}^2a(t)\bigg[\Delta(D_tS_1\dP-D_t\dP)\bigg(\chi\big(\f{r}{R(t)}\big)\f{G}{A}
+2\ss_{i=1}^3Z_i^2\dP\bigg)\bigg]dtdx\no\\
&\q+\int_{\O_T}R(t)^{\mu}\hat{c}^2a(t)\ss_{j=1}^3\p_j\big(D_tS_1\dP-D_t\dP\big)\bigg(\p_j\big(\chi\big(\f{r}{R(t)}\big)\f{G}{A}\big)
+2\ss_{i=1}^3\p_jZ_i^2\dP\bigg)dtdx.
\end{align}

And it follows from (\ref{apriori-better-bound}) that
\begin{align}
&|\chi\big(\f{r}{R(t)}\big)\f{G}{A}|\le C\ve\bigg(R(t)^{-3(\g-1)}\big(|S_1^2\dP|+|Z^2\dP|\big)+|Z\dP|\bigg),\\
&|\nabla_x\bigg(\chi\big(\f{r}{R(t)}\big)\f{G}{A}\bigg)|\le C\ve\bigg(R(t)^{-\f{3(\g-1)}2-1}\big(|S_1^2\dP|+|Z^2\dP|\big)
+R(t)^{-3(\g-1)}\big(|\nabla_xS_1^2\dP|+|\nabla_xZ^2\dP|\big)\no\\
&\qquad\qquad\qquad \qquad \qquad +\ss_{0\le l\le1}R(t)^{l-1}|\nabla_x^lZ\dP|+|\nabla_xS_0\dP|\bigg).
\end{align}

Combining (5.71)-(5.73), it follows
\begin{align}
&\bigg|\int_{B_T}\f{R(t)^{\mu}}{\sqrt{1+L^2}}\hat{c}^2a(t)D_t(S_1^2-S_1)\dP\cdot\p_r(S_1^2-S_1)\dP dS\bigg|\no\\
&\le C\ve^2+C\ve\int_{\Omega_T}\bigg(\ss_{l=0}^3 R(t)^{\mu-1-\dl+2l}(\nabla_{t,x}^lD_t\dot\Phi)^2
+\ss_{l=0}^2 R(t)^{\mu-1-3(\g-1)+2l}(\nabla_{t,x}^l\nabla_x\dot\Phi)^2\bigg)dtdx.
\end{align}
Then as in (5.61), thanks to (5.63) and (\ref{estimate-S_0^m}) in the case of $S^2=ZS_0$, it follows
\ben
&&\int_{\O_T}\big(|\hat{c}^2\Delta S_1\dP|+|\hat{c}^2\Delta\dP|\big)\cdot|\mathcal{M}\big(S_1^2+S_1\big)\dP| dtdx\no\\
&\le&C\ve^2+C\ve\bigg(\int_{\Omega_T}R(t)^{\mu+3-\dl}(\nabla_{t,x}^2D_t\dot\Phi)^2
+R(t)^{\mu+3-3(\g-1)}(\nabla_{t,x}^2\nabla_x\dot\Phi)^2 dtdx\bigg)^{\f12}.
\een
Next we deal with $\int_{\O_T}(S_1^2-S_1)\dot f\cdot\mathcal{M}(S_1^2-S_1)\dP dtdx$.
Similar to (5.12), we have
\beq
(S_1^2-S_1)\dot f=I_1+I_2+I_3,
\eeq
where
\begin{align}
&I_1=\ss_{i=1}^3\tf_{0i}\p_iD_t(S_1^2-S_1)\dP+\ss_{1\le i\neq j\le3}\tf_{ij}\p_{ij}^2(S_1^2-S_1)\dP
+\ss_{i=1}^3\tf_{ii}\p_{i}^2(S_1^2-S_1)\dP,\no\\
&I_2=\ss_{i=1}^3\tf_{0i}[(S_1^2-S_1),\p_iD_t]\dP+\ss_{0\le i\neq j\le3}\tf_{ij}[(S_1^2-S_1),\p_{ij}^2]\dP+\ss_{i=1}^3\tf_{ii}[(S_1^2-S_1),\p_{i}^2]\dP,\no\\
&I_3=\ss_{1\le l\le 2}C_{l1}\bigg\{\ss_{l_1+l_2=l,l_1\ge1}\tilde{C}_{l_1l_2}\bigg(\ss_{i=1}^3(S_1^{l_1}\tf_{0i})S_1^{l_2}(\p_iD_t\dP)
 +\ss_{1\le i,j\le3}(S_1^{l_1}\tf_{ij})S_1^{l_2}(\p_{ij}^2\dP)\bigg)\bigg\}+(S_1^2-S_1)\tf_0.\no
\end{align}
Similar to  (5.14), one has
\begin{align}
&I_1\cdot\mathcal{M}(S_1^2-S_1)\dP\no\\
&=D_t\bigg(-\f12a(t)R(t)^{\mu}\ss_{i,j=1}^3\big(\tf_{ij}\p_i(S_1^2-S_1)\dP\p_j(S_1^2-S_1)\dP\big)\bigg)\no\\
&\quad+\ss_{i=1}^3\p_i\bigg(\f12a(t)R(t)^{\mu}\tf_{0i}(D_t(S_1^2-S_1)\dP)^2
+a(t)R(t)^{\mu}D_t(S_1^2-S_1)\dP\ss_{j=1}^3\tf_{ij}\p_j(S_1^2-S_1)\dP\bigg)\no\\
&\quad-\f12a(t)R(t)^{\mu}(D_t(S_1^2-S_1)\dP)^2\ss_{i=1}^3\p_if_{0i}
-a(t)LR(t)^{\mu-1}\ss_{i,j=1}^3\big(\tf_{ij}\p_i(S_1^2-S_1)\dP\p_j(S_1^2-S_1)\dP\big)\no\\
&\quad-a(t)R(t)^{\mu}D_t(S_1^2-S_1)\dP\ss_{i,j=1}^3\big(\p_i\tf_{ij}\p_j(S_1^2-S_1)\dP\big)\no\\
&\quad+\f12\ss_{i,j=1}^3\big(D_t(a(t)R(t)^{\mu}\tf_{ij})\p_i(S_1^2-S_1)\dP\p_j(S_1^2-S_1)\dP\big).
\end{align}
On the other hand, by the expressions of $\tf_{ij},\tf_{0i}$ and (\ref{tangent-bd-condition}),  direct observation yields on $B_T$
\begin{align}
&\ss_{i=1}^3x_i\cdot\bigg(\f12a(t)R(t)^{\mu}\tf_{0i}(D_t(S_1^2-S_1)\dP)^2
+a(t)R(t)^{\mu}D_t(S_1^2-S_1)\dP\ss_{j=1}^3\tf_{ij}\p_j(S_1^2-S_1)\dP\bigg)\no\\
&=-a(t)R(t)^{\mu}(D_t(S_1^2-S_1)\dP)^2\ss_{i=1}^3\big(x_i\p_i\dP\big)
-a(t)R(t)^{\mu}D_t(S_1^2-S_1)\dP\ss_{i,j=1}^3\big(x_i\p_i\dP\p_j\dP\big)\no\\
&\quad -(\g-1)a(t)R(t)^{\mu}D_t(S_1^2-S_1)\dP\big(D_t\dP+\f12\ss_{k=1}^3(\p_k\dP)^2\big)\ss_{i=1}^3x_i\p_i(S_1^2-S_1)\dP\no\\
&=-(\g-1)a(t)R(t)^{\mu+1}\big(D_t\dP+\f12\ss_{k=1}^3(\p_k\dP)^2\big)\p_r(D_t S_1-D_t)\dP\cdot \bigg(\chi\big(\f{r}{R(t)}\big)\f{G}{A}+2\ss_{i=1}^3Z_i^2\dP\bigg).
\end{align}

Similar to (5.71)-(5.74), we obtain
\begin{align}
&\bigg|\int_{B_T}\ss_{i=1}^3x_i\cdot\bigg(\f12a(t)R(t)^{\mu}\tf_{0i}(D_t(S_1^2-S_1)\dP)^2
+a(t)R(t)^{\mu}D_t(S_1^2-S_1)\dP\ss_{j=1}^3\tf_{ij}\p_j(S_1^2-S_1)\dP\bigg)dS\bigg|\no\\
&\le
C\ve^2+C\ve\int_{\Omega_T}\bigg(\ss_{l=0}^3 R(t)^{\mu-1-\dl+2l}(\nabla_{t,x}^lD_t\dot\Phi)^2
+\ss_{l=0}^2 R(t)^{\mu-1-3(\g-1)+2l}(\nabla_{t,x}^l\nabla_x\dot\Phi)^2\bigg)dtdx.
\end{align}

Thus, as in (5.44), it follows
\begin{align}
&|\int_{\Omega_T}I_1\cdot\mathcal{M}(S_1^2-S_1)\dP dtdx|\no\\
&\le C\ve^2+C\ve\bigg(R(T)^{\mu}\int_{S_T}(D_t(S_1^2-S_1)\dP)^2dS+R(T)^{\mu-3(\g-1)}\int_{S_T}(\nabla_x(S_1^2-S_1)\dP)^2dS\no\\
&\quad+\int_{\Omega_T}\bigg(\ss_{l=0}^3 R(t)^{\mu-1-\dl+2l}(\nabla_{t,x}^lD_t\dot\Phi)^2
+\ss_{l=0}^2 R(t)^{\mu-1-3(\g-1)+2l}(\nabla_{t,x}^l\nabla_x\dot\Phi)^2\bigg)dtdx\bigg).
\end{align}

Then, similar to  (5.45) and (5.51), it follows
\begin{align}
&\int_{\Omega_T}\big(|I_2|+|I_3|\big)\cdot|\mathcal{M}(S_1^2-S_1)\dP| dtdx\no\\
&\le C\ve^2+C\ve\ss_{l=0}^2\int_{\Omega_T}\bigg( R(t)^{\mu-1-\dl+2l}(\nabla_{t,x}^lD_t\dot\Phi)^2
+ R(t)^{\mu-1-3(\g-1)+2l}(\nabla_{t,x}^l\nabla_x\dot\Phi)^2\bigg)dtdx.
\end{align}

Substituting (5.74)-(5.76) and (5.80)-(5.81) into (5.70), we
complete the proof of Lemma 5.3. $\hfill\square$

\begin{remark}
Under the assumptions of Theorem 5.1, as the estimation in the
proof of  Lemma 5.2, we have
\begin{align}
&R(T)^{\mu}\int_{S_T}(D_tS^2S_1\dP)^2dS+R(T)^{\mu-3(\g-1)}\int_{S_T}(\nabla_xS^2S_1\dP)^2dS\no\\
&\quad +\int_{\Omega_T}\bigg(R(t)^{\mu-1-\dl}(D_tS^2S_1\dP)^2+R(t)^{\mu-1-3(\g-1)}(\nabla_xS^2S_1\dP)^2\bigg)dtdx\no\\
&\le
C\ve^2+C\ve\int_{\Omega_T}\ss_{l=0}^2\bigg(R(t)^{\mu-1-\dl+2l}(\nabla_{t,x}^lD_t\dot\Phi)^2
+R(t)^{\mu-1-3(\g-1)+2l}(\nabla_{t,x}^l\nabla_x\dot\Phi)^2\bigg)dtdx\no\\
&\quad+C\ve\int_{\Omega_T}\bigg(R(t)^{\mu+5-\dl}(\nabla_{t,x}^3D_t\dot\Phi)^2
+R(t)^{\mu+5-3(\g-1)-\dl}(\nabla_{t,x}^3\nabla_x\dot\Phi)^2\bigg)dtdx\no\\
&\quad+C\ve\bigg(\int_{\Omega_T}\bigg(R(t)^{\mu+5-\dl}(\nabla_{t,x}^3D_t\dot\Phi)^2
+R(t)^{\mu+5-3(\g-1)-\dl}(\nabla_{t,x}^3\nabla_x\dot\Phi)^2\bigg)dtdx\bigg)^{\f12}.
\end{align}
\end{remark}

Next, we will give the estimates on $\nabla_{t,x}^4\dP$ in Lemma 5.4-Lemma 5.6.

\begin{lemma} {\bf(Estimate on $D_tSS_1^2\dP$ and $\nabla_xSS_1^2\dP$)} Under the assumptions of Theorem 5.1,
we have
\ben
&&R(T)^{\mu}\int_{S_T}(D_tSS_1^2\dP)^2dS+R(T)^{\mu-3(\g-1)}\int_{S_T}(\nabla_xSS_1^2\dP)^2dS\no\\
&&+\int_{\Omega_T}\bigg(R(t)^{\mu-1-\dl}(D_tSS_1^2\dP)^2+R(t)^{\mu-1-3(\g-1)}(\nabla_xSS_1^2\dP)^2\bigg)dtdx\no\\
&\le&C\ve^2+C\ve\ss_{l=0}^2\int_{\Omega_T}\bigg( R(t)^{\mu-1-\dl+2l}(\nabla_{t,x}^lD_t\dot\Phi)^2
+ R(t)^{\mu-1-3(\g-1)+2l}(\nabla_{t,x}^l\nabla_x\dot\Phi)^2\bigg)dtdx\no\\
&&+C\ve\int_{\Omega_T}\bigg( R(t)^{\mu+5-\dl}(\nabla_{t,x}^3D_t\dot\Phi)^2
+ R(t)^{\mu+5-3(\g-1)-\dl}(\nabla_{t,x}^3\nabla_x\dot\Phi)^2\bigg)dtdx\no\\
&&+C\ve\bigg(\ss_{l=0}^2\int_{\Omega_T}\bigg( R(t)^{\mu-1-\dl+2l}(\nabla_{t,x}^lD_t\dot\Phi)^2
+ R(t)^{\mu-1-3(\g-1)+2l}(\nabla_{t,x}^l\nabla_x\dot\Phi)^2\bigg)dtdx\no\\
&&+\int_{\Omega_T}\bigg( R(t)^{\mu+5-\dl}(\nabla_{t,x}^3D_t\dot\Phi)^2
+ R(t)^{\mu+5-3(\g-1)-\dl}(\nabla_{t,x}^3\nabla_x\dot\Phi)^2\bigg)dtdx\bigg)^{\f12}.
\een
\end{lemma}

{\bf Proof.} Applying $S$ to (5.66) yields
\beq
S_1(SS_1^2\dP-SS_1\dP)-2\ss_{i=1}^3SZ_i^2\dP-\chi\big(\f{r}{R(t)}\big)S\big(\f{G}{A}\big)=0\qquad
\text{on $B_T$}.
\eeq
In addition,  direct computation yields
\beq
\mathcal{L} Z(S_1^2-S_1)\dP=Z(S_1^2-S_1)\mathcal{L}\dP+4\hat{c}^2\Delta ZS_1\dP-2\hat{c}^2\Delta Z\dP,
\eeq
and
\ben
\mathcal{L} S_0(S_1^2-S_1)\dP&=&S_0(S_1^2-S_1)\mathcal{L}\dP-2L(S_1^2-S_1)\mathcal{L}\dP+S_0\bigg(4\hat{c}^2\Delta S_1\dP-2\hat{c}^2\Delta \dP\bigg)\no\\
&&-8L\hat{c}^2\Delta S_1\dP+4L\hat{c}^2\Delta\dP+3(\g-1)\hat{c}^2\Delta(S_1^2-S_1)\dP.
\een
Similar to (5.69), we obtain
\ben
&&\int_{\Omega_T}\mathcal{L}S(S_1^2-S_1)\dP\cdot \mathcal{M}S(S_1^2-S_1)\dP dtdx\no\\
&=&-\int_{B_T}\f{R(t)^{\mu}}{\sqrt{1+L^2}}\hat{c}^2a(t)D_tS(S_1^2-S_1)\dP\cdot\p_rS(S_1^2-S_1)\dP dS\no\\
&&+\bigg(\int_{S_T}-\int_{S^0}\bigg)\f12R(t)^{\mu}a(t)\big(D_tS(S_1^2-S_1)\dP+\hat{c}^2|\nabla_xS(S_1^2-S_1)\dP|^2\big)dS\no\\
&&+\int_{\O_T}\biggl\{\f{\dl}2R(t)^{\mu-1-\dl}\big(D_tS(S_1^2-S_1)\dP\big)^2+\f{\g}2R(t)^{\mu-1-3(\g-1)}\big(L(5-3\g)a(t)\no\\
&&\quad -R(t)a'(t)\big)|\nabla_xS(S_1^2-S_1)\dP|^2\biggr\}dtdx.
\een
As in (5.71), we have
\begin{align}
&\int_{B_T}\f{R(t)^{\mu}}{\sqrt{1+L^2}}\hat{c}^2a(t)D_tS(S_1^2-S_1)\dP\cdot\p_rS(S_1^2-S_1)\dP dS\no\\
&=\int_{\O_T}R(t)^{\mu}\hat{c}^2a(t)\bigg[\Delta(D_tSS_1\dP-D_t S\dP)\bigg(\chi\big(\f{r}{R(t)}\big)S\big(\f{G}{A}\big)+2\ss_{i=1}^3SZ_i^2\dP\bigg)\bigg]dtdx\no\\
&\quad+\int_{\O_T}R(t)^{\mu}\hat{c}^2a(t)\ss_{j=1}^3\p_j\big(D_t SS_1\dP-D_t S\dP\big)\bigg(\p_j\big(\chi\big(\f{r}{R(t)}\big)S\f{G}{A}\big)+2\ss_{i=1}^3\p_jSZ_i^2\dP\bigg)dtdx.
\end{align}
Note that, for $C^1$-smooth functions $f$ and $g$,
\ben
&&\int_{\O_T}(x_i\p_jf-x_j\p_if)gdtdx=\int_{\O_T}\bigg(\p_j(x_ifg)-\p_i(x_jfg)\bigg)-f(x_i\p_jg-x_j\p_ig)dtdx\no\\
&=&-\int_{\Omega_T}f(x_i\p_jg-x_j\p_ig)dtdx.\no
\een
Then
\beq
\int_{\O_T}Zf\cdot gdtdx=-\int_{\O_T}f\cdot Zgdtdx.
\eeq
Hence,  if we choose $S=Z$, then it follows from (5.88) that
\ben
&&\int_{\O_T}R(t)^{\mu}\hat{c}^2a(t)\bigg[\Delta(D_tZS_1\dP-ZD_t\dP)
\bigg(\chi\big(\f{r}{R(t)}\big)Z\big(\f{G}{A}\big)+2\ss_{i=1}^3ZZ_i^2\dP\bigg)\bigg]dtdx\no\\
&=&-\int_{\O_T}R(t)^{\mu}\hat{c}^2a(t)\bigg[\Delta(D_tS_1\dP-D_t\dP)
\bigg(\chi\big(\f{r}{R(t)}\big)Z^2\big(\f{G}{A}\big)+2\ss_{i=1}^3Z^2Z_i^2\dP\bigg)\bigg]dtdx.
\een
In addition, it follows from (\ref{apriori-better-bound}) that
\begin{align}
&\big|\chi\big(\f{r}{R(t)}\big)Z\big(\f{G}{A}\big)\big|
\le C\ve \bigg(R(t)^{-\f{3(\g-1)}2}|S_1^2\dP|+R(t)^{-3(\g-1)}(|ZS_1^2\dP|+|Z^3\dP|)+|Z\dP|+|Z^2\dP|\bigg),\\
&\big|\chi\big(\f{r}{R(t)}\big)Z^2\big(\f{G}{A}\big)\big|\le C\ve\bigg(|S_1^2\dP|+\ss_{l=1}^4|Z^l\dP|+R(t)^{-\f{3(\g-1)}2}|ZS_1^2\dP|+R(t)^{-3(\g-1)}|Z^2S_1^2\dP|\no\\
&\q\q\q\q\q\q\q\qquad+ R(t)|D_tZ^2\dP|\bigg)+R(t)^{3(\g-1)-1}|Z^2\dP|^2.
\end{align}
Combining (5.88) and (5.90)-(5.92), we have  that from (5.88)
\ben
&&|\int_{B_T}\f{R(t)^{\mu}}{\sqrt{1+L^2}}\hat{c}^2a(t)D_tZ(S_1^2-S_1)\dP\cdot\p_rZ(S_1^2-S_1)\dP dS|\no\\
&\le&C\ve\ss_{l=0}^2\int_{\Omega_T}\bigg( R(t)^{\mu-1-\dl+2l}(\nabla_{t,x}^lD_t\dot\Phi)^2
+ R(t)^{\mu-1-3(\g-1)+2l}(\nabla_{t,x}^l\nabla_x\dot\Phi)^2\bigg)dtdx\no\\
&&+C\ve\int_{\Omega_T}\bigg( R(t)^{\mu+5-\dl}(\nabla_{t,x}^3D_t\dot\Phi)^2
+ R(t)^{\mu+5-3(\g-1)-\dl}(\nabla_{t,x}^3\nabla_x\dot\Phi)^2\bigg)dtdx\no\\
&&+\int_{\O_T}R(t)^{\mu-1}|Z^2\dP|^2\cdot|\Delta(D_tS_1\dP-D_t\dP)|dtdx.
\een
By H\"older inequality and Lemma 2.4, we have
\begin{align}
&|R(t)^{\mu-1}(Z^2\dP)^2 \Delta(D_tS_1\dP-D_t\dP)|_{L^1(\O_T)}\no\\
&\le C|R(t)^{\mu-1}(\nabla_xZ\dP)^2\ss_{l=2}^3\big((R(t)\nabla_x)^lD_t\dP\big)|_{L^1(\O_T)}\no\\
&\le
C|R(t)^{\f{2\mu+5-3\g-\dl}4}\nabla_xZ\dP|_{L^4(\O_T)}^2|R(t)^{\f{\mu-1-\dl}2}
\ss_{l=2}^3\big((R(t)\nabla_x)^lD_t\dP\big)|_{L^2(\O_T)}\no\\
&\le C\ve\ss_{l=0}^2\bigg(|R(t)^{\f{\mu-1-\dl+2l}2}\nabla_{t,x}^lD_t\dot\Phi|_{L^2}
+|R(t)^{\f{\mu-1-3(\g-1)+2l}2}\nabla_{t,x}^l\nabla_x\dot\Phi|_{L^2}\bigg)
 |R(t)^{\f{\mu-1-\dl}2}\ss_{l=2}^3\big((R(t)\nabla_x)^lD_t\dP\big)|_{L^2}\no\\
&\qquad\qquad \qquad \qquad \qquad \qquad \qquad \qquad \qquad \qquad \qquad \qquad \qquad
\qquad  \text{(Applying Lemma 2.4 (iii) for $\dP$)}\no\\
&\le C\ve\bigg(\ss_{l=0}^3|R(t)^{\f{\mu-1-\dl+2l}2}\nabla_{t,x}^lD_t\dot\Phi|_{L^2}^2
+\ss_{l=0}^2
|R(t)^{\f{\mu-1-3(\g-1)+2l}2}\nabla_{t,x}^l\nabla_x\dot\Phi|_{L^2}^2\bigg),\no
\end{align}
that gives
\begin{align}
&\int_{\O_T}R(t)^{\mu-1}|Z^2\dP|^2\cdot|\Delta(D_tS_1\dP-D_t\dP)|dtdx \no\\
&\le C\ve\int_{\O_T}\bigg(\ss_{l=0}^3 R(t)^{\mu-1-\dl+2l}(\nabla_{t,x}^lD_t\dot\Phi)^2+\ss_{l=0}^2R(t)^{\mu-1-3(\g-1)+2l}(\nabla_{t,x}^l\nabla_x\dot\Phi)^2\bigg)dtdx.
\end{align}
Substituting (5.94) into (5.93) yields
\ben
&&|\int_{B_T}\f{R(t)^{\mu}}{\sqrt{1+L^2}}\hat{c}^2a(t)D_tS(S_1^2-S_1)\dP\cdot\p_rS(S_1^2-S_1)\dP dS|\no\\
&\le&C\ve\ss_{l=0}^2\int_{\Omega_T}\bigg( R(t)^{\mu-1-\dl+2l}(\nabla_{t,x}^lD_t\dot\Phi)^2
+ R(t)^{\mu-1-3(\g-1)+2l}(\nabla_{t,x}^l\nabla_x\dot\Phi)^2\bigg)dtdx\no\\
&&+C\ve\int_{\Omega_T}\bigg( R(t)^{\mu+5-\dl}(\nabla_{t,x}^3D_t\dot\Phi)^2
+ R(t)^{\mu+5-3(\g-1)-\dl}(\nabla_{t,x}^3\nabla_x\dot\Phi)^2\bigg)dtdx.
\een
Substituting (5.85) and (5.95) into (5.87), it follows
\ben
&&R(T)^{\mu}\int_{S_T}(D_tZ(S_1^2-S_1)\dP)^2dS+R(T)^{\mu-3(\g-1)}\int_{S_T}(\nabla_xZ(S_1^2-S_1)\dP)^2dS\no\\
&&+\int_{\Omega_T}\bigg(R(t)^{\mu-1-\dl}(D_tZ(S_1^2-S_1)\dP)^2+R(t)^{\mu-1-3(\g-1)}(\nabla_xZ(S_1^2-S_1)\dP)^2\bigg)dtdx\no\\
&\le&C\ve^2+\int_{\O_T}\bigg(Z(S_1^2-S_1)\mathcal{L}\dP+4\hat{c}^2\Delta ZS_1\dP
-2\hat{c}^2\Delta Z\dP\bigg)\cdot\mathcal{M}Z(S_1^2-S_1)\dP dtdx\no\\
&&+C\ve\ss_{l=0}^2\int_{\Omega_T}\bigg( R(t)^{\mu-1-\dl+2l}(\nabla_{t,x}^lD_t\dot\Phi)^2
+ R(t)^{\mu-1-3(\g-1)+2l}(\nabla_{t,x}^l\nabla_x\dot\Phi)^2\bigg)dtdx\no\\
&&+C\ve\int_{\Omega_T}\bigg( R(t)^{\mu+5-\dl}(\nabla_{t,x}^3D_t\dot\Phi)^2
+ R(t)^{\mu+5-3(\g-1)-\dl}(\nabla_{t,x}^3\nabla_x\dot\Phi)^2\bigg)dtdx.
\een

Similar to (5.75), we have
\ben
&&\int_{\O_T}\big(|\hat{c}^2\Delta ZS_1\dP|+|\hat{c}^2\Delta\dP|\big)\cdot|\mathcal{M}Z\big(S_1^2-S_1\big)\dP| dtdx\no\\
&\le&C\ve^2+C\ve\bigg(\int_{\Omega_T}\big(R(t)^{\mu+5-\dl}(\nabla_{t,x}^3D_t\dot\Phi)^2
+R(t)^{\mu+5-3(\g-1)-\dl}(\nabla_{t,x}^3\nabla_x\dot\Phi)^2\big) dtdx\bigg)^{\f12}.
\een

As (5.80)-(5.81) and (5.95), we have
\ben
&&|\int_{\Omega_T}Z(S_1^2-S_1)\mathcal{L}\dP\cdot\mathcal{M}Z(S_1^2-S_1)\dP dtdx|\no\\
&\le&C\ve^2+C\ve\bigg(R(T)^{\mu}\int_{S_T}(D_tZ(S_1^2-S_1)\dP)^2dS+R(T)^{\mu-3(\g-1)}\int_{S_T}(\nabla_xZ(S_1^2-S_1)\dP)^2dS\no\\
&&+\int_{\Omega_T}\bigg(\ss_{l=0}^2 R(t)^{\mu-1-\dl+2l}(\nabla_{t,x}^lD_t\dot\Phi)^2
+\ss_{l=0}^2 R(t)^{\mu-1-3(\g-1)+2l}(\nabla_{t,x}^l\nabla_x\dot\Phi)^2\bigg)dtdx\bigg)\no\\
&&+C\ve\int_{\Omega_T}\bigg( R(t)^{\mu+5-\dl}(\nabla_{t,x}^3D_t\dot\Phi)^2
+ R(t)^{\mu+5-3(\g-1)-\dl}(\nabla_{t,x}^3\nabla_x\dot\Phi)^2\bigg)dtdx\no\\
&&+C\ve\bigg(\int_{\Omega_T}\bigg( R(t)^{\mu+5-\dl}(\nabla_{t,x}^3D_t\dot\Phi)^2
+ R(t)^{\mu+5-3(\g-1)-\dl}(\nabla_{t,x}^3\nabla_x\dot\Phi)^2\bigg)dtdx\bigg)^{\f12}.
\een
Thus, together with (5.96)-(5.98), this
yields Lemma 5.4 for the case $S=Z$.

Next, we deal with the case of $S=S_0$. From (5.13), we have
\beq\int_{\O_T}D_tf\cdot gdtdx=\int_{\O_T}D_t(fg)-fD_tgdtdx=\bigg(\int_{S_T}-\int_{S^0}\bigg)fg dS
-\int_{\O_T}\bigg(\f{3L}{R(t)}fg+fD_tg\bigg)dtdx.
\eeq
If we choose $S=S_0$, then it follows from (5.88) and $[S_0,R(t)^2\Delta]=0$ that
\ben
&&\int_{\O_T}R(t)^{\mu}\hat{c}^2a(t)\bigg[\Delta(D_tS_0S_1\dP-D_t S_0\dP)\bigg(\chi\big(\f{r}{R(t)}\big)S_0\big(\f{G}{A}\big)+2\ss_{i=1}^3S_0Z_i^2\dP\bigg)\bigg]dtdx\no\\
&=&\int_{\O_T}D_t\bigg(R(t)^2\Delta(S_0S_1\dP-S_0\dP)\bigg)\cdot R(t)^{\mu-2}\hat{c}^2a(t)\bigg(\chi\big(\f{r}{R(t)}\big)S_0\big(\f{G}{A}\big)+2\ss_{i=1}^3S_0Z_i^2\dP\bigg)dtdx\no\\
&=&\bigg(\int_{S_T}-\int_{S^0}\bigg) R(t)^{\mu}\hat{c}^2a(t)\Delta(S_0S_1\dP-S_0\dP)\bigg(\chi\big(\f{r}{R(t)}\big)S_0\big(\f{G}{A}\big)+2\ss_{i=1}^3S_0Z_i^2\dP\bigg)dS\no\\
&&-\int_{\O_T}\bigg(R(t)^2\Delta(S_0S_1\dP-S_0\dP)\bigg)\cdot D_t\bigg(R(t)^{\mu-2}\hat{c}^2a(t)\big(\chi\big(\f{r}{R(t)}\big)S_0\big(\f{G}{A}\big)+2\ss_{i=1}^3S_0Z_i^2\dP\big)\bigg)dtdx\no\\
&&-\int_{\O_T}3LR(t)^{\mu}\hat{c}^2a(t)\Delta(S_0S_1\dP-S_0\dP)\bigg(\chi\big(\f{r}{R(t)}\big)S_0\big(\f{G}{A}\big)+2\ss_{i=1}^3S_0Z_i^2\dP\bigg)dtdx.
\een
And it follows from (\ref{apriori-better-bound}) that
\begin{align}
&|\chi\big(\f{r}{R(t)}\big)S_0\big(\f{G}{A}\big)|\le C\ve \bigg(R(t)^{-3(\g-1)}\ss_{l=0}^1(|S_0^lS_1^2\dP|+|S_0^lZ^2\dP|)+\ss_{l=0}^1(|S_0^lZ\dP|)\bigg),\\
&|\chi\big(\f{r}{R(t)}\big)S_0^2\big(\f{G}{A}\big)|\le C\ve \bigg(R(t)^{-3(\g-1)}\ss_{l=0}^2(|S_0^lS_1^2\dP|+|S_0^lZ^2\dP|)+\ss_{l=0}^2|S_0^lZ\dP|+R(t)|S_0^2D_t\dP|\bigg).
\end{align}
Thus, combining (5.88) and (5.100)-(5.102), we have
\ben
&&|\int_{B_T}\f{R(t)^{\mu}}{\sqrt{1+L^2}}\hat{c}^2a(t)D_tS_0(S_1^2-S_1)\dP\cdot\p_rS_0(S_1^2-S_1)\dP dS|\no\\
&\le&C\ve\ss_{l=0}^2\int_{\Omega_T}\bigg( R(t)^{\mu-1-\dl+2l}(\nabla_{t,x}^lD_t\dot\Phi)^2
+ R(t)^{\mu-1-3(\g-1)+2l}(\nabla_{t,x}^l\nabla_x\dot\Phi)^2\bigg)dtdx\no\\
&&+C\ve\int_{\Omega_T}\bigg( R(t)^{\mu+5-\dl}(\nabla_{t,x}^3D_t\dot\Phi)^2
+ R(t)^{\mu+5-3(\g-1)-\dl}(\nabla_{t,x}^3\nabla_x\dot\Phi)^2\bigg)dtdx.
\een

Substituting (5.86) and (5.103) into (5.87) yields
\begin{align}
&R(T)^{\mu}\int_{S_T}(D_tS_0(S_1^2-S_1)\dP)^2dS+R(T)^{\mu-3(\g-1)}\int_{S_T}(\nabla_xS_0(S_1^2-S_1)\dP)^2dS\no\\
&\q+\int_{\Omega_T}\bigg(R(t)^{\mu-1-\dl}(D_tS_0(S_1^2-S_1)\dP)^2+R(t)^{\mu-1-3(\g-1)}(\nabla_xS_0(S_1^2-S_1)\dP)^2\bigg)dtdx\no\\
&\le C\ve^2+C\ve\ss_{l=0}^2\int_{\Omega_T}\bigg( R(t)^{\mu-1-\dl+2l}(\nabla_{t,x}^lD_t\dot\Phi)^2
+ R(t)^{\mu-1-3(\g-1)+2l}(\nabla_{t,x}^l\nabla_x\dot\Phi)^2\bigg)dtdx\no\\
&\q+C\ve\int_{\Omega_T}\bigg( R(t)^{\mu+5-\dl}(\nabla_{t,x}^3D_t\dot\Phi)^2
+ R(t)^{\mu+5-3(\g-1)-\dl}(\nabla_{t,x}^3\nabla_x\dot\Phi)^2\bigg)dtdx\no\\
&\q+\int_{\O_T}\bigg(S_0(S_1^2-S_1)\mathcal{L}\dP-2L(S_1^2-S_1)\mathcal{L}\dP+S_0\big(4\hat{c}^2\Delta S_1\dP-2\hat{c}^2\Delta \dP\big)\no\\
&\q\q-8L\hat{c}^2\Delta S_1\dP+4L\hat{c}^2\Delta\dP+3(\g-1)\hat{c}^2\Delta(S_1^2-S_1)\dP\bigg)\cdot\mathcal{M}S_0(S_1^2-S_1)\dP dtdx.
\end{align}

Similar to (5.75), it follows
\begin{align}
&\int_{\O_T}\bigg(|S_0(\hat{c}^2\Delta S_1\dP)|+|S_0(\hat{c}^2\Delta \dP)|+|\hat{c}^2\Delta S_1\dP|+|\hat{c}^2\Delta\dP|+|\hat{c}^2\Delta(S_1^2-S_1)\dP|\bigg)
\cdot|\mathcal{M}S_0(S_1^2-S_1)\dP| dtdx\no\\
&\le C\ve^2+C\ve\bigg(\ss_{l=0}^2\int_{\Omega_T}\bigg( R(t)^{\mu-1-\dl+2l}(\nabla_{t,x}^lD_t\dot\Phi)^2
+ R(t)^{\mu-1-3(\g-1)+2l}(\nabla_{t,x}^l\nabla_x\dot\Phi)^2\bigg)dtdx\no\\
&\q+\int_{\Omega_T}\bigg( R(t)^{\mu+5-\dl}(\nabla_{t,x}^3D_t\dot\Phi)^2
+ R(t)^{\mu+5-3(\g-1)-\dl}(\nabla_{t,x}^3\nabla_x\dot\Phi)^2\bigg)dtdx\bigg)^{\f12}.
\end{align}
As  (5.80)-(5.81) and (5.95), we have
\ben
&&|\int_{\Omega_T}\big(|S_0(S_1^2-S_1)\mathcal{L}\dP|+|(S_1^2-S_1)\mathcal{L}\dP|\big)\mathcal{L}\dP\cdot\mathcal{M}S_0(S_1^2-S_1)\dP dtdx|\no\\
&\le&C\ve^2+C\ve\bigg(R(T)^{\mu}\int_{S_T}(D_tZ(S_1^2-S_1)\dP)^2dS+R(T)^{\mu-3(\g-1)}\int_{S_T}(\nabla_xZ(S_1^2-S_1)\dP)^2dS\no\\
&&+\int_{\Omega_T}\bigg(\ss_{l=0}^2 R(t)^{\mu-1-\dl+2l}(\nabla_{t,x}^lD_t\dot\Phi)^2
+\ss_{l=0}^2 R(t)^{\mu-1-3(\g-1)+2l}(\nabla_{t,x}^l\nabla_x\dot\Phi)^2\bigg)dtdx\bigg)\no\\
&&+C\ve\int_{\Omega_T}\bigg( R(t)^{\mu+5-\dl}(\nabla_{t,x}^3D_t\dot\Phi)^2
+ R(t)^{\mu+5-3(\g-1)-\dl}(\nabla_{t,x}^3\nabla_x\dot\Phi)^2\bigg)dtdx\no\\
&&+C\ve\bigg(\int_{\Omega_T}\bigg( R(t)^{\mu+5-\dl}(\nabla_{t,x}^3D_t\dot\Phi)^2
+ R(t)^{\mu+5-3(\g-1)-\dl}(\nabla_{t,x}^3\nabla_x\dot\Phi)^2\bigg)dtdx\bigg)^{\f12}.
\een
Together with (5.104)-(5.106), Lemma 5.4 is proved for the case of $S=S_0$
and it then completes the proof of the lemma. $\hfill\square$

\begin{lemma}{\bf (The fourth order radial derivative estimates)} Under the assumptions of Theorem 5.1,
we have
\ben
&&R(T)^{\mu-\dl}\int_{S_T}(D_tS_1^3\dP)^2dS+R(T)^{\mu-3(\g-1)-\dl}\int_{S_T}(\nabla_xS_1^3\dP)^2dS\no\\
&&+\int_{\Omega_T}\bigg(R(t)^{\mu-1-\dl}(D_tS_1^3\dP)^2+R(t)^{\mu-1-3(\g-1)-\dl}(\nabla_xS_1^3\dP)^2\bigg)dtdx\no\\
&\le&C\ve^2+C\ve\int_{\Omega_T}\ss_{l=0}^2 \bigg(R(t)^{\mu-1-\dl+2l}(\nabla_{t,x}^lD_t\dot\Phi)^2
+R(t)^{\mu-1-3(\g-1)+2l}(\nabla_{t,x}^l\nabla_x\dot\Phi)^2\bigg)dtdx\no\\
&&+C\ve\int_{\Omega_T}\bigg(R(t)^{\mu+5-\dl}(\nabla_{t,x}^lD_t\dot\Phi)^2
+R(t)^{\mu+5-3(\g-1)-\dl}(\nabla_{t,x}^l\nabla_x\dot\Phi)^2\bigg)dtdx\no\\
&&+C\ve\bigg(\ss_{l=0}^2\int_{\Omega_T}\bigg( R(t)^{\mu-1-\dl+2l}(\nabla_{t,x}^lD_t\dot\Phi)^2
+ R(t)^{\mu-1-3(\g-1)+2l}(\nabla_{t,x}^l\nabla_x\dot\Phi)^2\bigg)dtdx\no\\
&&+\int_{\Omega_T}\bigg( R(t)^{\mu+5-\dl}(\nabla_{t,x}^3D_t\dot\Phi)^2
+ R(t)^{\mu+5-3(\g-1)-\dl}(\nabla_{t,x}^3\nabla_x\dot\Phi)^2\bigg)dtdx\bigg)^{\f12}.
\een
\end{lemma}

\begin{remark}
Together with Lemma 5.1-5.5 and Remark 5.3-5.5, we obtain
\ben
&&\ss_{0\le l_1+l_2\le2}\bigg(R(T)^{\mu}\int_{S_T}(D_tS^{l_1}S_1^{l_2}\dP)^2dS
+R(T)^{\mu-3(\g-1)}\int_{S_T}(\nabla_xS^{l_1}S_1^{l_2}\dP)^2dS\bigg)\no\\
&&+\ss_{l_1+l_2=3,l_1\ge1}\bigg(R(T)^{\mu}\int_{S_T}(D_tS^{l_1}S_1^{l_2}\dP)^2dS
+R(T)^{\mu-3(\g-1)}\int_{S_T}(\nabla_xS^{l_1}S_1^{l_2}\dP)^2dS\bigg)\no\\
&&+R(T)^{\mu-\dl}\int_{S_T}(D_tS_1^3\dP)^2dS+R(T)^{\mu-3(\g-1)-\dl}\int_{S_T}(\nabla_xS_1^3\dP)^2dS\no\\
&&+\ss_{0\le l_1+l_2\le2}\int_{\Omega_T}\bigg(R(t)^{\mu-1-\dl}(D_tS^{l_1}S_1^{l_2}\dP)^2
+R(t)^{\mu-1-3(\g-1)}(\nabla_xS^{l_1}S_1^{l_2}\dP)^2\bigg)dtdx\no\\
&&+\ss_{l_1+l_2=3,l_1\ge1}\int_{\Omega_T}\bigg(R(t)^{\mu-1-\dl}(D_tS^{l_1}S_1^{l_2}\dP)^2
+R(t)^{\mu-1-3(\g-1)}(\nabla_xS^{l_1}S_1^{l_2}\dP)^2\bigg)dtdx\no\\
&&+\int_{\Omega_T}\bigg(R(t)^{\mu-1-\dl}(D_tS_1^3\dP)^2
+R(t)^{\mu-1-3(\g-1)-\dl}(\nabla_xS_1^3\dP)^2\bigg)dtdx\no\\
&\le&C\ve^2+C\ve\int_{\Omega_T}\ss_{l=0}^2\bigg(R(t)^{\mu-1-\dl+2l}(\nabla_{t,x}^lD_t\dot\Phi)^2
+R(t)^{\mu-1-3(\g-1)+2l}(\nabla_{t,x}^l\nabla_x\dot\Phi)^2\bigg)dtdx\no\\
&&+C\ve\int_{\Omega_T}\bigg(R(t)^{\mu+5-\dl}(\nabla_{t,x}^lD_t\dot\Phi)^2
+R(t)^{\mu+5-3(\g-1)-\dl}(\nabla_{t,x}^l\nabla_x\dot\Phi)^2\bigg)dtdx\no\\
&&+C\ve\bigg(\int_{\Omega_T}\ss_{l=0}^2\bigg(R(t)^{\mu-1-\dl+2l}(\nabla_{t,x}^lD_t\dot\Phi)^2
+R(t)^{\mu-1-3(\g-1)+2l}(\nabla_{t,x}^l\nabla_x\dot\Phi)^2\bigg)dtdx\no\\
&&+\int_{\Omega_T}\bigg(R(t)^{\mu+5-\dl}(\nabla_{t,x}^lD_t\dot\Phi)^2
+R(t)^{\mu+5-3(\g-1)-\dl}(\nabla_{t,x}^l\nabla_x\dot\Phi)^2\bigg)dtdx\bigg)^{\f12}.
\een
\end{remark}

{\bf Proof of Lemma 5.5.} Set
\begin{align}
N_1&=S_1^3-S_1^2-2\ss_{i=1}^3Z_i^2\dP,\\
N_2&=\chi\big(\f{r}{R(t)}\big)\f1A\bigg(-\f{\g-1}{R(t)^2}\ss_{i=1}^3(Z_i\dP)^2\bigg(S_1^2
+\ss_{i=1}^3Z_i^2\bigg)-\f{3(\g-1)L}{R(t)}\ss_{i=1}^3Z_i\dP Z_i-\f{2L}{R(t)}\ss_{i,j=1}^3C_{ij}Z_i\dP Z_j\no\\
&-\f{4}{R(t)^2}\ss_{i,j=1}^3Z_i\dP Z_j\dP Z_iZ_j+\f{1}{R(t)^2}\ss_{i,j=1}^3C_{ij}Z_i\dP Z_j\dP S_1^2
-\f{4}{R(t)^2}\ss_{i,j,k=1}^3C_{ijk}Z_i\dP Z_j\dP Z_k\bigg).
\end{align}
Then it follows from (5.66) that on $B_T$
\beq
N\dP\equiv (N_1-N_2)\dP=0.
\eeq
On the other hand,
\beq
\mathcal{L}N\dP=N\mathcal{L}\dP+[\mathcal{L},N_1]\dP-[\mathcal{L},N_2]\dP,
\eeq
where
\beq
[\mathcal{L},N_1]\dP=-6\hat{c}^2\Delta S_1^2\dP+8\hat{c}^2\Delta S_1\dP-4\hat{c}^2\Delta\dP,
\eeq
and
\ben
|[\mathcal{L},N_2]\dP|&\le& C\ve\chi\big(\f{r}{R(t)}\big)\bigg(R(t)^{-3(\g-1)-2}\ss_{0\le l\le1}\big(|S_0^lS_1^2\dP|+|S_0^lZ^2\dP|\big)+R(t)^{-2}\ss_{0\le l\le2}|S_0^lZ\dP|\no\\
&&+R(t)^{-3(\g-1)+\dl-2}\ss_{0\le l\le1}R(t)^l\big(|\nabla_x^lS_1^2\dP|+|\nabla_x^lZ^2\dP|\big)+R(t)^{-1}|S_0^2D_t\dP|\no\\
&&+R(t)^{-\f{3(\g-1)}2+\f{\dl}2-2}\ss_{1\le l\le2}R(t)^l|\nabla_x^lZ\dP|+R(t)^{-3(\g-1)+1}|\nabla_x^2D_t\dP|\bigg).
\een
Note that the term $R(t)^{-1}S_0^2D_t\dP$ is on the right hand side of (5.114). Thus,  we can not choose
the multiplier $\mathcal{M}N\dP$ as in Theorem 4.1
to derive the energy estimate of $N\dP$ from (5.112).
Otherwise, we can not control the term
$\int_{\O_T}R(t)^{\mu-1}S_0^2D_t\dP\cdot D_tN_1\dP dtdx$ in  $\int_{\O_T}[\mathcal{L},N_2]\dP\cdot\mathcal{M}N\dP dtdx$.


Instead, we will use a new multiplier $\tilde{\mathcal{M}} N\dP=R(t)^{\mu-\dl}D_tN\dP$. It follows  that
\begin{align}
&\int_{\Omega_T}\mathcal{L}N\dP\cdot \tilde{\mathcal{M}}N\dP dtdx\no\\
&=\bigg(\int_{S_T}-\int_{S^0}\bigg)\f12R(t)^{\mu-\dl}\big((D_tN\dP)^2+\hat{c}^2|\nabla_xN\dP|^2\big)dS\no\\
&\quad+\int_{\O_T}\biggl\{\f{\dl}2R(t)^{\mu-1-\dl}\big(D_tN\dP\big)^2+\f{\g}2R(t)^{\mu-1-3(\g-1)-\dl}L(5-3\g)|\nabla_xN\dP|^2
\biggr\}dtdx.
\end{align}
This yields
\ben
&&R(T)^{\mu-\dl}\int_{S_T}(D_tN\dP)^2dS+R(T)^{\mu-3(\g-1)-\dl}\int_{S_T}|\nabla_xN\dP|^2dS\no\\
&&+\int_{\O_T}\bigg(R(t)^{\mu-1-\dl}(D_tN\dP)^2+R(t)^{\mu-1-3(\g-1)-\dl}|\nabla_xN\dP|^2\bigg)dtdx\no\\
&\le& C\ve^2+C\int_{\O_T}\bigg(N\mathcal{L}\dP+[\mathcal{L},N_1]\dP-[\mathcal{L},N_2]\dP\bigg)\cdot\tilde{\mathcal{M}}N\dP dtdx.
\een
It follows from (5.101) and (5.114) that
\ben
&&\int_{\O_T}|[\mathcal{L},N_2]\dP|\cdot |\tilde{\mathcal{M}}N\dP|dtdx\no\\
&\le&C\ve^2+C\ve\int_{\Omega_T}\ss_{l=0}^2\bigg(R(t)^{\mu-1-\dl+2l}(\nabla_{t,x}^lD_t\dot\Phi)^2
+R(t)^{\mu-1-3(\g-1)+2l}(\nabla_{t,x}^l\nabla_x\dot\Phi)^2\bigg)dtdx\no\\
&&+C\ve\int_{\Omega_T}\bigg(R(t)^{\mu+5-\dl}(\nabla_{t,x}^lD_t\dot\Phi)^2
+R(t)^{\mu+5-3(\g-1)-\dl}(\nabla_{t,x}^l\nabla_x\dot\Phi)^2\bigg)dtdx,
\een
where we have used
\beq\int_{\O_T}R(t)^{\mu-3-\dl}(\chi\big(\f{r}{R(t)}\big)Z\dP)^2dtdx\le C\ve^2,\eeq
because $|Z\dP|\le C\ve R(t)^{1-3(\g-1)}$ holds for $r>\f13R(t)$.

Thus, by (5.83) and (5.101), we have
\ben
&&\int_{\O_T}|[\mathcal{L},N_1]\dP|\cdot |\tilde{\mathcal{M}}N\dP|dtdx\no\\
&\le&\int_{\O_T}|[\mathcal{L},N_1]\dP|\cdot R(t)^{\mu-\dl}\big(|D_tS_1^3\dP|+|D_tS_1^2\dP|
+\ss_{i=1}^3 |D_t Z_i^2\dP|+|D_tN_2\dP|\big)dtdx\no\\
&\le&C\ve^2+C\ve\int_{\Omega_T}\ss_{l=0}^2\bigg(R(t)^{\mu-1-\dl+2l}(\nabla_{t,x}^lD_t\dot\Phi)^2
+R(t)^{\mu-1-3(\g-1)+2l}(\nabla_{t,x}^l\nabla_x\dot\Phi)^2\bigg)dtdx\no\\
&&+C\ve\int_{\Omega_T}\bigg(R(t)^{\mu+5-\dl}(\nabla_{t,x}^lD_t\dot\Phi)^2
+R(t)^{\mu+5-3(\g-1)-\dl}(\nabla_{t,x}^l\nabla_x\dot\Phi)^2\bigg)dtdx\no\\
&&+C\ve\bigg( \int_{\Omega_T}\bigg(R(t)^{\mu+5-\dl}(\nabla_{t,x}^lD_t\dot\Phi)^2
+R(t)^{\mu+5-3(\g-1)-\dl}(\nabla_{t,x}^l\nabla_x\dot\Phi)^2\bigg)dtdx\bigg)^{\f12}.
\een

Next, for $\int_{\O_T}N\mathcal{L}\dP\cdot \tilde{\mathcal{M}}N\dP dtdx$,
 direct computation yields
\beq
N\mathcal{L}\dP=M_1+M_2,
\eeq
where
\ben
&M_1&=\ss_{i=1}^3\tf_{0i}\p_iD_tN\dP+\ss_{i,j=1}^3\tf_{ij}\p_{ij}N\dP,\no\\
&M_2&=\ss_{i=1}^3[\tf_{0i}\p_iD_t,N]\dP+\ss_{i,j=1}^3[\tf_{ij}\p_{ij},N]\dP.\no
\een
Then similar to  (5.16), we have
\ben
&&|\int_{\Omega_T}M_1\cdot\tilde{\mathcal{M}}N\dP dtdx|\no\\
&\le&C\ve^2+C\ve\bigg(R(T)^{\mu-\dl}\int_{S_T}(D_tN\dP)^2dS+R(T)^{\mu-3(\g-1)-\dl}\int_{S_T}(\nabla_xN\dP)^2dS\no\\
&&+\int_{\O_T}\big(R(t)^{\mu-1-\dl}(D_tN\dP)^2+R(t)^{\mu-1-3(\g-1)-\dl}(\nabla_xN\dP)^2dtdx\bigg).
\een
Moreover, by the expression of $N$ and Lemma 5.1-5.4, Remark 5.2-5.3, we have
\ben
&&|\int_{\Omega_T}M_1\cdot\tilde{\mathcal{M}}N\dP dtdx|\no\\
&\le&C\ve^2+C\ve\int_{\Omega_T}\ss_{l=0}^2\bigg(R(t)^{\mu-1-\dl+2l}(\nabla_{t,x}^lD_t\dot\Phi)^2
+R(t)^{\mu-1-3(\g-1)+2l}(\nabla_{t,x}^l\nabla_x\dot\Phi)^2\bigg)dtdx\no\\
&&+C\ve\int_{\Omega_T}\bigg(R(t)^{\mu+5-\dl}(\nabla_{t,x}^lD_t\dot\Phi)^2
+R(t)^{\mu+5-3(\g-1)-\dl}(\nabla_{t,x}^l\nabla_x\dot\Phi)^2\bigg)dtdx.
\een
Then as for (5.18)-(5.19), (5.30), (5.45)-(5.47), (5.51)-(5.52) and (5.57), we obtain
\ben
&&|\int_{\Omega_T}M_2\cdot\tilde{\mathcal{M}}N\dP dtdx|\no\\
&\le&C\ve^2+C\ve\int_{\Omega_T}\ss_{l=0}^2\bigg(R(t)^{\mu-1-\dl+2l}(\nabla_{t,x}^lD_t\dot\Phi)^2
+R(t)^{\mu-1-3(\g-1)+2l}(\nabla_{t,x}^l\nabla_x\dot\Phi)^2\bigg)dtdx\no\\
&&+C\ve\int_{\Omega_T}\bigg(R(t)^{\mu+5-\dl}(\nabla_{t,x}^lD_t\dot\Phi)^2
+R(t)^{\mu+5-3(\g-1)-\dl}(\nabla_{t,x}^l\nabla_x\dot\Phi)^2\bigg)dtdx.
\een
Finally, substituting (5.117), (5.119) and (5.122)-(5.123) into (5.116)
completes the proof of the lemma.

\qquad $\hfill\square$

As shown in Remark 5.6, the operators $S_1$ and $Z$ are used in the  energy estimates. However,
$R(t)\p_i$ $(1\le i\le 3)$ are equivalent to $S_1$ and $Z$ only when $r>\f13R(t)$. Hence, we
still need to obtain estimates for $r\le \f13R(t)$. For this purpose, set $\nu(s)=1-\chi(s)$
supported in $[0, \f23]$, with  $\chi(s)$ defined in (5.67).

\begin{lemma}{\bf (Estimates near $r=0$)} Under the assumptions of Theorem 5.1, we have
\begin{align}
&\int_{S_T}\bigg(R(T)^{\mu}\nu\big(\f{r}{R(t)}\big)\big(D_t(R(t)\nabla_{x})^k\dot\Phi\big)^2
+R(T)^{\mu-3(\g-1)}\nu\big(\f{r}{R(t)}\big)\big(\nabla_x(R(t)\nabla_{x})^k\dot\Phi\big)^2\bigg)dS\no\\
&+\int_{\Omega_T}\bigg(R(t)^{\mu-1-\dl}\nu\big(\f{r}{R(t)}\big)\big(D_t(R(t)\nabla_{x})^k\dot\Phi\big)^2
+R(t)^{\mu-1-3(\g-1)+2k}\nu\big(\f{r}{R(t)}\big)\big(\nabla_x(R(t)\nabla_{x})^k\dot\Phi\big)^2\bigg)dtdx\no\\
&\le C\ve^2+C\ve\int_{\Omega_T}\ss_{l=0}^2\bigg(R(t)^{\mu-1-\dl+2l}(\nabla_{t,x}^lD_t\dot\Phi)^2
+R(t)^{\mu-1-3(\g-1)+2l}(\nabla_{t,x}^l\nabla_x\dot\Phi)^2\bigg)dtdx\no\\
&+C\ve\int_{\Omega_T}\bigg(R(t)^{\mu+5-\dl}(\nabla_{t,x}^lD_t\dot\Phi)^2
+R(t)^{\mu+5-3(\g-1)-\dl}(\nabla_{t,x}^l\nabla_x\dot\Phi)^2\bigg)dtdx,
\end{align}
and
\begin{align}
&\int_{S_T}\bigg(R(T)^{\mu-\dl}\nu\big(\f{r}{R(t)}\big)\big(D_t(R(t)\nabla_{x})^3\dot\Phi\big)^2
+R(T)^{\mu-\dl-3(\g-1)}\nu\big(\f{r}{R(t)}\big)\big(\nabla_x(R(t)\nabla_{x})^3\dot\Phi\big)^2\bigg)dS\no\\
&\quad +\int_{\Omega_T}\bigg(R(t)^{\mu-1-\dl}\nu\big(\f{r}{R(t)}\big)\big(D_t(R(t)\nabla_{x})^3\dot\Phi\big)^2
+R(t)^{\mu-\dl-1-3(\g-1)}\nu\big(\f{r}{R(t)}\big)\big(\nabla_x(R(t)\nabla_{x})^3\dot\Phi\big)^2\bigg)dtdx\no\\
&\le C\ve^2+C\ve\int_{\Omega_T}\ss_{l=0}^2\bigg(R(t)^{\mu-1-\dl+2l}(\nabla_{t,x}^lD_t\dot\Phi)^2
+R(t)^{\mu-1-3(\g-1)+2l}(\nabla_{t,x}^l\nabla_x\dot\Phi)^2\bigg)dtdx\no\\
&\quad +C\ve\int_{\Omega_T}\bigg(R(t)^{\mu+5-\dl}(\nabla_{t,x}^lD_t\dot\Phi)^2
+R(t)^{\mu+5-3(\g-1)-\dl}(\nabla_{t,x}^l\nabla_x\dot\Phi)^2\bigg)dtdx.
\end{align}
\end{lemma}

{\bf Proof.} By  direct computation, we have
\beq
\big(R(t)\p_i\big)\mathcal{L}=\mathcal{L}\big(R(t)\p_i\big).
\eeq
This, together with Theorem 4.1, yields for $0\le k\le2$,
\begin{align}
&\int_{S_T}\bigg(R(T)^{\mu}\nu\big(\f{r}{R(t)}\big)\big(D_t(R(t)\nabla_{x})^k\dot\Phi\big)^2
+R(T)^{\mu-3(\g-1)}\nu\big(\f{r}{R(t)}\big)\big(\nabla_x(R(t)\nabla_{x})^k\dot\Phi\big)^2\bigg)dS\no\\
&+\int_{\Omega_T}\bigg(R(t)^{\mu-1-\dl}\nu\big(\f{r}{R(t)}\big)\big(D_t(R(t)\nabla_{x})^k\dot\Phi\big)^2
+R(t)^{\mu-1-3(\g-1)+2k}\nu\big(\f{r}{R(t)}\big)\big(\nabla_x(R(t)\nabla_{x})^k\dot\Phi\big)^2\bigg)dtdx\no\\
&\le
C\bigg(\int_{\O_T}R(t)^{\mu-1-3(\g-1)}|\nu'\big(\f{r}{R(t)}\big)|\big(\nabla_x(R(t)\nabla)^k\dP\big)^2dtdx\bigg)^{\f12}\no\\
&\quad \times\bigg(\int_{\O_T}R(t)^{\mu-1-\dl}|\nu'\big(\f{r}{R(t)}\big)|\big(D_t(R(t)\nabla)^k\dP\big)^2dtdx\bigg)^{\f12}
\no\\
&\quad +C\ve^2+\int_{\O_T}\mathcal{L}(R(t)\nabla)^k\dP\cdot\nu\big(\f{r}{R(t)}\big)\mathcal{M}(R(t)\nabla)^k\dP dtdx.
\end{align}
Note that the function $\nu'\big(\ds\f{r}{R(t)}\big)$ has a compact support away from $r=0$ so that the
first term on the right hand side of (5.127) can be estimated as in
Remark 5.6.

On the other hand, by a similar argument for $\int_{\O_T}\mathcal{L}Z^k\dP\cdot \mathcal{M}Z^k\dP dtdx$ in (5.43),
we can obtain
\begin{align}
&|\int_{\O_T}\mathcal{L}(R(t)\nabla)^k\dP\cdot\nu\big(\f{r}{R(t)}\big)\mathcal{M}(R(t)\nabla)^k\dP dtdx|\no\\
&\le C\ve^2+C\ve\int_{\Omega_T}\ss_{l=0}^2\bigg(R(t)^{\mu-1-\dl+2l}(\nabla_{t,x}^lD_t\dot\Phi)^2
+R(t)^{\mu-1-3(\g-1)+2l}(\nabla_{t,x}^l\nabla_x\dot\Phi)^2\bigg)dtdx\no\\
&\quad+C\ve\int_{\Omega_T}\bigg(R(t)^{\mu+5-\dl}(\nabla_{t,x}^lD_t\dot\Phi)^2
+R(t)^{\mu+5-3(\g-1)-\dl}(\nabla_{t,x}^l\nabla_x\dot\Phi)^2\bigg)dtdx.
\end{align}
Then (5.124) follows from (5.127)-(5.128) and Remark 5.6.
In addition, similar to (5.124), we obtain (5.125).
And this completes the proof of the lemma. $\hfill\square$.

Based on Lemma 5.1-Lemma 5.6 and Remark 5.1-Remark 5.6, we are ready to prove Theorem 5.1.

{\bf Proof of Theorem 5.1.}
By Remark 5.6 and Lemma 5.6, for sufficiently small $\ve>0$, we have
\begin{align}
&\ss_{k=0}^2\int_{S_T}\bigg(R(T)^{\mu+2k}(\nabla_{t,x}^kD_t\dot\Phi)^2+R(T)^{\mu-3(\g-1)+2k}(\nabla_{t,x}^k\nabla_x\dot\Phi)^2\bigg)dS\no\\
&\q+\int_{S_T}\bigg(R(T)^{\mu+6-\dl}(\nabla_{t,x}^3D_t\dot\Phi)^2+R(T)^{\mu-3(\g-1)+6-\dl}(\nabla_{t,x}^k\nabla_x\dot\Phi)^2\bigg)dS\no\\
&\q+\ss_{k=0}^2\int_{\Omega_T}\bigg(R(t)^{\mu-1-\dl+2k}(\nabla_{t,x}^kD_t\dot\Phi)^2
+R(t)^{\mu-1-3(\g-1)+2k}(\nabla_{t,x}^k\nabla_x\dot\Phi)^2\bigg)dtdx\no\\
&\q+\int_{\Omega_T}\bigg(R(t)^{\mu+5-\dl}(\nabla_{t,x}^3D_t\dot\Phi)^2
+R(t)^{\mu+5-3(\g-1)-\dl}(\nabla_{t,x}^3\nabla_x\dot\Phi)^2\bigg)dtdx\no\\
&\le C\ve^2+C\ve\bigg(\int_{\Omega_T}\ss_{l=0}^2\bigg(R(t)^{\mu-1-\dl+2l}(\nabla_{t,x}^lD_t\dot\Phi)^2
+R(t)^{\mu-1-3(\g-1)+2l}(\nabla_{t,x}^l\nabla_x\dot\Phi)^2\bigg)dtdx\no\\
&\q+\int_{\Omega_T}\bigg(R(t)^{\mu+5-\dl}(\nabla_{t,x}^lD_t\dot\Phi)^2
+R(t)^{\mu+5-3(\g-1)-\dl}(\nabla_{t,x}^l\nabla_x\dot\Phi)^2\bigg)dtdx\bigg)^{\f12}.
\end{align}
If
\begin{align}
&\int_{\Omega_T}\ss_{l=0}^2\bigg(R(t)^{\mu-1-\dl+2l}(\nabla_{t,x}^lD_t\dot\Phi)^2
+R(t)^{\mu-1-3(\g-1)+2l}(\nabla_{t,x}^l\nabla_x\dot\Phi)^2\bigg)dtdx\no\\
&\quad +\int_{\Omega_T}\bigg(R(t)^{\mu+5-\dl}(\nabla_{t,x}^lD_t\dot\Phi)^2
+R(t)^{\mu+5-3(\g-1)-\dl}(\nabla_{t,x}^l\nabla_x\dot\Phi)^2\bigg)dtdx\no\\
&\le C\ve^2,
\end{align}
then (5.2) and (5.3) can be derived directly.

If
\begin{align}
&\int_{\Omega_T}\ss_{l=0}^2\bigg(R(t)^{\mu-1-\dl+2l}(\nabla_{t,x}^lD_t\dot\Phi)^2
+R(t)^{\mu-1-3(\g-1)+2l}(\nabla_{t,x}^l\nabla_x\dot\Phi)^2\bigg)dtdx\no\\
&\quad +\int_{\Omega_T}\bigg(R(t)^{\mu+5-\dl}(\nabla_{t,x}^lD_t\dot\Phi)^2
+R(t)^{\mu+5-3(\g-1)-\dl}(\nabla_{t,x}^l\nabla_x\dot\Phi)^2\bigg)dtdx\no\\
&\ge C\ve^2,
\end{align}
then it follows from (5.129) that
\ben
&&\int_{\Omega_T}\ss_{l=0}^2\bigg(R(t)^{\mu-1-\dl+2l}(\nabla_{t,x}^lD_t\dot\Phi)^2
+R(t)^{\mu-1-3(\g-1)+2l}(\nabla_{t,x}^l\nabla_x\dot\Phi)^2\bigg)dtdx\no\\
&&+\int_{\Omega_T}\bigg(R(t)^{\mu+5-\dl}(\nabla_{t,x}^lD_t\dot\Phi)^2
+R(t)^{\mu+5-3(\g-1)-\dl}(\nabla_{t,x}^l\nabla_x\dot\Phi)^2\bigg)dtdx\no\\
&\le& C\ve\bigg(\int_{\Omega_T}\ss_{l=0}^2\bigg(R(t)^{\mu-1-\dl+2l}(\nabla_{t,x}^lD_t\dot\Phi)^2
+R(t)^{\mu-1-3(\g-1)+2l}(\nabla_{t,x}^l\nabla_x\dot\Phi)^2\bigg)dtdx\no\\
&&+\int_{\Omega_T}\bigg(R(t)^{\mu+5-\dl}(\nabla_{t,x}^lD_t\dot\Phi)^2
+R(t)^{\mu+5-3(\g-1)-\dl}(\nabla_{t,x}^l\nabla_x\dot\Phi)^2\bigg)dtdx\bigg)^{\f12},
\een
which implies
\begin{align}
&\int_{\Omega_T}\ss_{l=0}^2\bigg(R(t)^{\mu-1-\dl+2l}(\nabla_{t,x}^lD_t\dot\Phi)^2
+R(t)^{\mu-1-3(\g-1)+2l}(\nabla_{t,x}^l\nabla_x\dot\Phi)^2\bigg)dtdx\no\\
&\quad+\int_{\Omega_T}\bigg(R(t)^{\mu+5-\dl}(\nabla_{t,x}^lD_t\dot\Phi)^2
+R(t)^{\mu+5-3(\g-1)-\dl}(\nabla_{t,x}^l\nabla_x\dot\Phi)^2\bigg)dtdx\no\\
&\le C\ve^2.\no
\end{align}
Substituting this into (5.129) derives (5.2) and (5.3),
and then it   completes the proof of Theorem 5.1. $\hfill\square$

\section{Proofs of Theorem 1.1 and 1.2.}

To complete the proof of Theorem 1.2,  as in [14], the following
estimate is needed.\medskip

{\bf Lemma 6.1.}  {\it For $1\leq t\leq T_0$ and $k_0\ge 4$, we have that for
a smooth function $\vp(t,x)\in H^{k_0}(\O_T)$
\beq
\ds\sum_{0\le l\le k_0-4}|t^l\na^{l+1}\vp(t,x)|^2
\le C_0 t^{-3}\int_{S_T}
\ds\sum_{0\le l\le k_0-1} |t^l\na^{l+1}\vp(t,x)|^2 dx.
\eeq}

{\bf Proof.}
For any $t_1\in [1,T_0]$,  set
$$(t',x')=\f{1}{t_1}(t,x).$$
Then
\beq\na_{t,x}^k\vp=\f{1}{t_1^k}\na_{t',x'}^k\vp, \quad \forall\
k\in\Bbb N.
\eeq
Define $D_{*}=\{(t',x'): t'=1, |x'|\leq
L\}$. Then by the Sobolev imbedding theorem and by noting that $D_{*}$
satisfies the uniform interior cone condition,
$$|\na_{t',x'}\vp|^2(1,x')|\leq C\int_{D_{*}}\ds\sum_{0\leq l\leq
3}|\na_{t',x'}^{l+1}\vp|^2(1,x')dx'.$$
In view of (6.2), one has
\begin{align*}
|\na_{t,x}\vp|^2(t_1,x)=&\f{1}{t_1^2}|\na_{t',x'}\vp|^2(1,x')\\
\leq &\f{C}{t_1^2}\int_{D_{*}}\ds\sum_{0\leq l\leq
3}|\na_{t',x'}^{l+1}\vp|^2(1,x')dx'\\
=&\f{C}{t_1^2}\int_{S_T}\ds\sum_{0\leq
l\leq 3}|t_1^{l+1}\na_{t,x}^{l+1}\vp|^2(t_1,x) \f{1}{t_1^3} dx\\
=&\f{C}{t_1^3}\int_{S_T}\ds\sum_{0\leq
l\leq 3}|t_1^{l}\na_{t,x}^{l+1}\vp|^2(t_1,x)  dx.\end{align*}

This yields (6.1) for $l=0$. The cases of $1\le l\le k_0-4$ can be
estimated similarly and this completes the proof of the lemma.
 $\hfill\square$

\medskip

We now first prove Theorem 1.2 as follows.

{\bf Proof of Theorem 1.2.} It follows from Lemma 6.1 that, for $0\le t\le T$,
\beq\left\{\begin{aligned}
&\sum\limits_{0\le l\le 1}|R(t)^l\nabla_x^{l+1}\dot\Phi|^2\le
CR(t)^{-3}\int_{S_t}\sum\limits_{0\le l\le 3}
|R(t)^l\nabla_x^{l+1}\dot\Phi|^2dS,\\
& |R(t)D_t^{2}\dP|^2\le
CR(t)^{-3}\int_{S_t}\sum\limits_{0\le l\le 3}
|R(t)^{1+l}\nabla_x^lD_t^{2}\dP|^2dS,\\
& \sum\limits_{0\le l\le
1}|R(t)^l\nabla_x^lD_t\dP|^2\le
CR(t)^{-3}\int_{S_t}\sum\limits_{0\le l\le 3}
|R(t)^{l}\nabla_x^{l}D_t\dP|^2dS.
\end{aligned}\right.
\eeq
On the other hand, (5.2) gives
\beq\left\{\begin{aligned}
&\int_{S_t}\sum\limits_{0\le l\le 2}|R(t)^l\nabla_x^{l+1}\dot\Phi|^2dS \le C\ve^2 R(t)^{-\mu+3(\g-1)},\\
&\int_{S_t}\sum\limits_{0\le l\le 2}|R(t)^{l}\nabla_x^{l}D_t\dP|^2dS \le C\ve^2 R(t)^{-\mu},
\end{aligned}\right.
\eeq
and
\beq\left\{\begin{aligned}
&\int_{S_t}\sum\limits_{0\le l\le 2}|R(t)^{1+l}\nabla_x^{l+2}\dot\Phi|^2dS \le C\ve^2 R(t)^{-\mu+3(\g-1)+\dl},\\
&\int_{S_t}\sum\limits_{0\le l\le 2}|R(t)^{1+l}\nabla_x^{1+l}D_t\dP|^2dS \le C\ve^2 R(t)^{-\mu+\dl},\\
&\int_{S_t}\sum\limits_{0\le l\le 3}|R(t)^{1+l}\nabla_x^lD_t^{2}\dP|^2dS \le C\ve^2 R(t)^{-\mu+\dl}.
\end{aligned}\right.
\eeq
Hence, we obtain
\beq\left\{\begin{aligned}
&|\nabla_x\dot\Phi|\le C\ve R(t)^{-3(\g-1)+\f{\dl}2},\quad |R(t)^{l-1}D_t^l\dot\Phi|\le C\ve R(t)^{-3(\g-1)},(l=1,2)\\
&|R(t)\nabla_xD_t\dot\Phi|\le C\ve R(t)^{-3(\g-1)+\f{\dl}2},\quad |R(t)\nabla_x^2\dot\Phi|\le C\ve R(t)^{-\f{3(\g-1)-\dl}2}.
\end{aligned}\right.
\eeq
On the other hand, by the stream line equation
\beq
\f{dx_i(t)}{dt}=\f{Lx_i}{R(t)},\q x_i(0)=x_i^0,\q(i=1,2,3),
\eeq
we have
\beq
x_i(t)=x_i^0R(t)\q (i=1,2,3),
\eeq
where  $(x_1^0,x_2^0,x_3^0)$ is the initial point.
Then integrating along the stream line, we have   for $1<\g<\ds\f{4}3$ that
$$
|R(t)\nabla_x\dP(t,x(t))|\le|\nabla_x\dP(0,x(0))|+\int_0^t|D_t(R(t)\nabla_x\dP)|dt\le C\ve\big(1+R(t)^{1-3(\g-1)+\f{\dl}2}\big),$$
which implies
\beq
|\nabla_x\dP|\le C\ve\big(R(t)^{-1}+R(t)^{-3(\g-1)+\f{\dl}2}\big) \le C\ve R(t)^{-3(\g-1)+\f{\dl}2}.
\eeq
Thus, the first inequality in (\ref{apriori-bound}) is proved.

On the other hand, it follows from Remark 5.6 and Theorem 5.1 that
\begin{align}
&\ss_{0\le l_1+l_2\le2}\bigg(R(t)^{\mu}\int_{S_t}(D_tS^{l_1}S_1^{l_2}\dP)^2dS
+R(t)^{\mu-3(\g-1)}\int_{S_t}(\nabla_xS^{l_1}S_1^{l_2}\dP)^2dS\bigg)\no\\
&\quad+\ss_{l_1+l_2=3,l_1\ge1}\bigg(R(t)^{\mu}\int_{S_t}(D_tS^{l_1}S_1^{l_2}\dP)^2dS
+R(t)^{\mu-3(\g-1)}\int_{S_t}(\nabla_xS^{l_1}S_1^{l_2}\dP)^2dS\bigg)\no\\
&\le C\ve^2.
\end{align}
Noticing that for $r>\ds\f13R(t)$, $r\sim R(t)$ holds. Then (6.10) implies
\begin{align*}
&\ss_{0\le l\le2}\bigg(R(t)^{\mu+2l}\int_{S_t\bigcap\{r>\f13R(t)\}}\big(\nabla_x^lZD_t\dP\big)^2dS
+R(t)^{\mu-3(\g-1)+2l}\int_{S_t\bigcap\{r>\f13R(t)\}}\big(\nabla_x^{l+1}Z\dP\big)^2dS\bigg)\no\\
&\le C\ve^2.
\end{align*}
Together with Lemma 6.1, this yields for $r>\f13R(t)$,
\ben
&&|ZD_t\dP|^2\le CR(t)^{-3}\int_{S_t\bigcap\{r>\f13R(t)\}}\ss_{0\le l\le 2}|R(t)^l\nabla_x^lZD_t\dP|^2dS
\le C\ve R(t)^{-\mu-3},\no\\
&&|\nabla_xD_t\dP|^2\le CR(t)^{-3}\int_{S_t\bigcap\{r>\f13R(t)\}}\ss_{0\le l\le 2}|R(t)^l\nabla_x^{l+1}Z\dP|^2dS
\le C\ve R(t)^{-\mu+3(\g-1)-3}.\no
\een
Subsequently, one has
$$
|ZD_t\dP|\le C\ve R(t)^{-3(\g-1)},\q |\nabla_xD_t\dP|\le C\ve R(t)^{-\f{3(\g-1)}2}.
$$
In addition, for $1<\g<\ds\f{4}3$ and $r>\ds\f13R(t)$, integrating along the stream yields
$$
|Z\dP(t,x(t))|\le|Z\dP(0,x(0))|+\int_0^t|D_t(Z\dP)|dt\le C\ve\big(1+R(t)^{1-3(\g-1)}\big)\le C\ve R(t)^{1-3(\g-1)}.$$

Note that the generic constant  $C$ appeared in this section
 depends only on the initial data. Then we can choose the constant
$M=2C$ in (5.1)-(5.2) for small $\ve>0$
so that (5.1) and (\ref{apriori-better-bound})
hold. In this case, by the Bernoulli law (\ref{Bernoulli}),
we have $c^2(\rho)=
c^2(\hat\rho)-(\g-1)D_t\dP-\f{\g-1}2|\nabla_x\dP|^2$,
which gives $CR(t)^{3(1-\g)}-C\ve R(t)^{3(1-\g)}<c^2(\rho)<CR(t)^{3(1-\g)}+C\ve R(t)^{3(1-\g)}$. Thus, one
obtains $c^2(\rho)\sim R(t)^{3(1-\g)}>0$ for any $t\ge 0$ and small $\ve>0$.
Therefore, the proof of Theorem
1.2 is completed by the local existence result in Theorem 3.1 and
continuation argument. $\hfill\square$

Finally,  Theorem 1.1 follows.

{\bf Proof of Theorem 1.1.} Under the assumptions of Theorem 1.1, it follows from Theorem 3.1 and
the smallness of $L$ that
(1.1)-(1.2) has a local solution $(\rho(t,x), u(t,x))$ that satisfies
$$\|\rho_0(1,x)-1\|_{H^4(S^1)}+\|u(1,x)\|_{H^4(S^1)}<\dl_0,
\eqno{(6.11)}$$
and
$$rot ~u(1,x)\equiv 0, $$
where $S^1=\{x: |x|=1+L\}$, $\dl_0>0$ is a small number depending only on $L$ and $\ve_0$. This, together with Theorem 1.2
for the time $t\ge 1$, yields Theorem 1.1.
$\hfill\square$

\medskip

{\bf Acknowledgement.} {\it The authors  wish to express their gratitude to Professor Yang Tong,
the City University of Hong Kong, for his
interests in this problem and many very fruitful discussions.
In particular,
Professor Yang Tong gave many suggestions and comments that led to substantial
improvement of the presentation.}

\end{document}